\documentstyle[
11pt,twoside,amsfonts,amssymb,amsmath,amsthm,latexsym]{article}
{\tiny }
\newcommand {\rel} {{\mathbb R}}
\newcommand {\com} {{\mathbb C}}
\newcommand {\CP} {{\mathbb CP}}
\newcommand {\nat} {{\mathbb N}}

\newcommand {\ganz} {{\mathbb Z}}
\newcommand {\quat} {{\mathbb H}}
\newcommand {\sphere} {{\mathbb S}}
\newcommand {\Will} {{\mathcal{W} }}
\newcommand {\Wil} {{\mathcal{E} }}

\newcommand {\Lift} {{\mathcal{L} }}

\newcommand {\PP} {{\mathcal{P} }}

\newcommand {\A} {{\mathcal{A} }}

\def    \Hn     {{ {\cal H}^2 }}
\def    \Hnm    {{ {\cal H}^{1} }}

\def    \Lno    {{ {\cal L}^1 }}

\begin{document}
	\thanks{Institute: Mathematics Department, Technion, Israel 
		Institute of Technology, 3200003 Haifa, Israel.\\
		Phone: 00972 549298486.
		E-mail: rubenj at technion.ac.il}
	
	\newtheorem{theorem}{Theorem}[section]
	\newtheorem{definition}{Definition}[section]
	\newtheorem{proposition}{Proposition}[section]
	\newtheorem{lemma}{Lemma}[section]
	\newtheorem{corollary}{Corollary}[section]
	\newtheorem{remark}{Remark}[section]
	\newtheorem{example}{Example}
	
	\author{Ruben Jakob}
	
	\title{The Willmore flow of Hopf-tori in the $3$-sphere}
	
	\maketitle
	
	\begin{abstract}
		In this article, the author investigates 
		flow lines of the classical Willmore flow, 
		which start to move in a smooth parametrization 
		of a Hopf-torus in $\sphere^3$. 
		We prove that any such flow line of the Willmore flow exists globally, in particular does not develop any singularities, and subconverges to some smooth Willmore-Hopf-torus in every $C^{m}$-norm. Moreover, 
		if in addition the Willmore energy 
		of the initial immersion $F_0$ is required to be smaller than or equal to the threshold $\frac{8\pi^2}{\sqrt{2}}$, 
		then the unique flow line of the Willmore flow, starting 
		to move in $F_0$, converges fully to a conformally transformed Clifford torus in every $C^{m}$-norm,
		up to time dependent, smooth reparametrizations.
		Key instruments for the proofs are the
		equivariance of the Hopf-fibration 
		$\pi:\sphere^3 \longrightarrow \sphere^2$ w.r.t. 
		the effect of the $L^2$-gradient of the Willmore energy applied to smooth Hopf-tori in $\sphere^3$ and to smooth closed regular curves in $\sphere^2$, a particular version of the Lojasiewicz-Simon gradient inequality, and a well-known classification and description of smooth, arc-length parametrized solutions of the Euler-Lagrange equation of the elastic energy functional in terms of Jacobi Elliptic Functions and Elliptic Integrals, dating back to the 80s. 
	\end{abstract}
	
MSC-codes: 53C42, 53E40, 35R01, 58J35, 11Z05

\section{Introduction}
\label{intro}
In this article, the author investigates
the long-term behaviour of the classical Willmore flow: 
\begin{equation}  \label{Willmore.flow}
	\partial_t f_t = -\frac{1}{2}  \,
	\Big{(} \triangle_{f_t}^{\perp} \vec H_{f_t} + Q(A^{0}_{f_t})
	(\vec H_{f_t}) \Big{)} \equiv - \nabla_{L^2} \Will(f_t),
\end{equation}
moving smooth families of smooth immersions
$f_t$ of a fixed compact torus into the standard 
$3$-sphere ``$\sphere^3$''. It is well-known, that 
the $L^2$-gradient flow of the Willmore functional
(\ref{Willmore.functional}) has a unique smooth 
short-time solution, starting to move in any prescribed smooth initial immersion $F_0:\Sigma \longrightarrow M$, for any Riemannian target manifold $M$ of dimension $\geq 3$. 
But there are only few results addressing the long-term behaviour and full convergence of this flow. Inspired by Kuwert's and Sch\"atzle's optimal convergence result on the classical Willmore flow, moving spherical immersions into $\rel^3$ -- see \cite{Kuwert.Schaetzle.Annals}, Theorem 5.2 -- and also by the recent article \cite{Dall.Acqua.Schaetzle.Mueller.2020}, 
treating the Willmore flow of tori of revolution in $\rel^3$, we are going to prove a global existence-, 
a subconvergence- and a full convergence-result in 
Theorem \ref{convergence} for flow lines of the Willmore flow 
in $\sphere^3$, starting to move in parametrizations of arbitrary smooth Hopf-tori in $\sphere^3$. More precisely, we are going to employ classical ideas and results due to Singer and Langer \cite{Langer.Singer}, \cite{Langer.Singer.1987} about the classification and quantitative analysis of elastica in 
simply connected space forms $M$ and some modern improvements 
and clarifications in \cite{Mandel.2018} and \cite{Mueller.Spener.2020} - treating the case in which $M$ 
is the hyperbolic plane - in order to accurately estimate 
the size of the gap between the Willmore energy of the Clifford torus and of any further Willmore-Hopf-torus in $\sphere^3$ 
within the framework of ``Elliptic Integrals and Elliptic Functions''. Such an estimate plays a key role in our ambitious pursuit 
of an optimal statement about ``full, smooth convergence'' 
of flow lines of the Willmore flow in $\sphere^3$ to the \underline{Clifford torus} - at least up to conformal equivalence; see our Theorem \ref{convergence}.  \\
The basic notion of this article is the ``Willmore energy
of a closed surface'':
\begin{equation} \label{Willmore.functional}
	\Will(f):= \int_{\Sigma} K^M_f + \frac{1}{4} \, 
	\mid \vec H_f \mid^2 \, d\mu_f,
\end{equation}
which is well-defined for $C^{2}$-immersions
$f:\Sigma \longrightarrow M$ mapping any closed smooth
orientable surface $\Sigma$ into an arbitrary smooth Riemannian
manifold $M$, where $K^M_f(x)$ denotes the
sectional curvature of $M$ w.r.t. the ``immersed tangent plane'' $Df_x(T_x\Sigma)$ in $T_{f(x)}M$. In those cases being relevant in this article, we will only have $K_f\equiv 0$ 
for $M=\rel^n$ or $K_f\equiv 1$ for $M=\sphere^n$, 
$n\geq 3$.\\
Now, given an immersion $f:\Sigma \longrightarrow \sphere^3$, we endow the torus $\Sigma$ with the pullback $f^*g_{\textnormal{euc}}$ of the Euclidean metric on $\sphere^3$, i.e. with coefficients 
$g_{ij}:=\langle \partial_i f, \partial_j f \rangle_{\rel^4}$, 
and we let $A_f$ denote the second fundamental form of the immersion $f$, defined on pairs of tangent vector fields $X,Y$ on $\Sigma$ by
\begin{eqnarray}  \label{second.fundam.form}
A_{f}(X,Y) \equiv A_{f,\sphere^3}(X,Y):= \nonumber  \\ 
=D_X(D_Y(f)) - P^{\textnormal{Tan}(f)}(D_X(D_Y(f)))
\equiv (D_X(D_Y(f)))^{\perp_{f}}
\end{eqnarray}
where $D_X(V)\lfloor_{x}$ denotes the projection of the derivative of a vector field $V:\Sigma \longrightarrow \rel^{4}$ in direction of the tangent vector field $X$ into the respective fiber $T_{f(x)}\sphere^3$ of $T\sphere^3$,
$P^{\textnormal{Tan}(f)}:\bigcup_{x \in \Sigma} \{x\} \times \textnormal{T}_{f(x)}\sphere^3 \longrightarrow
\bigcup_{x \in \Sigma} \{x\} \times \textnormal{T}_{f(x)}(f(\Sigma))
=:\textnormal{Tan}(f)$
denotes the bundle morphism which projects
the entire tangent space $\textnormal{T}_{f(x)}\sphere^3$ orthogonally onto its subspace $\textnormal{T}_{f(x)}(f(\Sigma))$, the tangent space of the immersion $f$ in $f(x)$, for every $x \in \Sigma$, and where $^{\perp_{f}}$ abbreviates the bundle morphism $\textnormal{Id}_{\textnormal{T}_{f(\cdot)}\sphere^n}-P^{\textnormal{Tan}(f)}$. Furthermore, $A^{0}_f$ denotes the tracefree part of $A_f$, i.e.
$$
A^{0}_f(X,Y) \equiv A^{0}_{f,\sphere^3}(X,Y) := 
A_{f,\sphere^3}(X,Y) - \frac{1}{2} \,g_f(X,Y)\, 
\vec H_{f,\sphere^3}
$$
and $\vec H_{f,\sphere^3} := \textnormal{Trace}(A_{f,\sphere^3}) 
\equiv g_f^{ij}\, A_{f,\sphere^3}(\partial_i,\partial_j)$ 
(``Einstein's summation convention'') 
denotes the mean curvature vector of $f$, which we shall always abbreviate by $\vec H_{f}$, if there cannot arise any confusion with the mean curvature vector $\vec H_{f,\rel^4}$ of $f$ interpreted as an immersion into $\rel^4$. Finally, $Q(A^0_f)$ operates 
on sections $\phi$ into the normal bundle of $f$, 
by assigning $Q(A^0_f)(\phi):=A^0_f(e_i,e_j) \langle A^0_f(e_i,e_j),\phi \rangle$, which is by definition again a section of the normal bundle of $f$ within $T\sphere^3$.
We know from \cite{Weiner}, Section 2, respectively from 
Proposition 3.1 in \cite{Ndiaye.Schaetzle.2014}, 
that the first variation of the Willmore functional 
$\delta \Will(f,\phi)$
in a smooth immersion $f:\Sigma \longrightarrow \sphere^3$, 
in direction of a smooth section $\phi$
of the normal bundle of $f$ within $T\sphere^3$, is given by:
\begin{eqnarray} \label{first_variation}
	\delta \Will(f,\phi)
	= \frac{1}{2} \int_{\Sigma}
	\langle \triangle^{\perp}_{f} \vec H_{f} + Q(A^{0}_{f})(\vec H_{f}), \phi \rangle\, d\mu_{f}
	=: \int_{\Sigma} \langle \nabla_{L^2} \Will(f),\phi \rangle \, d\mu_{f}, \quad
\end{eqnarray}
which justifies the second equality sign in (\ref{Willmore.flow}).\\
In contrast to Kuwert's and Sch\"atzle's convergence theorem in \cite{Kuwert.Schaetzle.Annals} and in contrast to the new results in \cite{Dall.Acqua.Schaetzle.Mueller.2020}, treating the Willmore flow of tori of revolution in $\rel^3$, we do not necessarily have to start the Willmore flow below the prominent $8\pi$-energy threshold, because our approach does not rely on the Li-Yau inequality \cite{Li.Yau}, relating the pointwise $2$-dimensional Hausdorff-density of the image of an immersion into some $\rel^n$ to its Willmore energy.
Instead, we will exploit a certain type of ``dimension reduction'' - worked out in Section \ref{The.link} below - by means of a careful analysis of the Hopf-fibration, we will prove in Section \ref{Proof.main.theorem} the absence of singularities along the considered flow lines of the Willmore flow in $\sphere^3$, 
and we will compute the critical levels of the elastic energy (\ref{elastic.energy.functional}) as precisely as possible in our appendix, Section \ref{Appendix}, in order to obtain in Section \ref{Proof.main.theorem} the surprisingly strong 
``energy threshold'' $\frac{8\pi^2}{\sqrt{2}}$ for full, 
smooth convergence of the considered flow lines to the \underline{Clifford torus} in $\sphere^3$, up to M\"obius-transformations of $\sphere^3$, which is the 
final result of Theorem \ref{convergence} below.   
\begin{remark}
We should point out here, at the end of the introduction, 
that the entire investigation of this article was motivated 
by the examination of another variant of the Willmore flow, 
namely of the evolution equation 
	\begin{equation} \label{Moebius.flow}
	\partial_t f_t = -\frac{1}{2} \mid A^0_{f_t} \mid^{-4} \, \Big{(} \triangle^{\perp}_{f_t} \vec H_{f_t} + Q(A^{0}_{f_t})(\vec H_{f_t}) \Big{)} 
	\equiv - \mid A^0_{f_t} \mid^{-4} \,\nabla_{L^2} \Will(f_t), \quad
	\end{equation} 
	for differentiable families of $C^4$-immersions 
	$f_t:\Sigma \longrightarrow M$, with $M=\rel^n$ or  
	$M=\sphere^n$. As already pointed out in the author's article \cite{Jakob_Moebius_2016}, the ``umbilic free condition'' $|A^0_{f_t}|^2>0$ on $\Sigma$, implies
	$\chi(\Sigma)=0$ for the Euler-characteristic of the 
	surface $\Sigma$, which forces the flow (\ref{Moebius.flow}) 
	to be only well-defined on families of sufficiently 
	smooth umbilic-free tori, immersed into $\rel^n$ or $\sphere^n$.
	In \cite{Jakob_Moebius_2016}, the author has proved 
	short-time existence and uniqueness of flow lines 
	of this flow, starting to move in a umbilic-free 
	immersion of some fixed smooth compact torus $\Sigma$ 
	into $\rel^n$ or $\sphere^n$, for $n\geq 3$. 
	The big advantage of this flow - compared to the classical Willmore flow - is its ``conformal invariance'', which means that any family $\{f_t\}_{t\in [0,T]}$ of $C^4$-immersions $f_t:\Sigma \longrightarrow \sphere^n$ without any umbilic points, i.e. with $\mid A^0_{f_t} \mid^2>0$ on $\Sigma$\, $\forall \, t\in[0,T]$, solves flow equation (\ref{Moebius.flow}) on $\Sigma \times [0,T]$, if and only if its composition $\Phi(f_t)$ with an arbitrary applicable M\"obius-transformation $\Phi$ of $\sphere^n$ solves the same flow equation on $\Sigma \times [0,T]$ again. This property allows us to project flow lines of this flow with original target manifold $\rel^3$ into $\sphere^3$, which is a compact Lie-group, diffeomorphic to $\textnormal{SU}(2)$, the unit sphere within the division algebra $\quat$ of quaternions and which can be fibered by great circles by means of the Hopf-fibration $\pi:\sphere^3 \longrightarrow \sphere^2$. These rather elementary algebraic and topological properties of the $3$-sphere finally turned out to be very useful, in order to investigate the 
	\underline{classical Willmore flow in the $3$-sphere}, 
	and the method of this article can actually 
	be interpreted as a straight-forward continuation of Pinkall's paper \cite{Pinkall} about concrete conformal embeddings of flat tori into $\sphere^3$ and about  
	Hopf-tori which are ``Willmore'' but not ``minimal'' in $\sphere^3$. Unfortunately, flow lines of 
	the \underline{conformally invariant Willmore flow} (\ref{Moebius.flow}) might in general develop singularities, even if they start moving in smooth parametrizations of Hopf-tori with Willmore energies below $8\pi$. The author has already started to address this apparently difficult problem in his preprint \cite{Ruben.MIWF.V}, which relies on the preparatory results of Sections \ref{main.results} and \ref{The.link} in this article.    
\end{remark}
\noindent

\section{Main results and main tools}  \label{main.results}
The basic idea of this article is the use of the Hopf-fibration
$$
\sphere^1 \hookrightarrow \sphere^3
\stackrel{\pi}{\longrightarrow} \sphere^2
$$
and its equivariance w.r.t. the first variation 
of the Willmore energy along Hopf-tori in $\sphere^3$ 
and closed curves in $\sphere^2$; 
see formula (\ref{DHopf.first.Variation}) below.
In order to work with the most effective formulation of the
Hopf-fibration, we shall consider $\sphere^3$ as the 
subset of the four-dimensional $\rel$-vector space $\quat$ of quaternions, whose elements have length $1$, i.e. $\sphere^3:=\{q \in \quat \,|\, \bar q\cdot q=1 \}$. 
We shall use the usual notation 
for the generators of the division algebra $\quat$, 
i.e. 1,$i$,$j$,$k$. We therefore decompose every quaternion 
in the way $q=q_1+i\, q_2 +j\, q_3 + k\,q_4$, 
for unique ``coordinates'' $q_1, q_2, q_3, q_4 \in \rel$,
and especially the conjugate of any quaternion $q$ can be 
written as $\bar q= q_1 - i\, q_2 - j\, q_3 - k\,q_4$.
Moreover, we identify
$$
\sphere^2 =\{q \in \textnormal{span}\{1,j,k\}\, 
| \,\bar q\cdot q=1 \}
=\sphere^3 \cap \textnormal{span}\{1,j,k\},
$$
and we shall use the particular involution $q \mapsto \tilde q$ of $\quat$, which fixes the generators $1$, $j$ and $k$, but sends $i$ to $-i$. Following Section 2 of \cite{Pinkall}, we employ this involution to write the Hopf-fibration in the elegant way
\begin{equation} \label{Hopf}
\pi :\quat \longrightarrow  \quat, \qquad q \mapsto \tilde q \cdot q.
\end{equation}
We gather its most important
properties in the following elementary lemma.
\begin{lemma}  \label{Hopf.properties}
\begin{itemize}
\item[1)] $\pi(\sphere^3) = \sphere^3 
\cap \textnormal{Span}\{1,j,k\}=\sphere^2$.
\item[2)] $\pi(e^{i\varphi} q )= \pi(q)$, $\forall \,
\varphi \in \rel$ and $\forall \,q\in \sphere^3$.
\item[3)] The group $\sphere^3$ acts isometrically on $\sphere^3$ by right (and left) multiplication, and it acts isometrically on $\sphere^2$, which means that every $r \in \sphere^3$ induces the rotation
\begin{equation} \label{Spin1}
q \mapsto \tilde r \cdot q \cdot r, \quad \textnormal{for} \,\,\,q\in \sphere^2. \,\,\,
\end{equation}
Moreover, there holds
\begin{equation}  \label{equivariant}
\pi(q \cdot r) = \tilde r\cdot \pi(q) \cdot r  
\qquad \forall  \, q,r\in \sphere^3,
\end{equation}
which means that right multiplication
on $\sphere^3$ translates equivariantly via the
Hopf-fibration to rotation in $\sphere^2$.
\item[4)] The differential of $\pi$ in any $q\in \quat$, applied to some $v \in \quat$, reads
\begin{equation}   \label{Hopf.differential}
D\pi_q(v) = \,\tilde v \cdot q + \tilde q \cdot v.
\end{equation} 
\end{itemize}
\end{lemma}
\noindent
\proof:
Assertions (1)-(3) follow immediately from the definition 
in (\ref{Hopf}) and can also be found in Section 2 of \cite{Pinkall}. In order to prove formula (\ref{Hopf.differential}), we choose a smooth curve 
$[t \mapsto \gamma(t)]$, satisfying $\gamma(0)=q$ and $\gamma'(0)=v$ and compute by means of the linearity 
of the involution $I: q \mapsto \tilde q$:
\begin{eqnarray*}
	D\pi_q(v)= \partial_t \pi(\gamma(t))|_{t=0}
	= \partial_t (\tilde \gamma(t) \,\gamma(t))|_{t=0}\\
	=  \, DI(\gamma(0)).\gamma'(0) \cdot \gamma(0)
	+ I(\gamma(0)) \cdot \gamma'(0)
	= I(v)\cdot q + \tilde q \cdot v.
\end{eqnarray*}
\qed
\noindent
By means of the Hopf-fibration we introduce Hopf-tori in the following definition. See also Section 2 of \cite{Pinkall} 
for further explanations:
\begin{definition} \label{Hopf.torus.immersion}
	\begin{itemize}
		\item[1)] Let $\gamma: [a,b] \longrightarrow \sphere^2$ be a regular, smooth and closed curve in $\sphere^2$, and let $\eta: [a,b] \to \sphere^3$ be a smooth lift of $\gamma$ w.r.t. $\pi$ into $\sphere^3$, i.e. a smooth map from $[a,b]$ into $\sphere^3$
		satisfying $\pi \circ \eta = \gamma$. We define
		\begin{equation}   \label{Hopf.Torus}
		X(s,\varphi):= e^{i\varphi}  \cdot \eta(s), \quad
		\forall\, (s,\varphi) \in [a,b] \times [0,2\pi],
		\end{equation}
		and note that $(\pi \circ X)(s,\varphi) = \gamma(s)$,
		$\forall \,(s,\varphi)\in [a,b] \times [0,2\pi]$.
		\item[2)] We call this map $X$ the 
		``Hopf-torus-immersion''
		and its image respectively $\pi^{-1}(\textnormal{trace}(\gamma))$
		the ``smooth Hopf-torus'' in $\sphere^3$ w.r.t. the smooth curve $\gamma$.
	\end{itemize}
\end{definition}
\noindent
Preparing ourselves for precise computations which involve Hopf-tori and lifts w.r.t. the Hopf-fibration $\pi$, we need the following lemma, already anticipating the elementary 
result of Lemma \ref{closed.lifts} in the appendix.
\begin{lemma} \label{Hopf.Tori}
	Let $\gamma: \rel/L\ganz \longrightarrow \sphere^2$ regularly parametrize a smooth closed curve in $\sphere^2$. Moreover, let $\eta:[0,L] \to \sphere^3$ be a smooth lift of $\gamma$ w.r.t. $\pi$, having constant speed $|\eta'| \equiv 1$ and intersecting the fibers of $\pi$ perpendicularly, see here Lemma \ref{closed.lifts} below. Then there holds:
	\begin{itemize}
	\item[1)] $\eta' = u\,\eta$ for some function
	$u: [0,L] \longrightarrow \textnormal{span}\{j,k\}$ satisfying $|u(s)| \equiv 1$.
	\item[2)] $\partial_{\varphi}X(s,\varphi) = i e^{i\varphi} \eta(s)$ and thus $\Re(\overline{\eta'}(s) 
	\,i e^{i\varphi}\, \eta(s)) \equiv 0$
	and $|\partial_{\varphi}X(s,\varphi)| \equiv 1$.
	\item[3)] $\partial_{s} X(s,\varphi) = e^{i\varphi}\, u(s) \, \eta(s)$, and thus also 
	$|\partial_{s} X(s,\varphi)|\equiv 1$ and 
	\begin{equation}   \label{g.12} 
	\Re\big{(}\overline{\partial_{s} X(s,\varphi)}\cdot \partial_{\varphi}X(s,\varphi)\big{)} \equiv 0, 
	\quad \forall \,(s,\varphi)\in [0,L] \times [0,2\pi]. 
	\end{equation}
	\item[4)] $\gamma'(s) = 2 \, \tilde \eta(s) \,u(s) \,\eta(s)$, $\forall\, s\in [0,L]$.
	Every horizontal smooth lift $\eta$ of $\gamma$ w.r.t. $\pi$ intersects each fiber of $\pi$ exactly 
	$m\geq 1$ times on $[0,L]$ with constant speed $1$, 
	if and only if $\gamma$ performs $m$ loops through 
	its trace on $[0,L]$ with constant speed $2$. In this case, the length of $\textnormal{trace}(\gamma)$ is $2L/m$.
	\end{itemize}
\end{lemma}
\noindent
\proof:
1) From $|\eta|^2 \equiv 1$ and $|\eta'|^2 \equiv 1$ on $[0,L]$ we infer that $\eta'(s)$ is a unit vector, 
which is perpendicular to $\eta(s)$, for every $s \in [0,L]$.
Therefore, we can obtain the vector $\eta'(s)$ by means of a rotation of the vector $\eta(s)$ about a right angle within $\sphere^3$. In $\quat$ this can be achieved by means of left multiplication with a unit-length
quaternion $u(s)$, being contained in $\textnormal{Span}\{i,j,k\}$. Hence, since we assume that 
$\eta$ is a horizontal lift of $\gamma$ w.r.t. $\pi$, we 
infer from the second statement of Lemma \ref{Hopf.properties} 
that $u(s) \in \textnormal{Span}\{j,k\}$, for every 
$s \in [0,L]$.\\
2) The assertions of the second part of the 
lemma follow immediately from formula (\ref{Hopf.Torus}) and  
from the assumption on the lift $\eta$ to be horizontal w.r.t. $\pi$.  \\
3) We compute:
$$\overline{\partial_{s} X(s,\varphi)}\cdot
\partial_{\varphi}X(s,\varphi)
= \bar \eta(s) \bar u(s)  
\,e^{-i\varphi} \,i e^{i\varphi} \eta(s) 
= \bar \eta(s) \bar u(s)  \,i\,\eta(s)
= \overline{\eta'}(s) i \eta(s).
$$
Using the formula\, $\Re(\overline{\eta'}(s) \,i e^{i\varphi}\, \eta(s)) \equiv 0$\, from the second part of this lemma, we obtain especially $\Re(\overline{\eta'}(s) i \eta(s)) \equiv 0$, which proves formula (\ref{g.12}).\\
4) Using formula (\ref{Hopf.differential}) and
the chain rule, we compute:
\begin{eqnarray} \label{gamma.speed}
	\gamma'(s)=(\pi \circ \eta)'(s)
	=D\pi_{\eta(s)}(\eta'(s))
	=(I(\eta'(s)) \, \eta(s))
	+(I(\eta(s)) \, \eta'(s)) \\
	= (\tilde \eta(s) \,\tilde u(s)\, \eta(s))
	+(\tilde \eta(s) \, u(s) \, \eta(s))
	=  2 \tilde \eta(s) \, u(s) \eta(s), \nonumber
\end{eqnarray}
because $u(s) \in \textnormal{Span}\{j,k\}$ 
$\forall \,s\in \rel/L\ganz$, and therefore
$\tilde u \equiv u$. In particular, there holds
$|\gamma'| \equiv 2$ on $\rel/L\ganz$.  
In fact, this statement is equivalent to 
``$|\eta'|\equiv 1$ on $\rel/L\ganz$'' on account 
of formula (\ref{gamma.speed}) combined with the uniqueness of horizontal smooth lifts w.r.t. $\pi$, as stated and proved in Lemma \ref{closed.lifts} below. 
\qed
\noindent
\begin{remark} \label{Remark.lift}
	We can conclude from point (4) of Lemma \ref{Hopf.Tori}
	that every horizontal smooth lift $\eta$ of a closed
	regular smooth curve $\gamma:[0,L] \longrightarrow \sphere^2$ w.r.t. $\pi$ covers the path $\gamma$ with exactly half the speed of $\gamma$. This fact will not cause any technical problems in the sequel.
	It only leads to a multiplication of the ``elastic energy''
	$\frac{1}{2} \int_{\sphere^1} 1+|\kappa_{\gamma}|^2 \,d\mu_{\gamma}$ by the factor $2$.
	The algebraic background of this phenomenon is the
	short exact sequence of groups 
	$$
	1 \rightarrow \ganz/2 \longrightarrow
	\textnormal{Spin}(3) \stackrel{\mathcal{U}}{\longrightarrow}
	\textnormal{SO}(3) \rightarrow 1
	$$
	where $\mathcal{U}$ is defined by
	\begin{equation} \label{Spin2}
	\textnormal{Spin}(3) \ni [q \mapsto q\cdot r]
	\mapsto [q\mapsto \tilde r \, q \, r\,] \in
	\textnormal{SO}(3), \quad \textnormal{for every} \,\,\, 
	r \in \sphere^3,
	\end{equation}
	and satisfies $\pi(q \cdot r)=\mathcal{U}(r).(\pi(q))$,  for every $q,r \in \sphere^3$,
	on account of formula (\ref{equivariant}).
	The above exact sequence only rephrases the 
	topological fact, that $\mathcal{U}$ is the 
	universal covering of $\textnormal{SO}(3)$. 
\qed
\end{remark}
\noindent
In order to estimate, ``how large'' the subset of all 
Hopf-tori within the set of all smoothly immersed 
tori in $\sphere^3$ actually is, we shall follow 
Section 2 of \cite{Pinkall} and introduce ``abstract Hopf-tori''.
\begin{definition} \label{Hopf.curve}
	Let $\gamma: [0,L/2] \to \sphere^2$ be a path
	with constant speed $2$ which traverses a
	simple closed smooth curve in $\sphere^2$ of length
	$L>0$ and encloses the area $A$
	of the domain on $\sphere^2$, ``lying on the left hand side'' when performing one loop through the trace of $\gamma$. We assign to $\gamma$ the lattice $\Gamma_{\gamma}$ which is generated by the vectors $(2\pi,0)$ and $(A/2,L/2)$, and we call the quotient $M_{\gamma}:=\com/\Gamma_{\gamma}$
	the ``abstract Hopf-torus'' corresponding to $\textnormal{trace}(\gamma)$.
\end{definition}
\begin{remark} \label{Surfaces.genus.1}
	Pinkall proved in Proposition 1 of \cite{Pinkall}, that for every simple closed smooth curve $\gamma$ in $\sphere^2$ its preimage $\pi^{-1}(\textnormal{trace}(\gamma))$ can 
	be isometrically mapped onto its corresponding
	``abstract Hopf-torus'' $M_{\gamma}:=\com/\Gamma_{\gamma}$,
	as defined in Definition \ref{Hopf.curve}.
	Taking the uniformization theorem into account, Pinkall derived from this result in Section 3 of \cite{Pinkall}, that every conformal class of a compact Riemann surface of genus one can be realized by a Hopf-torus, i.e. by the preimage $\pi^{-1}(\textnormal{trace}(\gamma))$ of some appropriate smooth, simple and closed curve $\gamma$, as introduced in Definition \ref{Hopf.torus.immersion}. 
	Therefore the subset of Hopf-tori is sufficiently large, such that the main results of this article can be considered, to be ``significant''.
\end{remark}
\noindent
In order to rule out inappropriate parametrizations
of Hopf-tori along the flow lines of the Willmore flow, we 
shall introduce ``topologically simple" maps between tori in the following definition.
\noindent
\begin{definition} \label{General.lift}
	Let $\Sigma_1$ and $\Sigma_2$ be two compact tori, i.e. two
	compact sets which are homeomorphic to the standard torus
	$\sphere^1 \times \sphere^1$.
	We term a continuous map $F:\Sigma_1 \longrightarrow \Sigma_2$ \underline{simple}, if it has mapping degree $\pm 1$, i.e. if the induced map
	$$
	(F_{*2}): H_2(\Sigma_1,\ganz) \stackrel{\cong}\longrightarrow
	H_{2}(\Sigma_2,\ganz)
	$$
	is an isomorphism between these two singular homology groups in degree $2$.	
\end{definition}
\begin{remark} \label{Fundamental.gruppe.Homologie}
\begin{itemize}
		\item[1)] We should point out here - especially regarding
		the proof of the third part of Theorem \ref{convergence} - that if $\Sigma_1$ is a smooth compact manifold
		of genus $1$, i.e. a ``smooth compact torus'' without self-intersections, and if $F:\Sigma_1 \longrightarrow \sphere^n$ is a smooth immersion of $\Sigma_1$ into $\sphere^n$, $n\geq 3$, such that $F$ maps $\Sigma_1$ simply onto its image $\Sigma_2:=F(\Sigma_1)$, then $F$ is a smooth diffeomorphism between $\Sigma_1$ and $\Sigma_2$, if and only if $\Sigma_2$ is a smooth compact manifold of genus $1$, as well. This follows immediately from our Definition \ref{General.lift},
		Propositions 4.5 and 4.7 in Chapter VIII of \cite{Dold} and the remark following Proposition 4.7 in Chapter VIII of \cite{Dold}.
		\item[2)] If $\gamma:\rel/L\ganz \longrightarrow \sphere^2$ is a smooth, regular closed path, which traverses its trace exactly once with speed $2$, and if $\eta:[0,L] \longrightarrow \pi^{-1}(\textnormal{trace}(\gamma))$ is a horizontal smooth lift of $\gamma$ w.r.t. $\pi$, as discussed in Lemma \ref{closed.lifts} below, then the standard parametrization $X:[0,L] \times [0,2\pi] \longrightarrow \pi^{-1}(\textnormal{trace}(\gamma))$
		of the Hopf-torus $\pi^{-1}(\textnormal{trace}(\gamma))$ in $\sphere^3$, given by 
		$X(s,\varphi):=e^{i\varphi} \, \eta(s)$, has the 
		effect of a typical simple parametrization, namely 
		to cover $\pi^{-1}(\textnormal{trace}(\gamma))$ 
		``essentially once''. Precisely this means, that $X^{-1}(\{z\})$ consists of exactly one element 
		for $\Hn$-almost every 
		$z \in \pi^{-1}(\textnormal{trace}(\gamma))$.
\end{itemize}
\end{remark}
\noindent
Now we can state the main results of this article: 	
\begin{theorem} \label{convergence} 
[Global existence and Subconvergence]
Let $\gamma_0: \sphere^1 \to \sphere^2$ be a 
smooth, closed and regular path in $\sphere^2$.
Let $F_0:\Sigma \longrightarrow \sphere^3$ be an arbitrary smooth immersion, which maps a compact smooth torus $\Sigma$ simply onto the Hopf-torus $\pi^{-1}(\textnormal{trace}(\gamma_0))$ in the sense of Definition \ref{General.lift}. 
Then the following statements hold:
\begin{itemize}
\item[I)] There is a unique smooth global solution $\{\PP(t,0,F_0)\}_{t \geq 0}$
of the Willmore flow (\ref{Moebius.flow}) on $[0,\infty) \times \Sigma$, starting in $F_0$ at time $t=0$. The immersions
$\PP(t,0,F_0)$ remain umbilic-free and map $\Sigma$ simply onto Hopf-tori \,$\forall \,t \geq 0$.
\item[II)] This global solution $\{\PP(t,0,F_0)\}$ of equation (\ref{Willmore.flow}) subconverges in every $C^{m}$-norm - up to reparametrizations - to smooth Willmore-Hopf-tori in $\sphere^3$. More precisely we have:
For every sequence $t_j \to \infty$, there are a subsequence of times $t_{j_k} \to \infty$, a sequence of diffeomorphisms
$\varphi_{k}:\Sigma \stackrel{\cong}\longrightarrow \Sigma$ and a smooth immersion $\hat F:\Sigma \longrightarrow \sphere^3$ mapping $\Sigma$ simply onto some smoothly immersed Willmore-Hopf-torus in $\sphere^3$, such that there holds
\begin{eqnarray}  \label{convergence.1.3}
\PP(t_{j_k},0,F_0) \circ \varphi_k \longrightarrow \hat F \qquad \textnormal{in} \,\,\,C^{m}(\Sigma,\rel^4),
\end{eqnarray}
as $k \to \infty$, $\forall \, m \in \nat_0$.
\item[III)]
If the initial immersion $F_0:\Sigma \longrightarrow \sphere^3$ maps the smooth torus $\Sigma$ simply onto a Hopf-torus in $\sphere^3$ with Willmore energy 
$\Will(F_0) \leq \frac{8\pi^2}{\sqrt{2}}$, then 
there is a smooth family of smooth diffeomorphisms
$\Psi_t:\Sigma \longrightarrow \Sigma$, such that the reparametrization $\{\PP(t,0,F_0) \circ \Psi_t\}$ of the flow line $\{\PP(t,0,F_0)\}$, starting in $F_0$ 
at time $t=0$, consists of smooth diffeomorphisms between 
the smooth torus $\Sigma$ and their images in $\sphere^3$ for sufficiently large $t$, and the family $\{\PP(t,0,F_0) \circ \Psi_t\}$ converges fully in $C^m(\Sigma,\rel^4)$, for each $m\in \nat_0$, to some smooth diffeomorphism $F^*$ between $\Sigma$ and a conformally transformed Clifford torus 
in $\sphere^3$.
\qed
\end{itemize}
\end{theorem}
\noindent
See here also Appendix B in \cite{Dall.Acqua.Schaetzle.Mueller.2020}, especially 
Definitions B.7 and B.9 for the precise meaning of 
the terminology ``convergence in $C^{m}(M,\rel^4)$''. 
For the proof of the $C^{m}$-subconvergence (\ref{convergence.1.3}) in Part II of Theorem \ref{convergence} we will need the following special version of a general compactness-theorem due to Breuning \cite{Breuning} for proper smooth immersions of closed manifolds into $\rel^n$, which is a generalization of Langer's compactness theorem \cite{Langer}, and which appeared already in \cite{Kuwert_Schaetzle_2001}, Theorem 4.2, without proof.
\begin{proposition}  
[\cite{Breuning}, Theorem 1.3] \label{Breuning}
Let $F_j : M \longrightarrow \rel^4$ be a sequence of smooth immersions, where $M$ is a two-dimensional compact manifold without boundary. If there is a ball 
$B_R(0) \subset \rel^4$ with $F_j(M) \subset B_R(0)$, for every $j \in \nat$, and if there are constants $C_1>0$ and $C_2(m)>0$ for each $m \in \nat_0$, such that
\begin{eqnarray*}
\A(F_j):=\int_{M} d\mu_{F_j} \leq C_1 \qquad \textnormal{and} \qquad
	\parallel (\nabla^{\perp_{F_j}})^m(A_{F_j}) \parallel_{L^{\infty}(M)} \leq C_2(m)
\end{eqnarray*}
holds for every $j\in \nat$ and each $m \in \nat_0$, then there exists a smooth immersion 
$\hat F: M \longrightarrow \rel^4$, 
some subsequence $\{F_{j_k}\}$, and some sequence 
of smooth diffeomorphisms $\varphi_k: 
M \stackrel{\cong}\longrightarrow M$, such that
$$
F_{j_k} \circ \varphi_{k} \longrightarrow
\hat F \qquad \textnormal{converge in} \,\,\, C^{m}(M,\rel^4),
$$
as $k \to \infty$, for every $m \in \nat_0$, 
in the sense of Definition B.9 in Appendix B of \cite{Dall.Acqua.Schaetzle.Mueller.2020}.
\qed 
\end{proposition}
\noindent 
Moreover, for the proof of the full convergence of the auxiliary flow (\ref{elastic.energy.flow}) to smooth parametrizations of great circles in $\sphere^2$, we will need the ``Lojasiewicz-Simon gradient inequality'' for the elastic energy functional $\Wil$, restricted to the non-linear subspace $C^4_{\textnormal{reg}}(\sphere^1,\sphere^2)$ of the Banach space $C^4(\sphere^1,\rel^3)$:   
\begin{proposition}[Lojasiewicz-Simon gradient inequality]  \label{Lojasiewicz}
Let $\gamma^*:\sphere^1 \longrightarrow \sphere^2$ be a closed curve of class $C^{4}_{\textnormal{reg}}(\sphere^1,\sphere^2)$, which is also a critical point of the elastic energy functional 
	\begin{equation}  \label{elastic.energy.functional} 
	\Wil(\gamma):= \int_{\sphere^1} 1 + |\vec{\kappa}_{\gamma}|^2 \, d\mu_{\gamma}.
	\end{equation} 
	There exist constants $\theta \in (0,\frac{1}{2}]$,
	$c \geq 0$ and $\sigma>0$, only depending on $\gamma^*$, 
	such that for every closed curve
	$\gamma\in C^{4}_{\textnormal{reg}}(\sphere^1,\sphere^2)$ 
	satisfying $\parallel \gamma - \gamma^* \parallel_{C^{4}(\sphere^1,\rel^3)}
	\leq \sigma$ there holds:
	$$
	|\Wil(\gamma) - \Wil(\gamma^*) |^{1-\theta} \leq c \,\,\Big{(} \int_{\sphere^1} |\nabla_{L^2} \Wil(\gamma)|^2 \, d\mu_{\gamma} \Big{)}^{1/2}.
	$$
\end{proposition}
\proof 
First of all, we recall the Lojasiewicz-Simon-gradient-inequality 
for the Willmore functional $\Will$ applied to immersions $F:\Sigma \longrightarrow \rel^4$ of a torus $\Sigma$ into $\rel^4$, which is proved in Theorem 3.1 in \cite{Chill_Schatz}:
For any critical immersion $F^*:\Sigma \longrightarrow \rel^4$ of the Willmore functional 
$\tilde \Will(F):=\frac{1}{4}\, 
\int_{\Sigma} |\vec H_{F,\rel^4}|^2 \, d\mu_F$ there exist constants $\theta \in (0,\frac{1}{2}]$, $c \geq 0$ 
and $\sigma>0$, such that for every immersion $F\in C^{4}(\Sigma,\rel^4)$ satisfying 
$\parallel F - F^* \parallel_{C^{4}(\Sigma,\rel^4)}
\leq \sigma$ there holds:
\begin{equation}  \label{Lojasiewicz.Simon.R4}
|\tilde \Will(F) - 
\tilde \Will(F^*)|^{1-\theta} 
\leq c \,\Big{(}\int_{\Sigma} |\nabla_{L^2} \tilde \Will(F)|^2 \, d\mu_{F}\Big{)}^{1/2}.
\end{equation} 
Now we recall from \cite{Ndiaye.Schaetzle.2014}, formula 
(2.6), that the $L^2$-gradient of the functional 
$\tilde \Will$ coincides with the $L^2$-gradient of the functional $\Will$ in $C^4$-immersions 
$F:\Sigma \longrightarrow \sphere^3 \subset \rel^4$, i.e. such immersions satisfy: 
\begin{equation}  \label{first.variation.coincides} 
	\triangle^{\perp,\rel^4}_F(\vec H_{F,\rel^4}) 
	+ Q(A^0_{F,\rel^4})(\vec H_{F,\rel^4}) 
	= \triangle^{\perp,\sphere^3}_F
	(\vec H_{F,\sphere^3})  
	+ Q(A^0_{F,\sphere^3})(\vec H_{F,\sphere^3}), 
\end{equation} 
although their second fundamental tensors and 
mean curvature vectors do not coincide, namely there holds: 
\begin{equation}  \label{compare.second.fund.form}
A_{F,\rel^4} = A_{F,\sphere^3} - F\, g_F, \qquad
\vec H_{F,\rel^4} = \vec H_{F,\sphere^3} - 2 F,
\end{equation}
for $C^2$-immersions $F:\Sigma \longrightarrow \sphere^3 \subset \rel^4$ by formulae (2.2) and (2.3) in \cite{Ndiaye.Schaetzle.2014}. Hence, we have due to 
$\langle \vec H_{F,\sphere^3},F \rangle_{\rel^4} \equiv 0$ and $|F|^2 \equiv 1$:    
$|\vec H_{F,\rel^4}|^2 = 
|\vec H_{F,\sphere^3}|^2 + 4$   
for $C^2$-immersions $F:\Sigma \longrightarrow \sphere^3 \subset \rel^4$, and we consequently obtain together with formulae (\ref{Lojasiewicz.Simon.R4}) and (\ref{first.variation.coincides}): 
\begin{eqnarray} 
	|\Will(F) - \Will(F^*)|^{1-\theta}
	= |\tilde \Will(F) - \tilde \Will(F^*)|^{1-\theta}           \nonumber \\ 
	\leq c \,\Big{(}\int_{\Sigma} |\nabla_{L^2} \tilde \Will(F)|^2 \, d\mu_{F}\Big{)}^{1/2}      
	= c \,\Big{(}\int_{\Sigma} |\nabla_{L^2}\Will(F)|^2 \, d\mu_{F}\Big{)}^{1/2}  \label{Lojasiewicz.revisited}   
\end{eqnarray} 
for two $C^4$-immersions $F,F^*:\Sigma \longrightarrow \sphere^3$ satisfying $\parallel F - F^* \parallel_{C^{4}(\Sigma,\rel^4)}\leq \sigma$, 
provided $F^*$ is Willmore in $\sphere^3$, i.e. 
satisfies the equation 
$$
\triangle^{\perp,\sphere^3}_{F*}
(\vec H_{F^*,\sphere^3})  
+ Q(A^0_{F^*,\sphere^3})(\vec H_{F^*,\sphere^3}) 
\equiv 0 \quad \textnormal{on} \,\,\, \Sigma, 
$$ 
and provided the number $\sigma>0$ is sufficiently small. 
Now we infer from formula (\ref{DHopf.first.Variation}),
that the Willmore-Hopf-tori in $\sphere^3$  
correspond exactly to the elastic curves in $\sphere^2$, i.e. to the critical points of $\Wil$ in 
$C^4_{\textnormal{reg}}(\sphere^1,\sphere^2)$.  
Moreover we note, that the Hopf-fibration 
$\pi:\sphere^3 \longrightarrow \sphere^2$ 
is a smooth submersion, and it can therefore be 
locally transformed into an orthogonal projection from 
$\rel^3$ onto $\rel^2$. 
Hence, fixing some elastic curve 
$\gamma^*:\sphere^1 \to \sphere^2$,   
we infer from estimate (\ref{Lojasiewicz.revisited}), 
combined with formulae (\ref{Will.Wil}) and 
(\ref{First.Var.Will.Wil}) below, the existence of some  $\sigma>0$, such that for any closed curve $\gamma \in C^4_{\textnormal{reg}}(\sphere^1,\sphere^2)$ with 
$\parallel \gamma - \gamma^* \parallel_{C^{4}(\sphere^1,\rel^3)}
\leq \sigma$ there holds:
\begin{eqnarray*}  
\pi^{1-\theta}\, |\Wil(\gamma) - \Wil(\gamma^*)|^{1-\theta} =
|\Will(F_{\gamma}) - \Will(F_{\gamma^*})|^{1-\theta}      \\    
\leq c \,\Big{(}\int_{\Sigma} |\nabla_{L^2}\Will(F_{\gamma})|^2 \, d\mu_{F_{\gamma}}\Big{)}^{1/2}             
= c \, 2 \sqrt \pi\,
\Big{(}\int_{\sphere^1} |\nabla_{L^2}\Wil(\gamma)|^2 \, d\mu_{\gamma}\Big{)}^{1/2},     
\end{eqnarray*}
where $F_{\gamma}$ and $F_{\gamma^*}$ are appropriately chosen, simple $C^4$-parametrizations of the Hopf-tori 
$\pi^{-1}(\textnormal{trace}(\gamma))$ and 
$\pi^{-1}(\textnormal{trace}(\gamma^*))$ 
in $\sphere^3$, corresponding to the closed 
$C^4$-curves $\gamma$ and $\gamma^*$.
\qed 
\noindent

\section{The link between the Hopf-fibration and the first variation of the Willmore functional} \label{The.link}

The central result of this section is formula (\ref{DHopf.first.Variation}), the ``Hopf-Willmore-identity'', which shows that the differential of the Hopf-fibration
takes the first variation of the Willmore functional
evaluated in an immersed Hopf-torus in $\sphere^3$ into the first variation of the elastic energy functional, evaluated in the corresponding, projected curve in $\sphere^2$. 
Expressed in a more algebraic language: the Hopf-fibration and its differential transform the pair $(\Will, \nabla_{L^2}\Will)$ naturally from its effect on smooth immersed Hopf-tori in $\sphere^3$ to its effect on smooth closed regular curves in $\sphere^2$. Consequently, the  Hopf-fibration transforms flow lines of the Willmore flow in $\sphere^3$ into flow lines of the classical ``elastic energy flow'' (\ref{elastic.energy.flow}) in $\sphere^2$.
To start out, we recall some basic differential geometric 
terms in the following definition.
As in the introduction, we endow the unit $3$-sphere with
the Euclidean scalar product of $\rel^{4}$, i.e. we
set $g_{\sphere^3}:=\langle
\,\cdot\,,\,\cdot\,\rangle_{\rel^{4}}$.
\begin{definition}  \label{definitions}
	\begin{itemize}
		\item[1)] For any fixed $C^2-$immersion $G:\Sigma \longrightarrow \sphere^3$ and any smooth chart $\psi$ of an arbitrary coordinate neighbourhood $\Sigma'$ of a fixed smooth compact torus $\Sigma$, we will denote the resulting partial derivatives on $\Sigma'$ by $\partial_i$, $i=1,2$, the coefficients
		$g_{ij}:=\langle \partial_{i} G, \partial_{j} G \rangle$
		of the first fundamental form of $G$ w.r.t. $\psi$ and the associated Christoffel-symbols
		$(\Gamma_G)^m_{kl}:=g^{mj} \,\langle \partial_{kl} G,\partial_jG \rangle$ of $(\Sigma',G^*(g_{\sphere^3}))$.
		\item[2)] For any vector field $V \in C^2(\Sigma,\rel^{4})$ we define the first covariant derivatives $\nabla^G_i(V) \equiv \nabla^G_{\partial_i}(V)$, $i=1,2$, w.r.t. $G$ as the projections of the usual partial derivatives $\partial_i(V)(x)$
		of $V:\Sigma \longrightarrow \rel^{4}$ into the respective tangent spaces $T_{G(x)}\sphere^3$ of the $3$-sphere, $\forall\, x \in \Sigma$, and the second covariant derivatives by
		\begin{equation} \label{double_nabla}
			\nabla_{kl}^G(V) \equiv \nabla_k^{G} \nabla_l^{G}(V)
			:=\nabla^G_k(\nabla^G_l(V)) - (\Gamma_G)^m_{kl} \,\nabla^G_m(V).
		\end{equation}
		Moreover, we define the projections of its first derivatives into the normal bundle of the immersed surface $G(\Sigma)$ within $T\sphere^3$ by
		$$
		\nabla_{i}^{\perp_G}(V) \equiv (\nabla^G_i(V))^{\perp_G}
		:=\nabla^G_i(V)-P^{\textnormal{Tan}(G)}(\nabla^G_i(V))
		$$
		and the ``normal second covariant derivatives'' of $V$ w.r.t. $G$ by
		$$
		\nabla_{k}^{\perp_G} \nabla_{l}^{\perp_G}(V):= \nabla_k^{\perp_G} (\nabla_l^{\perp_G}(V))-(\Gamma_G)^m_{kl} \,\nabla_{m}^{\perp_G}(V).
		$$
		\item[3)]
		For a smooth regular curve $\gamma:[a,b] \to \sphere^2$ and a smooth tangent vector field $W$ on $\sphere^2$ along $\gamma$ we will denote by $\nabla_{\partial_x\gamma}(W)$ 
		the classical covariant derivative of the vector field $W$ w.r.t. the tangent vector field $\partial_x\gamma$ along $\gamma$ - being projected into $T\sphere^2$ -
		moreover by $\nabla_{\frac{\gamma'}{|\gamma'|}}(W)$ the covariant derivative of the tangent vector field $W$ along $\gamma$ w.r.t. the unit tangent vector
		field $\frac{\partial_x\gamma}{|\partial_x\gamma|}$ 
		along $\gamma$, and by $\nabla^{\perp}_{\frac{\gamma'}{|\gamma'|}}(W)$ 
		the orthogonal projection of the tangent vector
		field $\nabla_{\frac{\gamma'}{|\gamma'|}}(W)$ into the normal bundle of the curve $\gamma$ within $T\sphere^2$.   
		\item[4)] For any fixed closed regular curve 
		$\gamma \in C^{\infty}_{\textnormal{reg}}
		(\sphere^1,\sphere^2)$
		we denote by $\Gamma(\gamma^*T\sphere^2)$ the $\rel$-vector space of smooth sections of the pullback bundle $\gamma^*T\sphere^2$, i.e. of 
		all vector fields 
		$\eta \in C^{\infty}(\sphere^1,\rel^3)$ meeting the additional condition:
		$$
		\eta(x) \in T_{\gamma(x)}\sphere^2 
		\quad \textnormal{in every} \,\,x \in \sphere^1. 
		$$ 
		\item[5)] Finally, we define the subspace 
		$\Gamma^{\perp}(\gamma^*T\sphere^2)$ of the vector space $\Gamma(\gamma^*T\sphere^2)$ by the requirement, to consist of all those sections 
		$\eta \in \Gamma(\gamma^*T\sphere^2)$, 
		which are normal along the prescribed path $\gamma$, i.e. which additionally satisfy 
		$\langle \eta(x), \partial_x\gamma(x) \rangle_{\rel^3}=0$ in every
		$x \in \sphere^1$.
		\qed
	\end{itemize} 
\end{definition}
\noindent
Preparing ourselves for the fundamental Proposition \ref{Hopf.Willmore.prop} below, we firstly need 
the following elementary proposition, which results 
from straight forward computations on the basis of 
Section 2.1 in \cite{Dall.Acqua.Pozzi.2018}.
\begin{proposition} \label{compute.the.operator}
	The $L^2$-gradient of the elastic energy
	$\Wil(\gamma):= \int_{\sphere^1} 1 + |\vec{\kappa}_{\gamma}|^2 \, d\mu_{\gamma}$, with 
	$\vec{\kappa}_{\gamma}$ as in (\ref{Curva.vector}) below, evaluated in an arbitrary closed curve $\gamma \in
	C^{\infty}_{\textnormal{reg}}(\sphere^1,\sphere^2)$, reads exactly:
	\begin{equation}  \label{first.variation}
	\nabla_{L^2} \Wil(\gamma)(x) =
	2 \,\Big{(}\nabla_{\frac{\gamma'}{|\gamma'|}}^{\perp}\Big{)}^2
	(\vec \kappa_{\gamma})(x)
	+ |\vec{\kappa}_{\gamma}|^2(x) \,\vec{\kappa}_{\gamma}(x)
	+ \vec{\kappa}_{\gamma}(x),  \quad \textnormal{for} \,\, x\in \sphere^1.
	\end{equation}
	Using the abbreviation 
	$\partial_s\gamma:=\frac{\partial_x \gamma}{|\partial_x \gamma|}$ for the partial derivative of $\gamma$ normalized by arc-length, the leading term on the right hand side of equation (\ref{first.variation}) reads: 
	\begin{eqnarray}  \label{split.operator}
	\Big{(}\nabla_{\frac{\gamma'}{|\gamma'|}}^{\perp}\Big{)}^2(\vec{\kappa}_{\gamma})(x)
	=\Big{(}\nabla_{\partial_s\gamma}^{\perp}\Big{)}^2
	\big{(}(\partial_{ss}\gamma)(x) - \langle \gamma(x), \partial_{ss}
	\gamma(x) \rangle \,\gamma(x) \big{)}       \nonumber \\
	= (\partial_s)^{4}(\gamma)(x)
	- \langle (\partial_s)^{4}(\gamma)(x),\gamma(x) \rangle \,\gamma(x)           \qquad \\
	- \langle (\partial_s)^{4}(\gamma)(x),\partial_s\gamma(x) \rangle \,\partial_s\gamma(x)
	+ |(\nabla_{\partial_s\gamma})^2(\gamma)(x)|^2 \,\partial_{ss}\gamma(x),      \nonumber
	\end{eqnarray}
    for \,$x\in \sphere^1$. The fourth normalized 
    derivative $(\partial_s)^{4}(\gamma)\equiv
	\Big{(} \frac{\partial_x}{|\partial_x \gamma|}\Big{)}^4(\gamma)$ is non-linear w.r.t. 
	$\gamma$, and at least its leading term can be computed 
	in terms of ordinary partial derivatives of $\gamma$:
	\begin{eqnarray}  \label{split.operator.x}
		(\partial_s)^4(\gamma)
		=\frac{(\partial_x)^4(\gamma)}{|\partial_x\gamma|^4}
		- \frac{1}{|\partial_x\gamma|^4}\,
		\Big{\langle} (\partial_x)^4(\gamma), \frac{\partial_x\gamma}{|\partial_x\gamma|} \Big{\rangle}\, \frac{\partial_x\gamma}{|\partial_x\gamma|}\nonumber \\
		+ \,C((\partial_x)^2(\gamma),\partial_x(\gamma))  \cdot  (\partial_x)^3(\gamma) \\
		+ \,\textnormal{rational expressions which only involve} \,\,(\partial_x)^2(\gamma)
		\,\, \textnormal{and} \,\, \partial_x(\gamma),
		\nonumber
	\end{eqnarray}
	where $C:\rel^6 \longrightarrow \textnormal{Mat}_{3,3}(\rel)$ is a
	$\textnormal{Mat}_{3,3}(\rel)$-valued function, whose components are rational functions in $(y_1,\ldots,y_6) \in \rel^6$.  
	\qed    	
\end{proposition}
\noindent
Now we are ready to collect some basic differential geometric formulae in Proposition \ref{Hopf.Willmore.prop} for an arbitrary simple immersion $F$ mapping a compact torus $\Sigma$ onto the Hopf-torus $\pi^{-1}(\textnormal{trace}(\gamma))$
corresponding to some closed smooth regular curve
$\gamma:\sphere^1 \to \sphere^2$; see here 
Definition \ref{General.lift}. We shall understand below
in detail, that the following proposition paves the path to the decisive Proposition \ref{correspond.flows} and that therefore Proposition \ref{Hopf.Willmore.prop} constitutes 
the technical foundation for the proofs of 
our main results in Theorem \ref{convergence}.
\begin{proposition}  \label{Hopf.Willmore.prop}
	Let $F:\Sigma \longrightarrow \sphere^3$ be an immersion which maps the compact torus $\Sigma$ simply onto some Hopf-torus in $\sphere^3$, and let $\gamma:\sphere^1 \longrightarrow \sphere^2$ be a smooth regular parametrization of the closed curve 
	$\textnormal{trace}(\pi \circ F)$. Let moreover
	\begin{equation}  \label{Curva.vector}
	\vec{\kappa}_{\gamma}:= - \frac{1}{|\gamma'|^2} \,
	(\overline{\gamma'} \cdot \nu'_{\gamma}) \cdot \nu_{\gamma}
	\end{equation}
	be the curvature vector along the curve $\gamma$,
	for a unit normal field $\nu_{\gamma}$ along the trace of
	$\gamma$, and $\kappa_{\gamma}:=\langle \vec{\kappa}_{\gamma},\nu_{\gamma} \rangle_{\rel^3}$ the signed curvature along $\gamma$.
	Then there is some $\varepsilon =\varepsilon(F,\gamma)>0$,
	such that for an arbitrarily fixed point $s^* \in \sphere^1$ the following differential-geometric 
	formulae hold for the immersion $F$:
	\begin{equation}  \label{second.fundam.}
	A_F(\eta_F(s)) = 
	N_F(\eta_F(s))\, \left(\begin{array}{c}   
	2 \kappa_{\gamma}(s)  \quad  1     \\
	1  \qquad \quad \, 0 
    \end{array}  \right)                    
    \end{equation}  
	where $\eta_F:\sphere^1\cap B_{\varepsilon}(s^*) \to \Sigma$ denotes an arbitrary horizontal smooth lift of
	$\gamma\lfloor_{\sphere^1\cap B_{\varepsilon}(s^*)}$ w.r.t. the fibration $\pi \circ F$, as introduced in Lemma \ref{closed.lifts} below, and $N_F$ denotes a fixed unit normal field along the immersion $F$. This implies
	\begin{equation} \label{H}
	\vec H_F(\eta_F(s)) = \textnormal{trace} A_F(\eta_F(s))
	= 2 \kappa_{\gamma}(s)\, N_F(\eta_F(s))
	\end{equation}
	$\forall \, s  \in \sphere^1\cap B_{\varepsilon}(s^*)$,
	for the mean curvature vector of $F$ and also
	\begin{equation}  \label{second.fundam.trace.free}
	A^0_F(\eta_F(s)) =  N_F(\eta_F(s))\, 
	\left(  \begin{array}{c}   
	\kappa_{\gamma}(s) \qquad    1           \\
	\quad	1  \quad   -\kappa_{\gamma}(s)
	\end{array}  \right)                     
    \end{equation}   
	and consequently $|A^0_F|^2(\eta_F(s)) = 2 (\kappa_{\gamma}(s)^2+1)$, and also 
	\begin{eqnarray}  \label{Q.trace.free.H}
	Q(A^0_F)(\vec H_F)(\eta_F(s)) 
	= 4 \,(\kappa_{\gamma}^3(s) + \kappa_{\gamma}(s)) \, N_F(\eta_F(s)),              \\ \label{normal.laplace.H}
	\triangle^{\perp}_F(\vec H_F)(\eta_F(s))
	=  8\, \Big{(} \nabla_{\frac{\gamma'}{|\gamma'|}}\Big{)}^2
	(\kappa_{\gamma})(s) \, N_F(\eta_F(s)),
	\end{eqnarray}
	and finally for the traced sum of all 
	covariant derivatives of $A_F$ of order $k \in \nat$:
	\begin{equation}  \label{all.derivatives}
	|(\nabla^{\perp_F})^k(A_F)(\eta_F(s))|^2
	=  2^{2+2k} \, \Big{|} \Big{(} \nabla^{\perp}_{\frac{\gamma'}{|\gamma'|}}\Big{)}^k
	(\vec \kappa_{\gamma})(s)\Big{|}^2
	\end{equation}
	$\forall \,s \in \sphere^1\cap B_{\varepsilon}(s^*)$.
	In particular, we derive
	\begin{eqnarray} \label{MIWF.Hopf.Tori}
	\nabla_{L^2} \Will(F)(\eta_F(s)) 
	= 2 \,\Big{(} 2\,\Big{(} \nabla_{\frac{\gamma'}{|\gamma'|}}\Big{)}^2
	(\kappa_{\gamma})(s) + \kappa_{\gamma}^3(s) +
	\kappa_{\gamma}(s) \Big{)} \, N_F(\eta_F(s)), \qquad
	\end{eqnarray}
	and the ``Hopf-Willmore-identity'':
	\begin{eqnarray} \label{DHopf.first.Variation}
	D\pi_{F(\eta_F(s))}.\Big{(} 
	\nabla_{L^2} \Will(F)(\eta_F(s)) \Big{)}   
	=  4  \,\Big{(} 2\,\Big{(} \nabla^{\perp}_{\frac{\gamma'}{|\gamma'|}}\Big{)}^2
	(\vec{\kappa}_{\gamma})
	+ |\vec{\kappa}_{\gamma}|^2 \vec{\kappa}_{\gamma}
	+ \vec{\kappa}_{\gamma} \Big{)}(s)   \nonumber   \\
	\equiv  4 \, \nabla_{L^2} \Wil(\gamma)(s)  \qquad  
	\end{eqnarray}
	$\forall \,s \in \sphere^1\cap B_{\varepsilon}(s^*)$,
	where there holds $\pi \circ F \circ \eta_F = \gamma$ on
	$\sphere^1\cap B_{\varepsilon}(s^*)$; see Lemma 
	\ref{closed.lifts}. Finally, we will prove that
	\begin{equation}  \label{Will.Wil}
	\Will(F) \equiv \int_{\Sigma} 1 + \frac{1}{4} \,
	|\vec H_F|^2 \, d\mu_{F}
	=\pi  \,\int_{\sphere^1} 1 + \,|\kappa_{\gamma}|^2 \, d\mu_{\gamma} \equiv \pi  \,\Wil(\gamma), 
	\end{equation}
	and 
	\begin{equation}  \label{First.Var.Will.Wil}
	\int_{\Sigma} |\nabla_{L^2}\Will(F)|^2 \, d\mu_{F} 
	= 4\pi  \,\int_{\sphere^1} |\nabla_{L^2}\Wil(\gamma)|^2 \, d\mu_{\gamma}.                               
	\end{equation}
\end{proposition}
\noindent
\proof: Let $2L$ denote the length of the curve
$\textnormal{image}(\pi \circ F) \subset \sphere^2$.
Without loss of generality we may require $\gamma$ 
to perform only one loop through its trace. 
Hence, we can parametrize
the trace of $\gamma$ on $\rel/L\ganz$ in such 
a way that it has constant speed $2$ on $\rel/L\ganz$. We therefore assume during 
the entire proof of this proposition that $\gamma$ is 
defined on $\rel/L\ganz$, satisfies 
$|\gamma'| \equiv 2$ on $\rel/L\ganz$ and
that $s^*$ is a fixed point in $[0,L]$.
We infer from Lemma \ref{closed.lifts} the existence
of horizontal smooth lifts $\eta:[0,L] \to
\pi^{-1}(\textnormal{trace}(\gamma))$ of $\gamma$	
w.r.t. $\pi$, and we infer from point (4)
of Lemma \ref{Hopf.Tori}, that every such lift
$\eta$ must have constant speed $|\eta'| \equiv 1$ on $[0,L]$ 
and intersects each fiber of $\pi$ exactly once.
Hence, also statements (1)-(3) of Lemma \ref{Hopf.Tori}
can be applied in the sequel. In order to simplify our computations, we choose the explicit parametrization
$X(s,\varphi):=e^{i\varphi} \,\eta(s)$, for
$(s,\varphi) \in [0,L] \times [0,2\pi]$, from formula
(\ref{Hopf.Torus}) of $\pi^{-1}(\textnormal{trace}(\gamma))$,
which gives us the opportunity to apply many ideas 
arising in Sections 2 and 3 of \cite{Pinkall}. 
Firstly, we infer from Lemma
\ref{Hopf.Tori}, that the vector
$N_X(s,\varphi) = i \,u(s) \,e^{-i \varphi} \,\eta(s)$
is orthogonal to $\partial_sX(s,\varphi)$ and 
$\partial_{\varphi}X(s,\varphi)$
within $T_{X(s,\varphi)}\sphere^3$ and  
therefore a unit normal field along $\pi^{-1}(\textnormal{trace}(\gamma))$.
\footnote{See here also formula (20) in \cite{Pinkall}.} 
Now, we apply formula (\ref{Hopf.differential})
with $q:=X(s,\varphi) \in \pi^{-1}(\textnormal{trace}(\gamma))$ and $v:=N_X(s,\varphi)$ and compute similarly to 
(\ref{gamma.speed}):
\begin{eqnarray}  \label{Hopf.unit.normal}
	D\pi_{X(s,\varphi)}(N_X(s,\varphi))
	=I(N_X(s,\varphi)) \cdot X(s,\varphi)
	+ I(X(s,\varphi)) \cdot N_X(s,\varphi)  \quad\\
	= I(i u(s) \,e^{-i \varphi} \eta(s)) \cdot
	(e^{i\varphi} \,\eta(s))
	+I(e^{i\varphi} \,\eta(s)) \cdot 
	(i u(s) \,e^{-i \varphi} \, \eta(s))
	\nonumber\\
	=- \tilde \eta(s)\, e^{i \varphi} \,\tilde u(s)\, i
	e^{i\varphi} \,\eta(s)
	+\tilde \eta(s)\,e^{-i\varphi}  \, i u(s) \,e^{-i \varphi} \, \eta(s)  \nonumber \\
	=- \tilde \eta(s)\,\,\tilde u(s) i\,\eta(s)
	+\tilde \eta(s) i u(s) \, \eta(s)
	= 2 \tilde \eta(s)\, i \, u(s) \,\eta(s),   \nonumber
\end{eqnarray}
where we used that $i \,u(s) = - u(s) \,i$ and 
$\tilde u(s)=u(s)$. Moreover, we compute:
\begin{equation}  \label{normal.section.gamma}
	\overline{\tilde{\eta}(s) \, i \,u(s) \, \eta(s)}
	\, \tilde\eta(s)\, u(s)\, \eta(s)
	= \bar \eta(s) \,i \,\bar u(s)  \,\bar{\tilde{\eta}}(s) \,
	\tilde\eta(s)\, u(s) \,\eta(s)
	= \bar\eta(s) \, i\, \eta(s),
\end{equation}
and we know that $\Re(\bar\eta i \eta)\equiv 0$. 
Combining this with the fact that
$\gamma'(s) = 2 \tilde\eta(s)\, u(s) \,\eta(s)$ 
is tangential at the curve 
$\gamma = \pi\circ \eta$ in its point $\gamma(s)$,  
formula (\ref{normal.section.gamma}) shows that
there holds $\tilde{\eta} \, i \,u \, \eta \in 
\Gamma^{\perp}(\gamma^*T\sphere^2)$, being here 
even a normal section along $\gamma$ of constant length 
one. Hence, (\ref{Hopf.unit.normal})
shows that the differential $D\pi$ of the Hopf-fibration
maps the one-dimensional normal bundle of
the immersion $X$ within $T\sphere^3$ isomorphically onto the normal bundle of the curve $\gamma=\pi \circ \eta$ 
within $T\sphere^2$. Moreover, we verify by means of formula (\ref{gamma.speed}) that
\begin{eqnarray}  \label{Hopf.partial.s.X}
	D\pi_{X(s,\varphi)}(\partial_s X(s,\varphi))
	=I(\partial_s X(s,\varphi)) \cdot X(s,\varphi)
	+ I(X(s,\varphi)) \cdot \partial_s X(s,\varphi) \qquad \\
	= I(e^{i \varphi} \, u(s) \,\eta(s)) \cdot
	(e^{i\varphi} \,\eta(s))
	+I(e^{i\varphi} \,\eta(s)) \cdot (e^{i \varphi} \, u(s) \,\eta(s))
	\nonumber\\
	= \tilde \eta(s)\, \tilde u(s)\, e^{-i \varphi} \,
	e^{i\varphi} \,\eta(s)
	+ \tilde \eta(s)\,e^{-i\varphi}  \, e^{i \varphi} \, u(s) \, \eta(s) 
	= 2 \tilde \eta(s)\,  u(s) \,\eta(s) =\gamma'(s)   \nonumber
\end{eqnarray}
and
\begin{eqnarray}  \label{Hopf.partial.varphi.X}
	D\pi_{X(s,\varphi)}(\partial_{\varphi} X(s,\varphi))
	=I(\partial_{\varphi} X(s,\varphi)) \cdot X(s,\varphi)
	+ I(X(s,\varphi)) \cdot \partial_{\varphi} X(s,\varphi) \qquad \\
	= I(i e^{i \varphi} \, \eta(s)) \cdot
	(e^{i\varphi} \,\eta(s))
	+I(e^{i\varphi} \,\eta(s)) \cdot (i\, e^{i \varphi} \,\eta(s))                                 \nonumber\\
	= - \tilde \eta(s)\, i\, e^{-i \varphi} \,
	e^{i\varphi} \,\eta(s)
	+ \tilde \eta(s)\,e^{-i\varphi}  \,i\, e^{i \varphi}  \, \eta(s) \equiv 0.  \nonumber
\end{eqnarray}
Formulae (\ref{Hopf.partial.s.X}) and (\ref{Hopf.partial.varphi.X}) show, that $D\pi$ maps the tangent bundle of the immersion $X$ in $T\sphere^3$ 
onto the tangent bundle of the closed curve $\gamma$ in $T\sphere^2$ with a necessarily one-dimensional kernel.
Again using Lemma \ref{Hopf.Tori}, $i\,u(s) = -u(s)\, i$ 
and $u^2 \equiv -1$, we obtain as in formula (21) of 
\cite{Pinkall}:
\begin{eqnarray}  \label{Nt}
\partial_s N_X(s,\varphi)
=-2 \kappa_{\gamma}(s) \,\partial_s X(s,\varphi) - \partial_{\varphi} X(s,\varphi), \qquad \nonumber \\
\partial_{\varphi}N_X(s,\varphi)
= - i u(s) i e^{-i \varphi} \, \eta(s)
= - u(s) e^{-i \varphi} \, \eta(s) = -\partial_s X(s,\varphi),
\qquad
\end{eqnarray}
where we have introduced the curvature function
$\kappa_{\gamma}(s)$ of the curve
$\gamma = \pi \circ \eta$ in the point $s\in [0,L]$ 
by the relation:
\begin{equation}  \label{curvature.gamma}
u'(s) = 2i \kappa_{\gamma}(s) \,u(s),
\end{equation}
as in formula (22) of \cite{Pinkall}.
Since $u$ has no zeroes, equation (\ref{curvature.gamma}) 
defines a unique smooth function on $[0,L]$,
and one can easily verify, using the formula
$
(\tilde{\eta} \, i \,u \, \eta)'
=\tilde \eta \,\tilde u \, i \,u\, \eta
+\tilde \eta \,i\, u' \,\eta
+ \tilde \eta \, i \,u^2 \,\eta
$
and the fact that $u^2 \equiv -1$, that this function
$\kappa_{\gamma}$ exactly coincides with the signed curvature
function of the curve $\gamma$, as defined below definition (\ref{Curva.vector}). Hence using (\ref{normal.section.gamma}), we can write 
$\vec{\kappa}_{\gamma} = \kappa_{\gamma}
\, \tilde{\eta} \, i \,u \, \eta$.
Moreover, on account of Lemma \ref{Hopf.Tori} we know
that $\{\partial_{s} X, \partial_{\varphi} X\}$
is an orthonormal frame along $\pi^{-1}(\textnormal{trace}(\gamma))$, implying that   
$(g_X)_{ij}=\delta_{ij}$, and thus (\ref{Nt}) 
yields the coefficients of the $2$nd fundamental form $A_X$ w.r.t. the unit normal $N_X$: 
\begin{equation}  \label{A.X}
A_X(s,\varphi) = N_X(s,\varphi)\, \left(\begin{array}{c}   
2 \kappa_{\gamma}(s)  \quad  1     \\
1  \qquad \quad \, 0 
\end{array}  \right).                     
\end{equation} 
This yields
\begin{eqnarray} \label{mean.curvature.form}
\vec H_X(s,\varphi) \equiv \textnormal{trace} A_X(s,\varphi)
=  2 \,\kappa_{\gamma}(s) \,N_X(s,\varphi)
\end{eqnarray}
for the mean curvature vector of $X$ and also
\begin{equation}  \label{A.0.X}
A^0_X(s,\varphi) =  N_X(s,\varphi)\,  \left(\begin{array}{c}   
\kappa_{\gamma}(s) \qquad    1  \\
\quad 1  \qquad   -\kappa_{\gamma}(s)
\end{array}  \right),                       
\end{equation}   
in particular
\begin{equation} \label{A0.squared}
|A^0_X|^2(s, \varphi) = 2 \,(\kappa_{\gamma}(s)^2+1),
\end{equation}
and \,$Q(A^0_X)(\vec H_X)(s,\varphi) =
(A^0_X)_{ij} \langle (A^0_X)_{ij},\vec H_X \rangle(s,\varphi)
= 4 (\kappa_{\gamma}^3 + \kappa_{\gamma})(s) \, N_X(s,\varphi)$. Recalling now that 
$N_X(s,\varphi)= i \,u(s)\,e^{-i \varphi}\,\eta(s)$ 
and formula (\ref{mean.curvature.form}), we have:
\begin{eqnarray}  \label{nabla.t.H}
	\nabla_s^{\perp_X} \vec H_X(s,\varphi)
	= \nabla_s^{\perp_X} (2\kappa_{\gamma}(s) \,i u(s) 
	e^{-i \varphi} \, \eta(s))  \qquad \nonumber \\
	=  \Big{(}
	2 \kappa_{\gamma}'(s) \, i \,u(s) \,e^{-i \varphi} \, \eta(s)
	+ 2 \kappa_{\gamma}(s) \, i u'(s)\, e^{-i \varphi} \, \eta(s)
	+  2 \kappa_{\gamma}(s) \, i u(s) \, e^{-i \varphi} \, u(s) \eta(s)\Big{)}^{\perp_X} \nonumber \\
	=  2\, \kappa_{\gamma}'(s) \, N_X(s,\varphi). \qquad
\end{eqnarray}
In order to achieve this result, one has to derive
from $|u|^2=1$, $|\eta|^2=1$, $\eta'=u\,\eta$ and
$u(s) \in \textnormal{Span}\{j,k\}$, $\forall \,s \in [0,L]$:
$$
0=\bar u'\,u + \bar u\, u'=\bar u'\,u + \overline{\bar u'\,u}
=2 \Re(\bar u'\,u ) = 2 \Re(\bar u\,u')
$$
and thus $(\bar u\,u')(s) = - i \,\det(u,u')(s)$, 
and then compute:
\begin{eqnarray*}
	\Re\big{(}\overline{N_X(s,\varphi)} \, i u'(s) 
	e^{-i \varphi} \, \eta(s)\big{)}
	=\Re\big{(}\bar \eta(s) e^{i\varphi}\,\bar u(s)(-i) \,i u'(s) e^{-i \varphi} \, \eta(s)\big{)}  \\
	= -\Re\bar{(}\bar \eta(s)\, e^{i\varphi}i 
	\det(u(s), u'(s))\,e^{-i \varphi} \, \eta(s)\big{)}      
	= -\det(u(s),u'(s))\, \Re(\bar \eta(s) i \eta(s)) \equiv 0,
\end{eqnarray*}
also using that $\Re(\bar \eta \,i \,\eta) \equiv 0$, and 
then:
$$
\Re\big{(}\overline{N_X(s,\varphi)} \,i u(s) e^{-i \varphi} \, u(s) \eta(s) \big{)}
=\Re\big{(}\bar \eta(s) u(s) \eta(s)\big{)}
=\Re\big{(}\bar \eta(s) \eta'(s)\big{)}\equiv 0
$$
for $(s,\varphi)\in [0,L] \times [0,2\pi]$.
Moreover, combining formulae (\ref{Nt}) and  
(\ref{mean.curvature.form}) we obtain:
\begin{eqnarray}  \label{nabla.phi.H}
\nabla_{\varphi}^{\perp_X} \vec H_X(s,\varphi)
= 2 \kappa_{\gamma}(s) \, 
\Big{(} \partial_{\varphi} N_X(s,\varphi) \Big{)}^{\perp_X}
= -2 \kappa_{\gamma}(s) \, \Big{(}\partial_s X(s,\varphi)\Big{)}^{\perp_X} \equiv 0. \qquad
\end{eqnarray}
This immediately implies also
\begin{eqnarray*}
\nabla_{\varphi}^{\perp_X} \nabla_{\varphi}^{\perp_X}
\vec H_X(s,\varphi) =0  \quad \textnormal{and} \quad   
\nabla_{\varphi}^{\perp_X} \nabla_{s}^{\perp_X}
\vec H_X(s,\varphi) = 0 =
\nabla_{s}^{\perp_X}\nabla_{\varphi}^{\perp_X} 
\vec H_X(s,\varphi)
\end{eqnarray*}
for $(s,\varphi)\in [0,L] \times [0,2\pi]$. 
Finally, we derive from (\ref{nabla.t.H}) by analogy:
\begin{eqnarray*}
\nabla_{s}^{\perp_X} \nabla_{s}^{\perp_X}
\vec H_X(s,\varphi)
= \nabla_{s}^{\perp_X}(2 \,\kappa_{\gamma}'(s) N_X(s,\varphi))
= 2 \kappa_{\gamma}''(s) \, N_X(s,\varphi),
\end{eqnarray*}
for $(s,\varphi)\in [0,L] \times [0,2\pi]$, 
and therefore also:
\begin{equation}  \label{Nabla.k.A.X}
(\nabla_{s}^{\perp_X})^k A_X(s,\varphi) =  N_X(s,\varphi)\,  \left(  \begin{array}{c}   
		2 \kappa_{\gamma}^{(k)}(s)  \quad  0 \, \\
		0  \qquad \quad\,\,\,\,  0
\end{array}  \right)                    
\end{equation}  
for $k \in \nat$, which immediately yields
$|(\nabla^{\perp_X})^k (A_X)(s,\varphi)|^2
=  4 \,(\kappa_{\gamma}^{(k)}(s))^2$ on $[0,L] \times [0,2\pi]$. Together with $(g_X)_{ij}=\delta_{ij}=(g_X)^{ij}$ 
on $[0,L] \times [0,2\pi]$, we instantly arrive at the
expression for the normal laplacian of $\vec H_X$: 
\begin{equation}  \label{pseudo.Laplace.H}
\triangle^{\perp}_X (\vec H_X)(s,\varphi)
=\nabla_{s}^{\perp_X} \nabla_{s}^{\perp_X} \vec H_X(s,\varphi)
= 2 \kappa_{\gamma}''(s) N_X(s,\varphi),
\end{equation} 
see here also Section 2 in \cite{Jakob_Moebius_2016}. 
Combining formulae (\ref{A0.squared}) 
and (\ref{pseudo.Laplace.H}), we obtain:
\begin{eqnarray}  \label{Willmore.1st.step}
	\nabla_{L^2} \Will(X)(s,\varphi)
	= \Big{(} \kappa_{\gamma}''(s) + 2 \,\kappa_{\gamma}^3(s)
	+ 2\, \kappa_{\gamma}(s) \Big{)} \, N_X(s,\varphi) 
\end{eqnarray}
for $(s,\varphi) \in [0,L] \times [0,2\pi]$.
Now, the curvature vector $\vec{\kappa}_{\gamma}$ and 
its covariant derivatives 
$\Big{(} \nabla_{\frac{\gamma'}{\gamma'}|}\Big{)}^k
(\vec{\kappa}_{\gamma})$ w.r.t. the unit tangent vector 
field $\frac{\gamma'}{|\gamma'|}$ along $\gamma$ is 
invariant w.r.t. smooth reparametrization 
of $\gamma$, unlike the usual derivatives $\kappa_{\gamma}^{(k)}$ of the signed curvature $\kappa_{\gamma}$ w.r.t. the original parameter 
$s\in [0,L]$. Since the path $\gamma =\pi\circ \eta$ 
was assumed to have constant speed $2$, we thus obtain here:
\begin{eqnarray}  \label{reparametrization}
\Big{(} \nabla^{\perp}_{\frac{\gamma'}{|\gamma'|}}\Big{)}^k
(\vec{\kappa}_{\gamma})(s)
= \frac{1}{2^k} \,\kappa_{\gamma}^{(k)}(s) \,
\tilde{\eta}(s) \, i \,u(s) \, \eta(s)
\equiv\frac{1}{2^k} \,\kappa_{\gamma}^{(k)}(s) \, \nu_{\gamma}(s),
\end{eqnarray}
for every $k \in \nat$ and $\forall \, s\in [0,L]$.
Moreover, taking Lemma \ref{closed.lifts} into account, we know that for every horizontal smooth lift $\eta$ of $\gamma$ w.r.t. $\pi$ and for every simple parametrization $F:\Sigma \longrightarrow \pi^{-1}(\textnormal{trace}(\gamma))$
there is a horizontal smooth lift $\eta_F:\rel/L\ganz \cap B_{\epsilon}(s^*) \to \Sigma$ of $\gamma\lfloor_{\rel/L\ganz \cap B_{\epsilon}(s^*)}$ w.r.t. $\pi \circ F$, i.e. $\eta_F$ satisfies $\eta = F \circ \eta_F$ on $\rel/L\ganz \cap B_{\epsilon}(s^*)$.
Hence, taking also into account the invariance of traced tensors w.r.t. smooth local coordinate transformations and also
the invariance of the curvature vector $\vec{\kappa}_{\gamma}$
w.r.t. smooth reparametrizations of the path $\gamma$,
we can immediately derive from formulae
(\ref{A.X})--(\ref{reparametrization}) the
assertions (\ref{second.fundam.})--(\ref{MIWF.Hopf.Tori}).
Moreover, we note that formula (\ref{Hopf.unit.normal}) describes a geometric property of the pair $(\pi,D\pi)$, which is thus independent of the chosen parametrization of the Hopf-torus $\pi^{-1}(\textnormal{trace}(\gamma))$.
Hence, a combination of formulae (\ref{MIWF.Hopf.Tori}), (\ref{Hopf.unit.normal}) and (\ref{reparametrization}) 
immediately yields assertion (\ref{DHopf.first.Variation}).\\
Finally, since the standard parametrization
$X:[0,L] \times [0,2\pi] \longrightarrow \pi^{-1}(\textnormal{trace}(\gamma))$
in formula (\ref{Hopf.Torus}) covers $\pi^{-1}(\textnormal{trace}(\gamma))$ only once,
as pointed out in Remark \ref{Fundamental.gruppe.Homologie}, 
since $F$ has mapping degree $\pm 1$, and since $\gamma$ has constant speed $2$ on $\rel/L\ganz $, we  
obtain from the classical area-formula, formula (\ref{mean.curvature.form}) and Lemma \ref{Hopf.Tori}:
\begin{eqnarray*}
\Will(F)=|\deg(F)| \, \Will(X)
=   \int_{0}^{2\pi} \int_{0}^L
1 + \frac{1}{4} |2 \kappa_{\gamma}(s) |^2 \, ds\, d\varphi 
\nonumber \\
= \,2 \pi \, \, \int_{0}^L 1+|\kappa_{\gamma}|^2\,ds
=  \pi\, \int_{\rel/L\ganz}
1 + |\kappa_{\gamma}|^2 \, d\mu_{\gamma}
\equiv \pi  \,\Wil(\gamma).
\end{eqnarray*}
Similarly, we can infer from formulae (\ref{Willmore.1st.step}) and (\ref{reparametrization}), 
that there holds: 
\begin{eqnarray*}
	\int_{\Sigma} |\nabla_{L^2}\Will(F)|^2 \, d\mu_{F}
	= |\deg(F)| \, 	\int_0^{2\pi} \,\int_{0}^L  
	|\nabla_{L^2}\Will(X)|^2(s,\varphi) \, ds\, d\varphi                                           \\
	= 2\pi \,\int_0^L  \Big{|} \kappa_{\gamma}''(s) +
	2 \kappa_{\gamma}^3(s)
	+ 2 \kappa_{\gamma}(s) \Big{|}^2 \, ds            \\
	= 4\pi  \,\int_{\rel/L\ganz} \Big{|} 
	\,2\,\Big{(}\nabla^{\perp}_{\gamma'/|\gamma'|}\Big{)}^2
	(\vec{\kappa}_{\gamma})
	+ |\vec{\kappa}_{\gamma}|^2 \,\vec{\kappa}_{\gamma} 
	+ \vec{\kappa}_{\gamma} \,\Big{|}^2 \, d\mu_{\gamma}, 
\end{eqnarray*} 
which is just the asserted formula (\ref{First.Var.Will.Wil})
on account of the invariance of the differential operator 
$\nabla_{L^2}\Wil(\gamma)$ w.r.t. reparametrization 
of the curve $\gamma$.
\qed
\noindent
\begin{proposition}  \label{correspond.flows}
	Let $T>0$ be arbitrarily chosen, and let $\gamma_t: \sphere^1 \to \sphere^2$ be a smooth family of closed smooth regular curves, for $t \in [0,T]$.
	Moreover, let $F_t:\Sigma \longrightarrow \sphere^3$
	be an arbitrary smooth family of smooth immersions, which
	map some compact smooth torus $\Sigma$ simply onto the Hopf-tori $\pi^{-1}(\textnormal{trace}(\gamma_t))\subset \sphere^3$, for every $t\in [0,T]$. Then the following statement holds:\\
	The family of immersions $\{F_t\}$ moves according to the
	Willmore flow equation (\ref{Willmore.flow}) on $[0,T] \times \Sigma$ - up to smooth, time-dependent reparametrizations $\Phi_t$ with $\Phi_0=\textnormal{Id}_{\Sigma}$ - if and only if there is a smooth family $\sigma_t:\sphere^1 \to \sphere^1$ of reparametrizations with $\sigma_0=\textnormal{Id}_{\sphere^1}$, such that
	the family $\{\gamma_t\circ \sigma_t\}$ satisfies the
	``elastic energy evolution equation''
	\begin{eqnarray}  \label{elastic.energy.flow}
		\partial_t \tilde \gamma_t =
		- \,\Big{(} 2 \, \Big{(}\nabla^{\perp}_{\frac{\tilde \gamma_t'}{|\tilde \gamma_t'|}} \Big{)}^2(\vec{\kappa}_{\tilde \gamma_t})
		+ |\vec{\kappa}_{\tilde \gamma_t}|^2 \vec{\kappa}_{\tilde \gamma_t}
		+ \vec{\kappa}_{\tilde \gamma_t} \Big{)}                     
		\equiv - \,\nabla_{L^2}\Wil(\tilde \gamma_t)
	\end{eqnarray}
	on $[0,T] \times \sphere^1$, where 
	$\nabla_{L^2} \Wil$ denotes the $L^2$-gradient of 
	$\Wil(\gamma) =
	\int_{\sphere^1} 1+|\vec{\kappa}_{\gamma}|^2 \, \,d\mu_{\gamma}$.
\end{proposition}
\noindent
\proof: The easier direction of the assertion 
follows immediately from the ``Hopf-Willmore-identity''
(\ref{DHopf.first.Variation}) and the ordinary chain rule.
For, suppose a smooth family $F_t:\Sigma \longrightarrow
\pi^{-1}(\textnormal{trace}(\gamma_t))$ of smooth immersions,
which parametrize $\pi^{-1}(\textnormal{trace}(\gamma_t))$
simply, solves the Willmore equation (\ref{Moebius.flow}) 
on $[0,T] \times \Sigma$.
In this case, we apply Lemma \ref{closed.lifts} and 
choose some smooth family of horizontal smooth lifts $\eta_{F_t}:\sphere^1 \cap B_{\epsilon_t}(s^*)\longrightarrow \Sigma$ of $\gamma_t\lfloor_{\sphere^1 \cap B_{\epsilon_t}(s^*)}$,
for some arbitrary $s^*\in \sphere^1$ and some $\epsilon_t>0$,
i.e. such that $\gamma_t = \pi \circ F_t \circ \eta_{F_t}$
on $\bigcup_{t \in [0,T]} \{t\} \times (\sphere^1 \cap B_{\epsilon_t}(s^*))$, and conclude by means of formula (\ref{DHopf.first.Variation}):
\begin{eqnarray*}
	\partial_t \gamma_t(s)
	= \partial_t (\pi \circ F_t \circ \eta_{F_t})(s)
	= D\pi_{F_t(\eta_{F_t}(s))}.(\partial_t(F_t \circ \eta_{F_t})(s))                                \\
	= D\pi_{F_t(\eta_{F_t}(s))}.((\partial_tF_t) \circ \eta_{F_t})(s)
	+ D\pi_{F_t(\eta_{F_t}(s))}.\Big{(} D_{\eta_{F_t}(s)}(F_t).(\partial_t(\eta_{F_t}))(s) \Big{)}                                        \\
	= D\pi_{F_t(\eta_{F_t}(s))}.
	\Big{(} -\,\nabla_{L^2} \Will(F_t)(\eta_{F_t}(s)) \Big{)}                                  \qquad \\ 
	+ D\pi_{F_t(\eta_{F_t}(s))}.\Big{(} D_{\eta_{F_t}(s)}(F_t).(\partial_t(\eta_{F_t})(s)) \Big{)}                                        \\
	=- \,\nabla_{L^2} \Wil(\gamma_t)(s) +
	D\pi_{F_t(\eta_{F_t}(s))}.\Big{(} D_{\eta_{F_t}(s)}(F_t).(\partial_t(\eta_{F_t})(s)) \Big{)} \\
\end{eqnarray*}
on $\bigcup_{t \in [0,T]} \{t\} \times (\sphere^1 \cap B_{\epsilon_t}(s^*))$.
Now, the vector $\partial_t(\eta_{F_t})(s)$ is contained in
the tangent space $T_{\eta_{F_t}(s)}\Sigma$ touching $\Sigma$
at the point $\eta_{F_t}(s)$, for every $s \in \sphere^1 \cap B_{\epsilon_t}(s^*)$ and every $t\in [0,T]$.
Then $D_{\eta_{F_t}(s)}(F_t).(\partial_t(\eta_{F_t})(s))$
is a tangent vector of $\pi^{-1}(\textnormal{trace}(\gamma_t))$
at its point $F_t(\eta_{F_t}(s))$, and formula (\ref{Hopf.partial.varphi.X}) shows that
$D\pi_{F_t(\eta_{F_t}(s))}.\Big{(} D_{\eta_{F_t}(s)}(F_t).(\partial_t(\eta_{F_t})(s)) \Big{)}$ is a tangent vector of $\gamma_t$ in its point
$(\pi \circ F_t \circ \eta_{F_t})(s)=\gamma_t(s)$
on $\bigcup_{t \in [0,T]} \{t\} \times (\sphere^1 \cap B_{\epsilon_t}(s^*))$.
Since the vector $\nabla_{L^2} \Wil(\gamma_t)(s)$ is contained in the normal space of the curve $\gamma_t$ within $T\sphere^2$ in its point $\gamma_t(s)$ and since 
$s^* \in \sphere^1$ was arbitrarily chosen, we conclude that the family $\{\gamma_t\}$ solves the equation
\begin{equation}  \label{flow.on.curves.perp}
\Big{(} \partial_t \gamma_t(s) \Big{)}^{\perp_{\gamma_t}}
=- \,\nabla_{L^2} \Wil(\gamma_t)(s)
\end{equation}
$\forall \,(t,s) \in [0,T] \times \sphere^1$,
where $(V_t)^{\perp_{\gamma_t}}$ denotes the projection of vector fields $V_t \in \Gamma(\gamma_t^*T\sphere^2)$ 
along $\gamma_t$ into the subspace $\Gamma^{\perp}(\gamma_t^*T\sphere^2)$, for every 
fixed $t \in [0,T]$.
Now, using the fact that the smooth family $\{\gamma_t\}$
solves equation (\ref{flow.on.curves.perp}) on $[0,T] \times \sphere^1$ and that it consists of smooth, closed and regular curves only, we can follow the lines of the author's article \cite{Jakob_Moebius_2016}, p. 1177, solving a certain ODE in order to construct an appropriate smooth family of diffeomorphisms $\sigma_t: \sphere^1 \to \sphere^1$, 
for $t \in [0,T]$, satisfying $\sigma_0=\textnormal{Id}_{\sphere^1}$, such that
the composition $\gamma_t\circ \sigma_t$ solves equation
(\ref{elastic.energy.flow}) on $[0,T] \times \sphere^1$,
just as asserted. Vice versa, if there holds
\begin{equation}  \label{Moebius.flow.curves}
\partial_t \gamma_t(s) = 
- \, \nabla_{L^2} \Wil(\gamma_t)(s),
\end{equation}
$\forall \, (t,s)\in [0,T] \times \sphere^1$,
for some family of smooth closed regular curves
$\gamma_t:\sphere^1 \to \sphere^2$, $t\in [0,T]$, then
- again using Lemma \ref{closed.lifts} - we first
choose some family of horizontal smooth lifts
$\eta_{F_t}:\sphere^1 \cap B_{\epsilon_t}(s^*) \longrightarrow \Sigma$ of $\gamma_t\lfloor_{\sphere^1 
\cap B_{\epsilon_t}(s^*)}$ w.r.t. $\pi \circ F_t$ on $[0,T]$ 
for an arbitrarily given smooth family of simple parametrizations $F_t:\Sigma \longrightarrow \pi^{-1}(\textnormal{trace}(\gamma_t))$,
and then we compute by means of 
$\gamma_t = \pi \circ F_t \circ \eta_{F_t}$ on 
$\sphere^1 \cap B_{\epsilon_t}(s^*)$, formula (\ref{DHopf.first.Variation}) and the chain rule:
\begin{eqnarray} \label{vice.versa}
	D\pi_{F_t(\eta_{F_t}(s))}.
	\Big{(} -\,\nabla_{L^2} \Will(F_t)(\eta_{F_t}(s)) \Big{)}                  \nonumber \qquad \\
	=- \,\nabla_{L^2} \Wil(\gamma_t)(s)                
	= \partial_t (\gamma_t)(s)
	=\partial_t (\pi \circ F_t \circ \eta_{F_t})(s)     \qquad  \nonumber \\
	=D\pi_{F_t(\eta_{F_t}(s))}.((\partial_tF_t) \circ \eta_{F_t})(s) + D\pi_{F_t(\eta_{F_t}(s))}.\Big{(} D_{\eta_{F_t}(s)}(F_t).(\partial_t(\eta_{F_t})(s)) \Big{)}  \qquad
\end{eqnarray}
on $\bigcup_{t \in [0,T]} \{t\} \times (\sphere^1 \cap B_{\epsilon_t}(s^*))$ and for every family of horizontal smooth lifts $\eta_{F_t}$ of $\gamma_t\lfloor_{\sphere^1 \cap B_{\epsilon_t}(s^*)}$ w.r.t. $\pi \circ F_t$. Now, we fix some point $x\in \Sigma$ arbitrarily. For every $t \in [0,T]$, there is at least one point $s_x(t) \in \sphere^1$ satisfying
\begin{equation} \label{s.x}
\gamma_t(s_x(t)) = (\pi \circ F_t)(x).
\end{equation}
Since every path $\gamma_t$ is required to be regular,
i.e. $\gamma_t'(s)\not =0$ $\forall s \in \sphere^1$
in every $t\in [0,T]$, and since $\{\gamma_t\}$ and $\{F_t\}$
are required to be smooth families of smooth maps,
the implicit function theorem guarantees that
the solutions $s_x(t)$ to equation (\ref{s.x})
can be chosen in such a way, that $[t \mapsto s_x(t)]$ is a
smooth function mapping $[0,T]$ into $\sphere^1$.
On account of (\ref{s.x}) and Lemma \ref{closed.lifts},
we can choose the horizontal smooth lift $\eta_{F_t}$ for every
$t \in [0,T]$ in such a way, that exactly 
$\eta_{F_t}(s_x(t))=x$ holds. Inserting this into formula
(\ref{vice.versa}) we obtain
\begin{eqnarray} \label{vice.versa.2}
	D\pi_{F_t(x)}.
	\Big{(} - \nabla_{L^2} \Will(F_t)(x) \Big{)}  \quad \\ \nonumber
	=  D\pi_{F_t(x)}.\Big{(}\partial_tF_t(x)\Big{)}
	+ D\pi_{F_t(x)}.\Big{(} D_{x}(F_t).(\partial_t(\eta_{F_t})(s_x(t))) \Big{)},
\end{eqnarray}
for every fixed $x\in \Sigma$ and for $t\in [0,T]$.
Now, as we have computed in formulae 
(\ref{Hopf.unit.normal})--(\ref{Hopf.partial.varphi.X}), 
the Hopf-differential $D\pi_{F_t}$ maps the one-dimensional
normal bundle of the immersion $F_t$ in $T\sphere^3$
isomorphically onto the normal bundle of the regular
curve $\gamma_t$ within $T\sphere^2$
and the tangent bundle of the immersion $F_t$ onto
the tangent bundle of $\gamma_t$ in $T\sphere^2$,
which implies in particular that
$D\pi_{F_t(\eta_{F_t}(s))}.\Big{(} D_{\eta_{F_t}(s)}(F_t).
(\partial_t(\eta_{F_t})(s)) \Big{)}$
is a tangent vector of the curve $\pi \circ F_t \circ \eta_{F_t}= \gamma_t$ in its point $\gamma_t(s)$.
Hence, since $\nabla_{L^2} \Will(F_t)$ is a section of the normal bundle of the immersion $F_t$, a comparison of
the normal components in (\ref{vice.versa.2}) particularly
yields the following modified flow equation:
\begin{equation} \label{MIWF.normal}
\Big{(} \partial_t F_t \Big{)}^{\perp_{F_t}}(x)
=- \, \nabla_{L^2} \Will(F_t)(x)
\end{equation}
$\forall (t,x) \in [0,T] \times \Sigma$, where
$(V_t)^{\perp_{F_t}}$ denotes the projections of vector fields
$V_t \in \Gamma(F_t^*T\sphere^3)$ along the immersion 
$F_t$ onto the subspace $\Gamma^{\perp}(F_t^*T\sphere^3)$ 
of smooth sections of the normal bundle of $F_t$,
for every $t \in [0,T]$. Now we can follow the lines of 
the author's article \cite{Jakob_Moebius_2016}, p. 1177,
solving a certain system of ODE's, in order to infer
from the facts that the smooth family $\{F_t\}$
solves equation (\ref{MIWF.normal}) on $[0,T] \times \Sigma$
and that it consists of smooth immersions only,
that one can construct an appropriate smooth family of
smooth diffeomorphisms $\Phi_t: \Sigma \to \Sigma$, for 
$t \in [0,T]$, satisfying $\Phi_0=\textnormal{Id}_{\Sigma}$, such that the composition $F_t\circ \Phi_t$ solves the Willmore flow equation (\ref{Willmore.flow}) on $[0,T] \times \Sigma$, which proves also the second direction of the assertion.   
\qed
\noindent

\section{Proof of Theorem \ref{convergence}}
\label{Proof.main.theorem}
{\bf Part I}
Short-time existence and uniqueness of the
Willmore flow (\ref{Willmore.flow}) is classical. 
In particular, for any smooth closed regular 
curve $\gamma_0:\sphere^1 \longrightarrow \sphere^2$ 
and for any smooth parametrization 
$F_0:\Sigma \longrightarrow \sphere^3$ of the 
corresponding smooth Hopf-torus 
$\pi^{-1}(\textnormal{trace}(\gamma_0))$ there is 
a unique smooth short-time solution 
$\PP(\,\cdot\,,0,F_0)$ of evolution equation (\ref{Willmore.flow}) on $[0,\varepsilon] \times \Sigma$ 
with $\PP(0,0,F_0) = F_0$ on $\Sigma$. 
Hence, there has to exist a unique maximal time
$T_{\textnormal{max}} \in (\varepsilon,\infty]$, such that this short-time solution of equation (\ref{Willmore.flow}) 
extends from $[0,\varepsilon] \times \Sigma$ to a unique 
smooth solution of the same evolution equation on $[0,T_{\textnormal{max}}) \times \Sigma$. Moreover, we know from Theorem 1.1 in \cite{Dall.Acqua.Pozzi.2018}, that there is a unique global, smooth solution $\{P(t,0,\gamma_0)\}_{t \geq 0}$ 
to equation (\ref{elastic.energy.flow}), satisfying the 
Cauchy problem
\begin{equation} \label{short.time.curves}
\partial_t \tilde \gamma_t =
- \,\nabla_{L^2} \Wil(\tilde \gamma_t), \quad
\textnormal{with} \,\,\, \tilde \gamma_0 = \gamma_0 \quad \textnormal{on} \,\,\, \sphere^1.
\end{equation}
Now, on account of the requirement that $F_0$ is a simple
map from $\Sigma$ onto $\pi^{-1}(\textnormal{trace}(\gamma_0))$
and recalling the uniqueness of the Willmore flow, we obtain from Proposition \ref{correspond.flows}, that the  
short-time solution $\{\PP(t,0,F_0)\}_{t\in [0,\varepsilon]}$ of the Willmore flow (\ref{Willmore.flow}), starting in $F_0$, consists of simple smooth maps from $\Sigma$ onto $\pi^{-1}(\textnormal{trace}(P(t,0,\gamma_0)))$, for every $t \in [0,\varepsilon]$. Hence, the unique short-time flow line
$\{\PP(t,0,F_0)\}_{t\in [0,\varepsilon]}$ of the Willmore flow corresponds to the unique solution 
$\{P(\,\cdot \,,0,\gamma_0)\}$
of the Cauchy problem (\ref{short.time.curves}), restricted 
to $[0, \varepsilon] \times \sphere^1$,
via the Hopf-fibration $\pi$, in the precise sense, that
there is a family of horizontal smooth lifts
$\eta_{\PP(t,0,F_0)}:\sphere^1\cap B_{\epsilon_t}(s^*) \to \Sigma$ of $P(t,0,\gamma_0)$ w.r.t. $\pi \circ \PP(t,0,F_0)$, satisfying 
\begin{equation}  \label{get.it}
P(t,0,\gamma_0) = \pi \circ \PP(t,0,F_0) 
\circ \eta_{\PP(t,0,F_0)}
\end{equation}
on $\bigcup_{t \in [0,\varepsilon]} \{t\} \times
(\sphere^1\cap B_{\epsilon_t}(s^*))$, for every fixed
$s^* \in \sphere^1$; see Lemma \ref{closed.lifts}
for the existence of such horizontal smooth lifts.
Now, since we know from Theorem 1.1 in \cite{Dall.Acqua.Pozzi.2018}, that the maximal solution 
$\{P(t,0,\gamma_0)\}_{t\geq 0}$ of
the initial value problem (\ref{short.time.curves})
is global and smooth, we immediately derive from Proposition
\ref{correspond.flows}, from formula (\ref{get.it}), and 
from the uniqueness of flow lines of the Willmore flow, 
that there has to hold ``$T_{\textnormal{max}} = \infty$'',
that formula (\ref{get.it}) holds
on $\bigcup_{t \in [0,\infty)} \{t\} \times 
(\sphere^1\cap B_{\epsilon_t}(s^*))$ and that  
the resulting global smooth solution 
$\{\PP(t,0,F_0)\}_{t \geq 0}$ 
of the Willmore flow equation is unique.
Moreover, using formula (\ref{second.fundam.trace.free}) and
the fact that every immersion $\PP(t,0,F_0)$ has to be 
a smooth and simple parametrization of the 
Hopf-torus $\pi^{-1}(\textnormal{trace}(P(t,0,\gamma_0)))$ 
- because identity (\ref{get.it}) continues to hold for 
any $t \geq 0$ - we see that the flow line $\{\PP(t,0,F_0)\}_{t \geq 0}$ of the 
Willmore flow consists of umbilic free immersions only, 
just as asserted in Part I of this theorem. \\\\
\noindent
{\bf Part II} Now we start to prove ``subconvergence'' of the 
Willmore flow, as $t \to \infty$.
First of all, given the global flow line
$\{F_t\}_{t\geq 0}:=\{\PP(t,0,F_0)\}_{t\geq 0}$ 
of the Willmore flow (\ref{Willmore.flow}) starting in some prescribed immersion $F_0$, which maps the torus $\Sigma$ simply onto some prescribed Hopf-torus
$\pi^{-1}(\textnormal{trace}(\gamma_0))$, we infer from 
the proof of Part I of this theorem, that there is a unique and global flow line $\{\gamma_t\}_{t \geq 0}:=
\{P(t,0,\gamma_0)\}_{t \geq 0}$ 
of evolution equation (\ref{elastic.energy.flow}), 
starting in the prescribed curve $\gamma_0$, which corresponds
to the global flow line $\{F_t\}$ of equation (\ref{Willmore.flow}) in the precise sense of identity (\ref{get.it}), holding for $t \geq 0$. 
Now we note that
\begin{eqnarray} \label{Pseudo.Gradient.flow}
\frac{d}{dt} \Wil(\gamma_t)
= - \int_{\sphere^1}  \,\Big{|} 2 \,\Big{(}\nabla^{\perp}_{\gamma_t'/|\gamma_t'|}\Big{)}^2
(\vec{\kappa}_{\gamma_t})
+ | \vec{\kappa}_{\gamma_t}|^2\, \vec\kappa_{\gamma_t} +
\,\vec{\kappa}_{\gamma_t} \Big{|}^2 \, d\mu_{\gamma_t}  \leq 0,    \nonumber
\end{eqnarray}
for any $t\geq 0$, i.e. that the function $[t \mapsto \Wil(P(t,0,\gamma_0))]$ is not increasing on $[0,\infty)$. In particular, there exists some constant $K(\gamma_0)>0$, such that both
\begin{eqnarray}  \label{Length.total.curva}
\textnormal{length}(P(t,0,\gamma_0)) \leq K(\gamma_0)  \,\,\, \textnormal{and} \,\,\,         
\int_{\sphere^1} |\vec{\kappa}_{P(t,0,\gamma_0)}|^2  \, d\mu_{P(t,0,\gamma_0)} \leq K(\gamma_0)
\end{eqnarray}
hold for every $t \in [0,\infty)$. Now, applying the 
elementary inequality
$$
\Big{(} \int_{\sphere^1} |\vec{\kappa}_{\gamma}| \,
d\mu_{\gamma}\Big{)}^2 \geq  4\pi^2 - \textnormal{length}(\gamma)^2
$$
holding for every closed smooth regular path
$\gamma:\sphere^1 \longrightarrow \sphere^2$, see \cite{Teufel},
Proposition 1, one can easily derive the lower bound
\begin{equation} \label{lower.bound.length}
\textnormal{length}(\gamma) \geq
\min\Big{\{} \pi, \frac{3\pi^2}{\Wil(\gamma)} \Big{\}},
\end{equation}	
see Lemma 2.9 in \cite{Dall.Acqua.Pozzi.2018}, for any
closed smooth regular path $\gamma:\sphere^1 \longrightarrow \sphere^2$. Combining estimate (\ref{lower.bound.length}) 
with statement (\ref{Pseudo.Gradient.flow}), we see that
the lengths of the curves $P(t,0,\gamma_0)$ can be uniformly 
bounded from below, as well:
\begin{equation} \label{lower.bound.8.3pi}
	\textnormal{length}(\textnormal{trace}(P(t,0,\gamma_0))) \geq \min\Big{\{} \pi, \frac{3 \pi^2}{\Wil(\gamma_0)} \Big{\}}=:\emph{l}(\gamma_0), \quad \forall \, t \in [0,\infty),
\end{equation}
which particularly rules out extinction of flow lines of the
flow (\ref{elastic.energy.flow}) at any time $t\in [0,\infty)$.
Now, using both (\ref{Length.total.curva}) and (\ref{lower.bound.8.3pi}), one can argue as in Section 4.2 in \cite{Dall.Acqua.Pozzi.2018}, using also Lemma 2.6 in \cite{Dall.Acqua.Pozzi.2018}, or also as in Steps 1-6 of the proof of Theorem 1.1 of \cite{Dall.Acqua.Pozzi.2014}, that there hold uniform curvature estimates for the flowing  
curves $\gamma_t = P(t,0,\gamma_0)$, namely:
\begin{equation}  \label{curvature.bound.infty}
	\parallel \Big{(} \nabla_{\frac{\gamma_t'}{|\gamma_t'|}}\Big{)}^m
	(\vec{\kappa}_{\gamma_t}) \parallel_{L^{\infty}(\sphere^1)} \leq C(\gamma_0,m), \quad \textnormal{for} \quad t \in [0,\infty),
\end{equation}
for some constant $C(\gamma,m)>0$, for each $m\in \nat_0$. From this fact and from the
compactness of $\sphere^2$ one can easily conclude, that
for every sequence $t_k \nearrow \infty$ there exists some subsequence $\{t_{k_j}\}$ and some regular smooth closed 
curve $\gamma_{\infty}:\sphere^1 \to \sphere^2$ such that
\underline{reparametrizations} $\tilde{\gamma}_{t_{k_j}}$ of $\gamma_{t_{k_j}}$ to \underline{arc-length} converge:
\begin{equation} \label{convergence.gamma.t}
\tilde{\gamma}_{t_{k_j}} \longrightarrow \gamma_{\infty} \not =\{\textnormal{point}\}
\quad \textnormal{in} \quad C^m(\sphere^1,\rel^3), 
\quad \forall \,m\in \nat_0,
\end{equation}
as $j\to \infty$. For ease of exposition, we rename the
subsequence $\{t_{k_j}\}$ into $\{t_j\}$ again.
We note here that $\textnormal{trace}(\gamma_{\infty})$ cannot be degenerated, i.e. a point on $\sphere^2$, because of the lower bound (\ref{lower.bound.8.3pi}) for the lengths of the traces of curves $\gamma_t=P(t,0,\gamma)$, for $t\in [0,\infty)$. 
Moreover, we can derive from line (\ref{convergence.gamma.t}) that $\A(F_{t_j}):=\int_{\Sigma} d\mu_{F_{t_j}} \leq \textnormal{const}$. $\forall \, j \in \nat$, and from formula (\ref{all.derivatives}) and estimate (\ref{curvature.bound.infty}) we obtain the uniform curvature bounds for the immersions $F_t=\PP(t,0,F_0)$:
\begin{equation} \label{estimate.2nd.fund.form}
	\parallel  (\nabla^{\perp_{F_t,\sphere^3}})^m
	(A_{F_{t},\sphere^3})
	\parallel_{L^{\infty}(\Sigma)} \leq \textnormal{Const}(F_0,m),
	\qquad \textnormal{for} \quad t \in [0,\infty),
\end{equation}
for each $m\in \nat_0$, where we indicated that 
we can only derive from this argument estimates for the 
covariant derivatives $(\nabla^{\perp_{F_t,\sphere^3}})^m$ of $A_{F_{t}}\equiv A_{F_{t},\sphere^3}$ in 
the normal bundle of $F_t$ within $T\sphere^3$. 
Now, arguing as in the proof of Lemma 2.1 in \cite{Ndiaye.Schaetzle.2014}, one can prove by induction, 
that for any smooth section $\Phi$ of the normal bundle
of some smooth immersion $F:\Sigma \longrightarrow \sphere^3 \subset \rel^4$ within $T\sphere^3$ there holds: 
\begin{equation}  \label{compare.covariant.derivatives} 
	(\nabla^{\perp_{F,\rel^4}})_{i_m}\circ \cdots 
	\circ (\nabla^{\perp_{F,\rel^4}})_{i_1}(\Phi)
	= (\nabla^{\perp_{F,\sphere^3}})_{i_m} \circ \cdots 
	\circ (\nabla^{\perp_{F,\sphere^3}})_{i_1}(\Phi)
\end{equation}     
for any order $m\in \nat$ and any $i_j\in \{1,2\}$. 
Moreover, we recall from formula 
(\ref{compare.second.fund.form}), that \, 
$A_{F,\rel^4} = A_{F,\sphere^3} - F\, g_F$. 
Since the first covariant derivative 
$\nabla^{F,\rel^4}(F)$ of the immersion $F$ is a 
section of the tangent bundle of $F$, its projection 
$\nabla^{\perp_{F,\rel^4}}(F)$ into the normal bundle of 
$F$ within $\rel^4$ vanishes identically on $\Sigma$, 
see also here the proof of Lemma 2.1 in \cite{Ndiaye.Schaetzle.2014}. 
Hence, every further derivative of $\nabla^{\perp_{F,\rel^4}}(F)$
vanishes as well. Combining this insight with the fact that 
$\nabla^{F,\rel^4}(g_F)=0$ - by definition of 
$\nabla^{F,\rel^4}$ - and with formulae (\ref{compare.second.fund.form}), (\ref{estimate.2nd.fund.form}) and (\ref{compare.covariant.derivatives}) - here applied to $\Phi = A_{F,\sphere^3}$ - we arrive at stronger uniform curvature bounds for $F_t=\PP(t,0,F_0)$:
\begin{eqnarray} \label{estimate.2nd.fund.form.2}
	\parallel  (\nabla^{\perp_{F_t,\rel^4}})^m
	(A_{F_t,\rel^4})	\parallel_{L^{\infty}(\Sigma)} 
	=\parallel  (\nabla^{\perp_{F_t,\rel^4}})^m
	(A_{F_t,\sphere^3})	\parallel_{L^{\infty}(\Sigma)}    \nonumber\\  
	= \parallel (\nabla^{\perp_{F_t,\sphere^3}})^m(A_{F_t,\sphere^3})\parallel_{L^{\infty}(\Sigma)} 	
	\leq \textnormal{Const}(F_0,m), \quad \textnormal{for} 
	\,\,\, t \in [0,\infty),
\end{eqnarray}
for each $m\in \nat$. Finally we have again by estimate 
(\ref{estimate.2nd.fund.form}): 
\begin{eqnarray}  \label{estimate.2nd.fund.form.0}
	|A_{F_t,\rel^4}(x)|^2 	
	=|A_{F_t,\sphere^3}(x) - F_t(x)\, g_{F_t}(x)|^2   \nonumber   \\
	\equiv g_{F_t}^{ij}(x) g_{F_t}^{lk}(x)  
	\langle (A_{F_t,\sphere^3} - F_t\, g_{F_t})_{il}(x),
	(A_{F_t,\sphere^3} - F_t\, g_{F_t})_{jk}(x) \rangle_{\rel^4}  \nonumber \\
	= |A_{F_t,\sphere^3}|^2(x) + 2  \leq 
	\textnormal{Const}(F_0,0)^2 + 2, \quad \forall \,x\in \Sigma,
\end{eqnarray}
for every $t \in [0,\infty)$. Hence, we may apply here Proposition \ref{Breuning} to the family of immersions 
$F_t=\PP(t,0,F_0):\Sigma \longrightarrow \sphere^3$ 
and obtain the existence of a subsequence 
$\{F_{t_{j_k}}\}$ and of $C^{\infty}$-diffeomorphisms
$\varphi_{k}:\Sigma \stackrel{\cong}\longrightarrow \Sigma$, 
such that
\begin{equation}  \label{u.j.to.infty}
F_{t_{j_k}} \circ \varphi_{k} 
\longrightarrow \hat F    \quad
\textnormal{converge in} \quad  C^{m}(\Sigma,\rel^4),
\quad \textnormal{as}\,\, k \to \infty,
\end{equation}
for each $m \in \nat_0$, in the sense of Definition B.9 in Appendix B of \cite{Dall.Acqua.Schaetzle.Mueller.2020},
where $\hat F$ is again a $C^{\infty}$-smooth immersion of $\Sigma$ into $\sphere^3$, satisfying the bounds (\ref{estimate.2nd.fund.form}), (\ref{estimate.2nd.fund.form.2}) 
and (\ref{estimate.2nd.fund.form.0}). Combining convergence 
(\ref{u.j.to.infty}) with convergence (\ref{convergence.gamma.t}) 
and with identity (\ref{get.it}), we see that there holds:
\begin{eqnarray} \label{rechts.links}
	\pi^{-1}(\textnormal{trace}(\gamma_{\infty}))  \longleftarrow
	\pi^{-1}(\textnormal{trace}(\gamma_{t_{j_k}})) 
	= \textnormal{image}(F_{t_{j_k}})        \nonumber  \\
	=\textnormal{image}(F_{t_{j_k}}
	\circ \varphi_{k}) \longrightarrow 
	\textnormal{image}(\hat F),\quad
\end{eqnarray}
in Hausdorff-distance, as $k \to \infty$. 
Equation (\ref{rechts.links}) particularly shows 
that there holds $\textnormal{image}(\hat F)
= \pi^{-1}(\textnormal{trace}(\gamma_{\infty}))$,  
and since we know already that $\gamma_{\infty}$ is a non-degenerate,
regular closed curve in $\sphere^2$, this means that the limit immersion $\hat F$ maps the torus $\Sigma$ onto a non-degenerate Hopf-torus in $\sphere^3$, namely $\pi^{-1}(\textnormal{trace}(\gamma_{\infty}))$, just as asserted 
in line (\ref{convergence.1.3}). 
Obviously, the immersion $\hat F$ is a simple map
from $\Sigma$ onto $\pi^{-1}(\textnormal{trace}(\gamma_{\infty}))$
in the sense of Definition \ref{General.lift},
because it is the uniform limit of the simple parametrizations
$F_{t_{j_k}} \circ \varphi_{k}$, taking the proved convergence
(\ref{convergence.1.3}) respectively 
(\ref{u.j.to.infty}) into account.  \\ 
Finally we prove that $\hat F$ is actually ``Willmore''. 
To this end, we firstly use formula (\ref{Pseudo.Gradient.flow}),
in order to integrate the derivative of the function
$t\mapsto \Wil(\gamma_t)$ over the interval $[0,T)$,
for any fixed $T>0$:
\begin{equation}  \label{finite.L1.on.infty}
\limsup_{T \to \infty}  \int_0^T \int_{\sphere^1}
|\nabla_{L^2} \Wil(\gamma_t)|^2  \,d\mu_{\gamma_t} dt
=\limsup_{T \to \infty} (\Wil(\gamma_0) - \Wil(\gamma_T))
\leq \Wil(\gamma_0) <\infty.
\end{equation}
Therefore, the limit
$\int_0^{\infty} \int_{\sphere^1} \,
|\nabla_{L^2} \Wil(\gamma_t)|^2  \,d\mu_{\gamma_t} dt$ 
exists and is finite. We can now apply Lemma 3.1 in \cite{Dall.Acqua.Pozzi.2014}, respectively Lemma 2.5 in \cite{Dall.Acqua.Pozzi.2018}, in 
order to derive from the uniform curvature-bounds in (\ref{curvature.bound.infty}), as well:
\begin{equation} \label{curvature.bound.derivative}
\parallel \nabla_{t} \circ
\Big{(}\nabla_{\frac{\gamma_t'}{|\gamma_t'|}}\Big{)}^m
(\vec{\kappa}_{\gamma_t}) \parallel_{L^{\infty}(\sphere^1)}
\leq  C^*(\gamma_0,m),  
\quad \textnormal{for} \,\,t\in [0,\infty),
\end{equation}
along the flow (\ref{elastic.energy.flow}), for some positive 
constant $C^*(\gamma_0,m)$, for each $m\in \nat_0$. Combining estimates (\ref{curvature.bound.derivative})
again with estimates (\ref{Length.total.curva}) and 
(\ref{curvature.bound.infty}), we infer via the usual 
chain- and product/quotient-rule:
$$
\Big{|} \frac{d}{dt} \Big{(} \int_{\sphere^1}
\,|\nabla_{L^2} \Wil(\gamma_t)|^2  \,d\mu_{\gamma_t} \Big{)} \Big{|} \leq C(\gamma_0), \quad \textnormal{for} \,\,t\in [0,\infty),
$$
for some constant $C(F_0)>0$. Together with (\ref{finite.L1.on.infty}) we achieve full convergence:
$$
\int_{\sphere^1}\,
|\nabla_{L^2} \Wil(\gamma_{t})|^2  \,d\mu_{\gamma_{t}}
\longrightarrow 0 \quad \textnormal{as} \,\,\,t \to \infty.
$$
In particular, there has to be some subsequence $\{t_{j_l}\}$
of the particular sequence $t_j \to \infty$ of
convergence (\ref{convergence.gamma.t}), such that
$$
|\nabla_{L^2} \Wil(\gamma_{t_{j_l}})|^2 \longrightarrow 0  \quad \textnormal{in}\,\,\, \Hnm\textnormal{-almost every point of} \,\,\, \sphere^1
$$
as $l \to \infty$. Hence, together with the smooth convergence 
in (\ref{convergence.gamma.t}) we conclude:
$$
0 \longleftarrow
\nabla_{L^2} \Wil(\tilde{\gamma}_{t_{j_l}})
\longrightarrow  \nabla_{L^2} \Wil(\gamma_{\infty})  
\quad \textnormal{in}\,\, 
\Hnm\textnormal{-almost every point of} \,\, \sphere^1.
$$
Therefore, the limit curve $\gamma_{\infty}$ is a smooth elastic curve, i.e. satisfies: $\nabla_{L^2} \Wil(\gamma_{\infty}) \equiv 0$ on $\sphere^1$. 
Hence, the ``Hopf-Willmore-formula'' (\ref{DHopf.first.Variation}) and formula (\ref{Hopf.unit.normal}) immediately imply, that 
every smooth and simple parametrization $F^*:\Sigma \longrightarrow \sphere^3$ of $\pi^{-1}(\textnormal{trace}(\gamma_{\infty}))$
satisfies $\nabla_{L^2} \Will(F^*) \equiv 0$ on $\Sigma$,
and especially the smooth limit immersion $\hat F$ 
in convergence (\ref{u.j.to.infty}) turns out to be ``Willmore'', which finishes the proof of Part II of the theorem.\\\\
{\bf Part III}
Instead of directly proving the assertion of this part of the
theorem about flow lines of the Willmore flow
(\ref{Willmore.flow}), which start to move with Willmore energy
smaller than or equal to $\frac{8\pi^2}{\sqrt{2}}$,
we will at first prove the full convergence of flow lines
$P(\,\cdot\,,0,\gamma_0)$ of the simpler flow
(\ref{elastic.energy.flow}), moving closed smooth curves
in $C^{\infty}_{\textnormal{reg}}(\sphere^1,\sphere^2)$,
which start in some fixed curve
$\gamma_0 \in C^{\infty}_{\textnormal{reg}}(\sphere^1,\sphere^2)$
whose trace is $(\pi \circ F_0)(\Sigma)$ in 
$\sphere^2$ and whose elastic energy
$\Wil(\gamma_0)$ is smaller than or equal to
$\frac{8\pi}{\sqrt{2}}$, on account of formula (\ref{Will.Wil}). We recall from line (\ref{Pseudo.Gradient.flow}), that the function $[t \mapsto \Wil(P(t,0,\gamma_0))]$ is not increasing on $[0,\infty)$, for every initial curve $\gamma_0 \in C^{\infty}_{\textnormal{reg}}(\sphere^1,\sphere^2)$.
Therefore, our initial condition 
``$\Wil(\gamma_0) \leq \frac{8\pi}{\sqrt{2}}$'' 
implies that $\lim_{t \to \infty} \Wil(P(t,0,\gamma_0))$ 
exists and satisfies:
$\lim_{t \to \infty} \Wil(P(t,0,\gamma_0)) \in
\big{[}2\pi, \frac{8\pi}{\sqrt{2}} \big{]}$.
Moreover, from Part II we recall the existence of some 
sequence $\{\gamma_{t_j}\}:=\{P(t_j,0,\gamma_0)\}$, with 
$t_j \nearrow \infty$, and of some sequence of smooth diffeomorphisms $\psi_j:\sphere^1 \stackrel{\cong}\longrightarrow \sphere^1$,
such that the reparametrizations 
$\gamma_{t_j} \circ \psi_j$ have \underline{constant speed} on 
$\sphere^1$ and converge in every $C^m(\sphere^1,\rel^3)$-norm to some elastic curve $\gamma_{\infty}$, again being parametrized with constant speed. We obtain therefore:
\begin{equation} \label{light}
2\pi \leq \Wil(\gamma_{\infty}) = 
\lim_{j \to \infty} \Wil(\gamma_j \circ \psi_j) 
=\lim_{t \to \infty} \Wil(P(t,0,\gamma_0)) 
\leq \frac{8\pi}{\sqrt{2}}.   
\end{equation}
{\bf Case 1}: $\Wil(\gamma_{\infty}) > 2\pi$.
In this case, formula (\ref{light}) contradicts
Proposition \ref{Energy.gap} below, stating that
there are no critical values of $\Wil$ in the 
interval $\big{(}2 \pi,\frac{8\pi}{\sqrt{2}} \big{]}$. \\
{\bf Case 2}: $\Wil(\gamma_{\infty}) = 2\pi$, i.e. in this second case the limit curve $\gamma_{\infty}$ is a smooth and 
regular parametrization of some great circle in $\sphere^2$ 
of constant speed $1$.     \\
We are going to prove the assertion of Part III of this theorem in this remaining second case in the following three steps.\\
{\bf Step 1}: As on p. 2187 in \cite{Dall.Acqua.Spener.2016} 
we firstly assume, that the function
$[t \mapsto \Wil(P(t,0,\gamma_0))]$ 
was not strictly monotonically decreasing 
for $t \in [0,\infty)$. In this case, there was 
some finite time $t^*\geq 0$, such that
$\partial_t(\Wil(P(t,0,\gamma_0)))\lfloor_{t=t^*}=0$.
Then we would have: 
$0=- \int_{\sphere^1}
|\nabla_{L^2}\Wil(P(t,0,\gamma_0))|^2 \,d\mu_{P(t,0,\gamma_0)}\lfloor_{t=t^*}$ 
on account of equation (\ref{elastic.energy.flow}),
implying that the path $P(t^*,0,\gamma_0)$ would 
parametrize an elastic curve with elastic energy $\Wil(P(t^*,0,\gamma_0)) \in 
\big{[}2\pi,\frac{8\pi}{\sqrt{2}} \big{]}$,
hence exactly with elastic energy $2 \pi$ on account of 
Proposition \ref{Energy.gap}. 
Again using the weak monotonicity of the function 
$[t \mapsto \Wil(P(t,0,\gamma_0))]$,
this would imply that $\Wil(P(t,0,\gamma_0))=2 \pi$ 
for every $t \geq t^*$, and therefore again on account of 
equation (\ref{elastic.energy.flow}):
$0=\partial_t(\Wil(P(t,0,\gamma_0)))
=-\int_{\sphere^1} |\nabla_{L^2}\Wil(P(t,0,\gamma_0))|^2 \,d\mu_{P(t,0,\gamma_0)}$, for $t \geq t^*$. 
Combining this again with evolution equation (\ref{elastic.energy.flow}), we see that the flow line 
$\{P(t,0,\gamma_0)\}_{t\geq 0}$ would not move at all 
for $t\geq t^*$ and thus would have got stuck in a 
smooth parametrization $\gamma^{*}$ of some great circle 
in $\sphere^2$. In particular we would obtain here: 
\begin{equation} \label{flat.convergence}
P(t,0,\gamma_0) \longrightarrow \gamma^{*}  \,\,\,
\textnormal{in} \,\, C^m(\sphere^1,\rel^3),
\,\, \textnormal{as} \,\, t \to \infty,
\,\, \textnormal{for every} \,\, m\in \nat_0.
\end{equation} 	    
{\bf Step 2}: We shall suppose throughout in steps 2 and 3 of 
the proof of this third part of the theorem,
that the function $[t \mapsto \Wil(P(t,0,\gamma_0))]$ 
is strictly monotonically decreasing on $[0,\infty)$.
Combining this assumption with estimates (\ref{Length.total.curva}), (\ref{lower.bound.8.3pi}) and (\ref{curvature.bound.infty}) and with our version of the Lojasiewicz-Simon-inequality, Proposition \ref{Lojasiewicz}, 
we will be able, to exactly copy the reasoning of 
Section 4.4 in \cite{Rupp.Spener.2020}, 
in order to arrive at the full smooth convergence - below in (\ref{convergence.p.reparam}) - of the  
\underline{constant speed}-reparametrization 
$\tilde P(\, \cdot \,,0,\gamma_0)$ of the global smooth 
flow line $P(\,\cdot\,,0,\gamma_0)$ of 
evolution equation (\ref{elastic.energy.flow}), 
instead of working out the much more general and thus more complicated technique of Lemma 4.1 in Section 4 of \cite{Chill_Schatz}, respectively of Theorem 1.2 in Section 5 of \cite{Dall.Acqua.Spener.2016}. The gist of Rupp's and Spener's argument in \cite{Rupp.Spener.2020} consists of the pragmatic, new idea, to simply work with the entire reparametrized flow line $\tilde P(\, \cdot \,,0,\gamma_0)$ - and not only with a convergent subsequence of it yielding some smooth and stationary limit curve, as we did above following e.g. Section 5 of \cite{Dall.Acqua.Spener.2016} - and to realize that this reparametrization can and should be written down explicitly, allowing for further precise investigation.  
In order to translate Definition 4.9 of \cite{Rupp.Spener.2020} to our situation, we should identify $\sphere^1$ with $[-\pi,\pi]/(-\pi\sim \pi)$ via the map 
$s\mapsto x=\frac{\log(s)}{i}$, where ``$\log$'' denotes 
the principal branch of the natural logarithm, and we should 
rather consider the curves 
$x\mapsto f_t(x):=P(t,0,\gamma_0)(\exp(ix))$, for $x\in[-\pi,\pi]$, instead of $s\mapsto P(t,0,\gamma_0)(s)$, 
for $s\in \sphere^1$, for each fixed time $t\in [0,\infty)$. 
Working with this slightly changed notation, we are able to 
define exactly as in Definition 4.9 of \cite{Rupp.Spener.2020}: 
\begin{equation}  \label{reparametrized.flow} 
\tilde f_t(x):= f_t(\psi_t(x)), 
\,\,\, \textnormal{for} \,\,(x,t)\in [-\pi,\pi]\times [0,\infty),
\end{equation}	  
where $\psi_t:[-\pi,\pi] \stackrel{\cong}\longrightarrow [-\pi,\pi]$ is the inverse of the smooth diffeomorphism
$$ 
\varphi_t(y):=\frac{2\pi}{\textnormal{length}(P(t,0,\gamma_0))}
\int_{-\pi}^y |\partial_xf_t(z)| \,dz, \,\,\, \textnormal{for} 
\,\, y \in [-\pi,\pi],
$$        
defined here for each fixed $t\in [0,\infty)$.  
The reparametrization in (\ref{reparametrized.flow})  
automatically yields the smooth reparametrization 
$\tilde P(t,0,\gamma_0)(s)
:=\tilde f_t\big{(}\frac{\log(s)}{i}\big{)}$, for  
$s \in \sphere^1$, of the original flow line 
$P(\,\cdot\,,0,\gamma_0)$. 
Now, Lemma 4.10 of \cite{Rupp.Spener.2020} and estimate (\ref{Length.total.curva}) yield the important comparison: 
\begin{equation}  \label{compare} 
\parallel \tilde f_t \parallel_{L^2([-\pi,\pi],\Lno)} 
\leq \sqrt{4\pi}\, \sqrt{\frac{1}{\emph{l}(\gamma_0)} 
+ 4 \,\Wil(f_0)} \,
\parallel f_t \parallel_{L^2([-\pi,\pi],\mu_{f_t})},	 
\end{equation}  
for each fixed $t \in [0,\infty)$, whose proof can be adopted 
here from \cite{Rupp.Spener.2020} without any alteration, taking 
also formula (2.14) in \cite{Dall.Acqua.Pozzi.2018} into account. For the convenience of the reader we shall now sketch the 
proof of Theorem 1.2 in \cite{Rupp.Spener.2020}, 
in order to arrive at the desired convergence (\ref{convergence.p.reparam}) below: 
First of all, we apply Proposition \ref{Lojasiewicz} 
to some small $C^4$-neighbourhood of our smooth elastic 
limit curve $\gamma_{\infty}$, stating that there are 
constants $\theta \in (0,\frac{1}{2}]$,
$c \geq 0$ and $\sigma>0$, only depending on $\gamma_{\infty}$, 
such that for every curve 
$\gamma\in C^{4}_{\textnormal{reg}}(\sphere^1,\sphere^2)$ 
satisfying $\parallel \gamma - \gamma_{\infty} \parallel_{C^{4}(\sphere^1,\rel^3)}\leq \sigma$ 
there holds:
\begin{equation} \label{Lojasi}
|\Wil(\gamma)-\Wil(\gamma_{\infty})|^{1-\theta} 
\leq c \,\,\Big{(} \int_{\sphere^1} 
|\nabla_{L^2} \Wil(\gamma)|^2 \, d\mu_{\gamma} \Big{)}^{1/2}.
\end{equation}
Moreover, we recall here from the beginning of the 
proof of this third part of the theorem, that there is an 
increasing sequence of times $\{t_j\}$, such that our reparametrized curves $\tilde P(t_j,0,\gamma_0)$ converge 
in every $C^m$-norm to $\gamma_{\infty}$, as $j\to \infty$. 
We can therefore choose $j_0$ that large, such that the suprema  
$$ 
s_j:=\sup \{s\geq t_j \,|\,\parallel \tilde P(t,0,\gamma_0)-
\gamma_{\infty}\parallel_{C^{4}(\sphere^1,\rel^3)}
<\sigma  \,\, \textnormal{for every} \,\, t \in [t_j,s]\,\}
$$
are well-defined and satisfy $s_j>t_j$ for each $j>j_0$,
where $\sigma$ denotes again the constant from inequality (\ref{Lojasi}). We finally recall, that here the function 
$[t \mapsto \Wil(\tilde f_t)=\Wil(P(t,0,\gamma_0))]$ 
is strictly monotonically decreasing and converges to  
$\Wil(\gamma_{\infty})$ as $t\to \infty$, 
on account of statement (\ref{light}). 
We can therefore also introduce the smooth, strictly 
monotonically decreasing and positive function
$G(t):=(\Wil(\tilde f_t)- \Wil(\gamma_{\infty}))^{\theta}$, for $t \in [0,\infty)$, with exponent $\theta$ from inequality (\ref{Lojasi}), and we compute by means of the chain rule, evolution equation (\ref{elastic.energy.flow}) and the invariance of 
the elastic energy w.r.t. smooth reparametrization:
\begin{eqnarray*}
	-\frac{d}{dt}G(t)= -\frac{d}{dt}\Big{(}
	\big{(}\Wil(\tilde f_t) - \Wil(f_{\infty})\big{)}^{\theta} \Big{)}       
    =-\frac{d}{dt}\Big{(} \big{(}\Wil(P(t,0,\gamma_0)) - \Wil(\gamma_{\infty})\big{)}^{\theta} \Big{)}  \\
	= \theta \,\big{(}\Wil(P(t,0,\gamma_0)) - 
	\Wil(f_{\infty})\big{)}^{\theta-1} \,
	\int_{\sphere^1}   \,
	|\nabla_{L^2}\Wil(P(t,0,\gamma_0))|^2 \, 
	d\mu_{P(t,0,\gamma_0)}                              \\
	=  \theta \, \big{(}\Wil(\tilde P(t,0,\gamma_0)) - 
	\Wil(\gamma_{\infty})\big{)}^{\theta-1}
	\parallel \nabla_{L^2} \Wil(P(t,0,\gamma_0))  \parallel_{L^2(\mu_{P_t})} \,
	\parallel  \partial_t P(t,0,\gamma_0)  \parallel_{L^2(\mu_{P_t})}                        \\
	=  \theta \, \big{(}\Wil(\tilde P(t,0,\gamma_0)) - 
	\Wil(\gamma_{\infty})\big{)}^{\theta-1}
	\parallel \nabla_{L^2}\Wil(\tilde P(t,0,\gamma_0))  \parallel_{L^2(\mu_{\tilde P_t})} \,
	\parallel  \partial_t P(t,0,\gamma_0)  \parallel_{L^2(\mu_{P_t})}                         \\
	\geq  \frac{\theta}{c} \,
	\parallel  \partial_t P(t,0,\gamma_0)  \parallel_{L^2(\sphere^1,\mu_{P(t,0,\gamma_0)})}, \,\, \textnormal{for} \,\, t \in [t_j,s_j) \,\,\, 
	\textnormal{and each} \,\,j>j_0, 
\end{eqnarray*}
where we could apply the Lojasiewicz-Simon-gradient inequality 
(\ref{Lojasi}) in the last line to the reparametrized curves 
$\tilde P(t,0,\gamma_0)$, being sufficiently close to the great circle $\gamma_{\infty}$ in $C^{4}(\sphere^1,\rel^3)$ 
for $t \in [t_j,s_j)$, by definition of the suprema $s_j$. 
Combining the above inequality 
$-\frac{d}{dt}G(t)\geq  \frac{\theta}{c} \,
\parallel  \partial_t P(t,0,\gamma_0)  \parallel_{L^2(\sphere^1,\mu_{P(t,0,\gamma_0)})}$ with 
inequality (\ref{compare}), we achieve as in formula (4.9) of 
\cite{Rupp.Spener.2020} the decisive estimate: 
\begin{equation}  \label{decisive.1} 
-\frac{d}{dt}G(t)\geq C(\theta,c,\Wil(\gamma_0),\emph{l}(\gamma_0)) \,
\parallel  \partial_t \tilde P(t,0,\gamma_0)  \parallel_{L^2(\sphere^1,\Hnm)}, \,\,
\textnormal{for} \,\, t\in [t_j,s_j),
\end{equation}  
for each $j>j_0$. 
Integration of inequality (\ref{decisive.1}) from $t_j$ 
to any $T\in (t_j,s_j)$ yields: 
\begin{eqnarray}  \label{decisive.2} 
\parallel \tilde P(T,0,\gamma_0) - \tilde P(t_j,0,\gamma_0) 
\parallel_{L^2(\sphere^1,\Hnm)}     \nonumber     \\
\leq \int_{t_j}^T \parallel \partial_t \tilde P(t,0,\gamma_0)  \parallel_{L^2(\sphere^1,\Hnm)}  \, dt 
\leq \frac{1}{C(\theta,c,\Wil(\gamma_0),\emph{l}(\gamma_0))} 
\, G(t_j),
\end{eqnarray} 
whose right-hand side converges to $0$ as $j\to \infty$ 
because of statement (\ref{light}).   
Now, as in the proof of Theorem 1.2 in \cite{Rupp.Spener.2020}
we can combine statement (\ref{decisive.2}) with the general subconvergence of the elastic energy flow 
(\ref{elastic.energy.flow}) from Section 4.2 in \cite{Dall.Acqua.Pozzi.2018}, which we had already used in 
the proof of the second part of this theorem,  
and with the fact that exactly the reparametrized curves  
$\tilde P(t_j,0,\gamma_0)$ converge smoothly to the 
great circle $\gamma_{\infty}$, in order to prove that 
some of the suprema $s_{j}$, say $s_{J}$ for some large $J>j_0$, 
cannot be finite. Hence, by definition of $s_{J}$ we can 
conclude hereby, that $\parallel \tilde P(t,0,\gamma_0) 
- \gamma_{\infty} \parallel_{C^4(\sphere^1)}<\sigma$ 
holds for every $t\geq t_J$, implying that also inequality (\ref{decisive.1}) holds for every $t\geq t_J$. 
Herewith we can finally infer, that the function 
$[t \mapsto \parallel \partial_t \tilde P(t,0,\gamma_0) 
\parallel_{L^2(\sphere^1,\Hnm)}]$ is of class 
$L^1([0,\infty),\rel)$, and thus the reparametrized 
flow line $\{\tilde P(t,0,\gamma_0)\}_{t\geq 0}$ converges 
fully in $L^2(\sphere^1,\Hnm)$ as $t \to \infty$, 
and its $L^2$-limit has to be the great circle 
$\gamma_{\infty}$, which we had started to work with at 
the beginning of the proof of Part III of the theorem. 
Now, combining this full $L^2$-convergence with 
the uniform curvature estimates (\ref{curvature.bound.infty}) - 
which hold for any smooth reparametrization of the original  
flow line $\{P(t,0,\gamma_0)\}_{t\geq 0}$ of flow (\ref{elastic.energy.flow}) - we finally infer that: 
\begin{equation}  \label{convergence.p.reparam}
\tilde P(t,0,\gamma_0) \longrightarrow   
\gamma_{\infty} \quad \textnormal{in} \,\,\, C^m(\sphere^1,\rel^3), \,\, \textnormal{as} \,\, 
t \to \infty, \,\, \textnormal{for every} \,\, 
m \in \nat_0.
\end{equation} 
Hence, just as in conclusion (\ref{flat.convergence}) of Step I 
of the considered second case within the proof of Part III 
of this theorem, we have arrived in (\ref{convergence.p.reparam}) at the full and smooth convergence of our flow line 
$\{P(t,0,\gamma_0)\}_{t\geq 0}$ of flow (\ref{elastic.energy.flow}) - up to the smooth reparametrization 
in (\ref{reparametrized.flow}) - to a smooth parametrization 
$\gamma_{\infty}$ of some great circle in $\sphere^2$. \\
{\bf Step 3}:
On account of the proof of Part I of Theorem \ref{convergence} we know, that the unique flow line
$\{\PP(\,\cdot\,,0,F_0)\}$ of the Willmore flow
(\ref{Moebius.flow}) corresponds to the global flow line
$\{P(t,0,\gamma_0)\}_{t\geq 0}$ of flow
(\ref{elastic.energy.flow}) via the Hopf-fibration in the
precise sense of formula (\ref{get.it}). 
Hence, every immersion $\PP(t,0,F_0):\Sigma \longrightarrow \sphere^3$ yields a smooth and 
simple map from the torus $\Sigma$ onto the Hopf-torus
$\pi^{-1}(\textnormal{trace}(P(t,0,\gamma_0)))$, for 
$t \geq 0$. Moreover, combining the convergence in line (\ref{convergence.p.reparam}) with the strict monotonicity 
of the function $[t \mapsto \Wil(P(t,0,\gamma_0))]$ and with formula (\ref{Will.Wil}), we can infer from the Li-Yau inequality \cite{Li.Yau}, that the immersions $\PP(t,0,F_0)$ are 
smooth diffeomorphisms between $\Sigma$ and their images $\pi^{-1}(\textnormal{trace}(P(t,0,\gamma_0)))$, for sufficiently large $t \geq t_0>>1$, showing in 
particular that the Hopf-tori $\pi^{-1}(\textnormal{trace}(P(t,0,\gamma_0)))$ are smooth 
compact manifolds of genus one, for $t \geq t_0>>1$.
Furthermore, since we know that $\Wil(\gamma_{\infty}) 
= 2 \pi$ and thus $\Will(\pi^{-1}(\textnormal{trace}(\gamma_{\infty}))) 
= 2 \,\pi^2$ by formula (\ref{Will.Wil}), we can conclude 
from Theorem A in \cite{Marques.Neves},
that the limit Hopf-torus
$\pi^{-1}(\textnormal{trace}(\gamma_{\infty}))$
is a conformal image
$M\big{(}\frac{1}{\sqrt 2}(\sphere^1 \times \sphere^1)\big{)}$ of the Clifford torus, for some suitable $M\in \textnormal{M\"ob}(\sphere^3)$. 
Moreover, by Lemma \ref{closed.lifts} there are smooth
horizontal lifts $\eta_t:\sphere^1 \setminus \{y_{t}^*\} \longrightarrow \sphere^3$, for arbitrarily 
chosen points $y^*_t \in \sphere^1$,
of the curves $\tilde P(t,0,\gamma_{0})$ from line 
(\ref{convergence.p.reparam}) w.r.t. the Hopf-fibration $\pi$,
and by Lemma \ref{Hopf.Tori} and Remark \ref{Fundamental.gruppe.Homologie}
the maps $X_t:(\sphere^1 \setminus \{y^*_t\}) \times [0,2\pi] \longrightarrow \sphere^3$ explicitly given by
\begin{equation}   \label{X.t}
X_t(s,\varphi):= e^{i \varphi} \, \eta_t(s), \quad
\textnormal{for} \,\,\,t \geq t_0,
\end{equation}
are isometric parametrizations of the Hopf-tori
$\pi^{-1}(\textnormal{trace}(P(t,0,\gamma_0)))$, 
covering $\pi^{-1}(\textnormal{trace}(P(t,0,\gamma_0)))$
exactly once up to a subset of vanishing $\Hn$-measure, for every fixed $t \geq t_0$. Now, choosing a smooth finite atlas for the torus $\Sigma$ and a subordinate smooth partition of unity, one can use restrictions of the parametrizations $X_t$ in (\ref{X.t}) to appropriate open subsets of $\sphere^1 \times [0,2\pi]$, in order to construct smooth immersions $Y_t:\Sigma \longrightarrow \sphere^3$, which map the torus $\Sigma$ simply onto $\pi^{-1}(\textnormal{trace}(P(t,0,\gamma_0)))$,
for every $t \geq t_0$. In particular, 
by Definition \ref{General.lift} the induced maps 
in singular homology
\begin{equation}  \label{hom.isomorphism}
(Y_t)_{*2}: H_{2}(\Sigma,\ganz) \stackrel{\cong}\longrightarrow
H_2(\pi^{-1}(\textnormal{trace}(P(t,0,\gamma_0))),\ganz)
\end{equation}
are isomorphisms, for every fixed $t \geq t_0$. 
Since we know already that the Hopf-tori $\pi^{-1}(\textnormal{trace}(P(t,0,\gamma_0)))$ and also $\Sigma$ are smooth compact manifolds of genus one, 
and that $Y_t$ are immersions of $\Sigma$ into
$\sphere^3$, for every $t \geq t_0$, we can conclude 
from statement (\ref{hom.isomorphism}) and
Remark \ref{Fundamental.gruppe.Homologie}, that
the maps $Y_t$ are actually smooth diffeomorphisms
between $\Sigma$ and $\pi^{-1}(\textnormal{trace}(P(t,0,\gamma_0)))$,
for every $t \geq t_0$. Now, recalling (\ref{reparametrized.flow}) we set $\tilde \gamma_t:=\tilde P(t,0,\gamma_{0})$, for
$t\geq 0$, we fix some $y^* \in \sphere^1$, we choose some point $s^* \in  \sphere^1 \setminus \{y^*\}$ arbitrarily, and we also recall from the first part of Lemma \ref{closed.lifts}, that for every fixed $t^* \geq t_0$ and
every fixed $q^*\in \pi^{-1}(\tilde \gamma_{t^*}(s^*)) \subset \sphere^3$ there is a unique horizontal, smooth lift $\eta_{t^*}^{(s^*,q^*)}:
\sphere^1 \setminus \{y^*\} \longrightarrow \pi^{-1}(\textnormal{trace}(\tilde \gamma_{t^*}))$ of
$\tilde \gamma_{t^*} \lfloor_{\sphere^1 \setminus \{y^*\}}$ w.r.t. the Hopf-fibration $\pi$, satisfying
$\eta_{t^*}^{(s^*,q^*)}(s^*)=q^*$, and that this
horizontal lift was obtained via the unique smooth
flow generated by the initial value problem (\ref{Initial.value}) below, here with generating
vector field $V_{\gamma}:= V_{\tilde \gamma_{t^*}}$. 
Now, since the flow line $\{\tilde \gamma_t\}_{t\geq t_0}$ 
is a smooth family of smooth closed regular paths, also the corresponding generating vector fields 
$[(q,t)\mapsto V_{\tilde \gamma_t}(q)]$ - to be substituted 
into line (\ref{Initial.value}) below - are smooth sections of $T\sphere^2$, which depend smoothly on the time $t$ as well. Interpreting the time $t$ as an additional real parameter of the generating vector field $V_{\tilde \gamma_t}$, Theorem 1.5.3 in \cite{Kong.2014} guarantees us, that the unique solution 
of initial value problem (\ref{Initial.value}) with smooth right hand side $[(q,t)\mapsto V_{\tilde \gamma_t}(q)]$ depends smoothly on the initial value $q^*$ and also on the additional 
parameter $t \geq t_0$. We can therefore construct a 
smooth family of horizontal lifts $\eta_t \equiv \eta_{t}^{(s^*,q_t)}:\sphere^1 \setminus \{y^*\} \longrightarrow \sphere^3$ of 
$\tilde \gamma_t \lfloor_{\sphere^1 \setminus \{y^*\}}$ 
w.r.t. $\pi$ with $\eta_{t}^{(s^*,q_t)}(s^*)= q_t$ 
in such a way that $q_t \longrightarrow q_{\infty}$, 
as $t\nearrow \infty$, for some point 
$q_{\infty} \in \pi^{-1}(\gamma_{\infty}(s^*))$,         
making here also use of the full convergence 
of $\{\tilde \gamma_t\}$ to the great circle-parametrization $\gamma_{\infty}$ in (\ref{convergence.p.reparam}). By 
successive differentiation - w.r.t. $\tau \in \rel$ -  
of equation (\ref{Initial.value}) for any fixed 
$t\geq t_0$, which reads here:  
$$ 
\frac{d\eta_{t}^{(s^*,q_t)}}{ds}(s) 
= V_{\tilde \gamma_t}(\eta_{t}^{(s^*,q_t)}(s)),  
\,\,\, \textnormal{for every} \,\, s\in \sphere^1 
\setminus \{y^*\},
$$ 
and again using convergence (\ref{convergence.p.reparam}) successively for every $m\in \nat_0$, one  
infers inductively that any sequence $t_j \nearrow \infty$ possesses a certain subsequence $\{t_{j_k}\}$, such that  
\begin{equation} \label{converg.eta.k}
\eta_{t_{j_k}} \equiv \eta_{t_{j_k}}^{(s^*,q_{t_{j_k}})} 
\longrightarrow \eta_{\infty} \,\,\, 
\textnormal{in} \,\,\, C^m_{\textnormal{loc}}
(\sphere^1 \setminus \{y^*\},\rel^4), \,\,\, 
\textnormal{for each} \,\,\, m \in \nat_0,
\end{equation} 
as $k\to \infty$, for some smooth limit function
$\eta_{\infty}:\sphere^1 \setminus \{y^*\}
\longrightarrow \sphere^3$ satisfying $\eta_{\infty}(s^*)=q_{\infty}$.
The limit function $\eta_{\infty}$ has to be a smooth horizontal lift of the smooth limit curve $\gamma_{\infty}\lfloor_{\sphere^1 \setminus \{y^*\}}$
w.r.t. $\pi$, because of: 
$$
\pi \circ \eta_{\infty} \longleftarrow 
\pi \circ \eta_{t_{j_k}} = \tilde \gamma_{t_{j_k}}  
\longrightarrow \gamma_{\infty} \,\,\, \textnormal{in} 
\,\,\, C^m_{\textnormal{loc}}(\sphere^1 \setminus \{y^*\},\rel^3), \,\, \textnormal{as}\,\, k\to \infty,
$$
for each $m\in \nat_0$, using here both (\ref{convergence.p.reparam}) and (\ref{converg.eta.k}). 
Now, the first part of Lemma \ref{closed.lifts} 
guarantees us, that the horizontal lift $\eta_{\infty}$ of 
$\gamma_{\infty}\lfloor_{\sphere^1 \setminus \{y^*\}}$
is uniquely determined by the additional equation $\eta_{\infty}(s^*)=q_{\infty}\in \pi^{-1}(\gamma_{\infty}(s^*))$. Combining this again with convergence (\ref{converg.eta.k}), 
thus the principal of subsequences guarantees us, that 
there holds actually: 
\begin{equation} \label{converg.eta.t}
\eta_t\equiv \eta_{t}^{(s^*,q_{t})} \longrightarrow 
\eta_{\infty} \,\,\,  \textnormal{in} \,\,\,C^m_{\textnormal{loc}}
(\sphere^1 \setminus \{y^*\},\rel^4), \,\,\, 
\textnormal{for each} \,\,\, m \in \nat_0,
\end{equation}  
as $t\to \infty$, where $y^*\in \sphere^1$ had been fixed arbitrarily. Since the parametrizations $Y_t :\Sigma \stackrel{\cong}\longrightarrow \pi^{-1}(\textnormal{trace}(P(t,0,\gamma_0)))
= \pi^{-1}(\textnormal{trace}(\tilde \gamma_t))$
had been defined via a fixed smooth atlas of
$\Sigma$ and by means of restrictions of the immersions
$X_t$ in (\ref{X.t}) to appropriate open subsets of
$\sphere^1 \times [0,2\pi]$, the full smooth convergence 
of $\{\eta_t\}$ in (\ref{converg.eta.t}) implies that:
\begin{equation}   \label{smooth.convergence}
Y_t \longrightarrow  F^* \,\,\, \textnormal{converge in}
\,\,\, C^m(\Sigma,\rel^4), \,\,\, \textnormal{as} \,\,\, 
t \to \infty,
\end{equation}
for each $m \in \nat_0$, in the sense of Definition B.9 in Appendix B of \cite{Dall.Acqua.Schaetzle.Mueller.2020}.
Here, the limit map $F^*$ can be constructed
exactly as the immersions $Y_t$, but using now the 
above horizontal lift $\eta_{\infty}$ of the great 
circle-parametrization $\gamma_{\infty}$ instead of the lifts $\eta_t$ in formula (\ref{X.t}). As the lift 
$\eta_{\infty}$ is horizontal, also the limit map 
$F^*$ turns out to be an immersion of $\Sigma$ into $\sphere^3$. Moreover, we note that the map $F^*$ maps $\Sigma$ onto $\pi^{-1}(\textnormal{trace}(\gamma_{\infty}))=
M\big{(}\frac{1}{\sqrt 2} (\sphere^1 \times \sphere^1)\big{)}$, because convergences  (\ref{convergence.p.reparam}) and 
(\ref{smooth.convergence}) imply: 
\begin{equation}   \label{X.t.1}
M\big{(}\frac{1}{\sqrt 2} 
(\sphere^1 \times \sphere^1)\big{)}
\longleftarrow  \pi^{-1}(\textnormal{trace}(P(t,0,\gamma_0)))
=\textnormal{image}(Y_t) \longrightarrow  \textnormal{image}(F^*)
\end{equation}
in Hausdorff-distance, as $t \to \infty$. Since the parametrizations $Y_t$ are simple maps, we infer from convergence (\ref{smooth.convergence}) and directly from Definition \ref{General.lift}, that the limit map $F^*$ in  (\ref{smooth.convergence}) is actually a simple map from $\Sigma$ onto $M\big{(}\frac{1}{\sqrt 2} (\sphere^1 \times \sphere^1)\big{)}$.
Now, combining this with the fact that both $\Sigma$ and
$M\big{(}\frac{1}{\sqrt 2}(\sphere^1 \times \sphere^1){)}$ are smooth compact manifolds of genus one, and that $F^*$ is an immersion of $\Sigma$ into $\sphere^3$, we can argue - just as above - by means of Remark \ref{Fundamental.gruppe.Homologie},
that $F^*$ is a smooth diffeomorphism between
the smooth tori $\Sigma$ and $M\big{(}\frac{1}{\sqrt 2}(\sphere^1 \times \sphere^1)\big{)}$. Moreover, recalling that both the immersions $\PP(t,0,F_0)$ 
and $Y_t$ are smooth diffeomorphisms between $\Sigma$ and their images $\pi^{-1}(\textnormal{trace}(P(t,0,\gamma_0)))$,
in every $t \geq t_0$ - for $t_0$ sufficiently large -
we conclude that there is a unique diffeomorphism
$\Theta :\Sigma \stackrel{\cong}\longrightarrow \Sigma$ such that $\PP(t_0,0,F_0) = Y_{t_0} \circ \Theta^{-1}$ holds on $\Sigma$, at time $t=t_0$. Now, since both families $\{\PP(t,0,F_0)\}_{t \geq t_0}$ and $\{Y_t \circ \Theta^{-1}\}_{t \geq t_0}$ parametrize the Hopf-tori $\pi^{-1}(\textnormal{trace}(P(t,0,\gamma_0)))$, for
$t\geq t_0$, Proposition \ref{correspond.flows} yields 
a smooth family of smooth diffeomorphisms 
$\Xi_t :\Sigma \longrightarrow \Sigma$ satisfying
$$
\PP(t,0,F_0) = Y_t \circ \Theta^{-1} \circ (\Xi_t)^{-1}, \,\,\, \textnormal{for every} \,\,\,  t \geq t_0.
$$
Hence, we infer from statement
(\ref{smooth.convergence}) the full, smooth convergence  
$$
\PP(t,0,F_0) \circ (\Xi_t \circ \Theta) = Y_t 
\longrightarrow F^* \,\,\, \textnormal{in} \,\,\, C^m(\Sigma,\rel^4), \,\,\,
\textnormal{as} \,\,\, t \to \infty,
$$
for each $m \in \nat_0$, in the sense of Definition B.9 in \cite{Dall.Acqua.Schaetzle.Mueller.2020}, where $F^*$ is 
the desired smooth diffeomorphism between $\Sigma$ 
and ``a Clifford torus'' $M\big{(}\frac{1}{\sqrt 2} (\sphere^1 \times \sphere^1)\big{)}$.
\qed

\section{Appendix}    \label{Appendix}

In this appendix we firstly prove existence
of local, horizontal and smooth lifts 
of some arbitrary smooth, closed path 
$\gamma:\sphere^1 \longrightarrow \sphere^2$ 
w.r.t. fibrations of the type $\pi \circ F$
for ``simple'' parametrizations $F:\Sigma \longrightarrow
\pi^{-1}(\textnormal{trace}(\gamma))\subset \sphere^3$,
in the sense of Definition \ref{General.lift}. 
Since $\gamma$ should be allowed to have self-intersections  
in view of our Theorem \ref{convergence}, we cannot 
blindly apply the general theory of smooth fiber bundles over smooth base manifolds. Instead, we have to construct here such local lifts elementarily, using Lemmata \ref{Hopf.properties} 
and \ref{Hopf.Tori} and the theory of ``ODEs''.     
\begin{lemma} \label{closed.lifts}
	Let $\gamma: \sphere^1 \longrightarrow \sphere^2$ be a smooth, closed and regular path in $\sphere^2$, and let $F:\Sigma \longrightarrow \sphere^3$ be a
	smooth immersion, which maps a smooth compact torus $\Sigma$ simply and smoothly onto the Hopf-torus
	$\pi^{-1}(\textnormal{trace}(\gamma)) \subset \sphere^3$.
	\begin{itemize}
		\item[1)] For every fixed $s^* \in \sphere^1$
		and $q^* \in \pi^{-1}(\gamma(s^*))\subset \sphere^3$
		there is a unique horizontal, smooth lift $\eta^{(s^*,q^*)}:
		\textnormal{dom}(\eta^{(s^*,q^*)}) \longrightarrow \pi^{-1}(\textnormal{trace}(\gamma))$,
		defined on a non-empty, open interval
		$\textnormal{dom}(\eta^{(s^*,q^*)}) \subset \sphere^1$, of $\gamma:\sphere^1 \longrightarrow \sphere^2$ w.r.t. the Hopf-fibration $\pi$, 
		such that $\textnormal{dom}(\eta^{(s^*,q^*)})$
		contains the point $s^*$ and such that $\eta^{(s^*,q^*)}$ attains the value $q^*$ in $s^*$; i.e. $\eta^{(s^*,q^*)}$ is a smooth path in the torus $\pi^{-1}(\textnormal{trace}(\gamma))$, which intersects the fibers of $\pi$ perpendicularly and satisfies:
		\begin{equation} \label{eta}
		(\pi \circ \eta^{(s^*,q^*)})(s) = \gamma(s)  \quad
		\forall \,s \in \textnormal{dom}(\eta^{(s^*,q^*)})  \,\,\, \textnormal{and} \,\,\, 
		\eta^{(s^*,q^*)}(s^*) = q^*,
		\end{equation}
	    and there is only one such function $\eta^{(s^*,q^*)}$ mapping the open interval $\textnormal{dom}(\eta^{(s^*,q^*)}) \subset 
	    \sphere^1$ into $\pi^{-1}(\textnormal{trace}(\gamma))$.
		\item[2)] There is some $\epsilon=\epsilon(F,\gamma)>0$, such that for every fixed $s^* \in \sphere^1$ and every 
		$x^*\in (\pi \circ F)^{-1}(\gamma(s^*))$
		there is a horizontal smooth lift $\eta_F^{(s^*,x^*)}$ of $\gamma\lfloor_{\sphere^1 \cap B_{\epsilon}(s^*)}$ w.r.t. the fibration
		$\pi \circ F :\Sigma \longrightarrow \textnormal{trace}(\gamma) \subset \sphere^2$, attaining the value $x^*$ in $s^*$, i.e. $\eta_F^{(s^*,x^*)}$ is a smooth path in the torus $\Sigma$ which intersects the fibers of $\pi \circ F$ perpendicularly and satisfies:
		\begin{equation} \label{eta_F}
		(\pi \circ F \circ \eta_F^{(s^*,x^*)})(s) = \gamma(s)  \quad \forall \,s \in \sphere^1 \cap B_{\epsilon}(s^*) \,\,\, \textnormal{and} \,\,\,
		\eta_F^{(s^*,x^*)}(s^*)=x^*.
		\end{equation}	
		In particular, for the above $\epsilon=\epsilon(F,\gamma)>0$ the function 
		$\eta_F \mapsto F \circ \eta_F$ maps
		the set $\Lift(\gamma\lfloor_{\sphere^1 
		\cap B_{\epsilon}(s^*)},\pi \circ F)$
		of horizontal smooth lifts of $\gamma\lfloor_{\sphere^1 \cap B_{\epsilon}(s^*)}$ w.r.t. $\pi \circ F$ 
		surjectively onto the set
		$\Lift(\gamma\lfloor_{\sphere^1 \cap B_{\epsilon}(s^*)},\pi)$ of horizontal
		smooth lifts of $\gamma\lfloor_{\sphere^1 \cap B_{\epsilon}(s^*)}$ w.r.t. $\pi$.
	\end{itemize}
\end{lemma}
\noindent
\proof: Without loss of generality, we may assume here, 
that $\gamma$ performs only one loop through its closed trace, 
and that it is defined on $\rel/L\ganz$ with $|\gamma'|_{g_{\sphere^2}} \equiv 2$, 
where $2L$ is the length of the trace of $\gamma$.
We consider the unique smooth vector field $V_{\gamma} \in \Gamma(T(\pi^{-1}(\textnormal{trace}(\gamma))))\subset
\Gamma(T\sphere^3)$, which intersects the fibers of the Hopf-fibration perpendicularly, satisfies 
$|V_{\gamma}|_{g_{\sphere^3}} \equiv 1$  
throughout on $\pi^{-1}(\textnormal{trace}(\gamma))$ and 
\begin{equation}  \label{assumption.2.X}
D\pi_{q^*}(V_{\gamma}(q^*)) = c^* \, \gamma'(0), 
\,\, \textnormal{for some} \,\, q^* \in \pi^{-1}(\gamma(0)) \,\, \textnormal{and some} \,\, c^* > 0.
\end{equation}  
Now, the unique flow $\Psi$ on $\pi^{-1}(\textnormal{trace}(\gamma))$ 
generated by the initial value problem
\begin{equation}  \label{Initial.value}
\frac{dv}{d\tau}(\tau) = V_{\gamma}(v(\tau)) \qquad \textnormal{and} \quad  v(0)=q,
\end{equation}
for any fixed point $q \in \pi^{-1}(\textnormal{trace}(\gamma))$, exists eternally, because the Hopf-torus $\pi^{-1}(\textnormal{trace}(\gamma))$
is a compact subset of $\sphere^3$, which can locally be parametrized by smooth charts, similarly to a smooth compact $2$-manifold without any boundary points. We can 
readily infer from the requirements on 
$V_{\gamma}$, that the flow lines of the resulting eternal flow
$\Psi:\rel \times \pi^{-1}(\textnormal{trace}(\gamma)) 
\longrightarrow \pi^{-1}(\textnormal{trace}(\gamma))$ intersect the fibers of $\pi$ perpendicularly, have constant speed $1$ and are mapped by $\pi$ onto $\textnormal{trace}(\gamma)$. 
Now we choose some arbitrary initial point  
$q_0 \in \pi^{-1}(\textnormal{trace}(\gamma))$, and we 
derive from the properties 
$$
|\Psi_{\tau}(q_0)|_{g_{\sphere^3}}\equiv 1 \quad
\textnormal{and} \quad
|\frac{d}{d\tau}\Psi_{\tau}(q_0)|_{g_{\sphere^3}} = |V_{\gamma}(\Psi_{\tau}(q_0))|_{g_{\sphere^3}} \equiv 1
$$ 
as in Lemma \ref{Hopf.Tori}, that there is some unit vector field $\tau \mapsto u(\tau)=u_{\gamma,q_0}(\tau) \in \textnormal{Span}\{j,k\} \subset \quat$, 
satisfying $V_{\gamma}(\Psi_{\tau}(q_0)) 
= u(\tau) \cdot \Psi_{\tau}(q_0) \,\, \textnormal{in} \,\,\quat$, for every $\tau \in \rel$, 
and we can consequently compute exactly as in formula (\ref{gamma.speed}): 
\begin{eqnarray} \label{D.pi.X}
	\frac{d}{d\tau}(\pi \circ \Psi_{\tau}(q_0)) 
	= (I(V_{\gamma}(\Psi_{\tau}(q_0))) \cdot \Psi_{\tau}(q_0))
	+(I(\Psi_{\tau}(q_0)) \cdot V_{\gamma}(\Psi_{\tau}(q_0)))         \nonumber         \\
	=  2 \, I(\Psi_{\tau}(q_0)) \cdot u(\tau) 
	\cdot \Psi_{\tau}(q_0) 
	\equiv 2 \,I(\Psi_{\tau}(q_0)) \cdot V_{\gamma}(\Psi_{\tau}(q_0)), \,\,\,
\end{eqnarray}
for every $\tau \in \rel$. Now, since $\pi$ maps the Hopf-torus $\pi^{-1}(\textnormal{trace}(\gamma))$
onto the trace of $\gamma$, we know a-priori that the trace of the path $[\tau \mapsto \pi \circ \Psi_{\tau}(q_0)]$
is contained in the trace of the curve $\gamma$. 
Moreover, equation (\ref{D.pi.X}) shows us:
\begin{equation}  \label{D.pi.X.2}
|\frac{d}{d\tau}(\pi \circ \Psi_{\tau}(q_0))| \equiv 2 
\equiv |\gamma'(\tau)|, \quad \forall \,\tau \in \rel, 
\end{equation}
and for every initial point 
$q_0 \in \pi^{-1}(\textnormal{trace}(\gamma))$. Hence, 
taking especially $q_0 \in \pi^{-1}(\gamma(0))$ we have 
$\pi \circ \Psi_{0}(q_0) = \gamma(0)$,
and therefore assumption (\ref{assumption.2.X}) and  
equation (\ref{D.pi.X.2}) prove, that there has to hold:
\begin{equation}  \label{eta.2}
\pi \circ \Psi_{\tau}(q_0) = \gamma(\tau), \quad 
\forall \, \tau \in \rel,
\end{equation}
and for every $q_0 \in \pi^{-1}(\gamma(0))$. 
Moreover, combining identity (\ref{eta.2}) with the group property ``$\Psi_{\tau_1 +\tau_2} = \Psi_{\tau_1} \circ \Psi_{\tau_2}$'' of the flow $\Psi$, it follows that for any initial point $q_0 \in \pi^{-1}(\textnormal{trace}(\gamma))$ - 
not only for $q_0 \in \pi^{-1}(\gamma(0))$ - the path 
$[\tau \mapsto \pi \circ \Psi_{\tau}(q_0)]$ 
parametrizes the trace of $\gamma$ exactly once, as $\tau$ increases from $0$ to $L$, which means that any such flow line 
$\{\Psi_{\tau}(q_0)\}_{\tau \in [0,L]}$ intersects each fiber of $\pi$ over the trace of $\gamma$ exactly once. 
Now, for any $s^* \in (0,L)$ and 
$q^* \in \pi^{-1}(\gamma(s^*))$ we consider the unique flow line of $\Psi$, which starts moving in the point $\Psi_{s^*}^{-1}(q^*) = \Psi_{-s^*}(q^*)$ at time $\tau = 0$. Applying now statement (\ref{eta.2}) to the initial point $q_0:= \Psi_{s^*}^{-1}(q^*)\in \pi^{-1}(\gamma(0))$ we can easily infer, that the function
$$
\eta^{s^*,q^*}(\tau):= \Psi(\tau,\Psi_{s^*}^{-1}(q^*)) \quad
\textnormal{for} \,\,\, \tau \in (0,L)
$$
satisfies the first property in (\ref{eta}),
where $\textnormal{dom}(\eta^{(s^*,q^*)})\subset \sphere^1$ can be any open and connected subset $\not=\sphere^1$, which contains the prescribed point $s^*$ -
here being identified with some open subinterval of $(\rel/L\ganz) \setminus \{0\}$ - and we 
also verify that
$$
\eta^{(s^*,q^*)}(s^*)=\Psi(s^*,\Psi_{s^*}^{-1}(q^*))
=\Psi_{s^*-s^*}(q^*)=q^*, 
$$
which is just the second property in (\ref{eta}).
Uniqueness of such a local, horizontal smooth lift of 
$\gamma$ follows easily from the above construction.       \\ 
2) First of all, since $F:\Sigma \longrightarrow \sphere^3$
is required to be an immersion whose image is $\pi^{-1}(\textnormal{trace}(\gamma))$
and since $\Sigma$ is compact,
there is some $\delta =\delta(F)>0$ such that for an arbitrarily fixed $q \in \pi^{-1}(\textnormal{trace}(\gamma))$ the preimage $F^{-1}(B_{\delta}(q))$ consists of finitely 
many disjoint open subsets of $\Sigma$, which are mapped diffeomorphically onto their images
in $B_{\delta}(q) \cap \pi^{-1}(\textnormal{trace}(\gamma))$
via $F$. Moreover, from the first part of the
lemma we infer the existence of horizontal smooth
lifts $\eta:\textnormal{dom}(\eta) \to \pi^{-1}(\textnormal{trace}(\gamma))$
of $\gamma\lfloor_{\textnormal{dom}(\eta)}$ w.r.t. $\pi$.
Now, we infer from Lemma \ref{Hopf.Tori} and from the
compactness of $\sphere^1$, that there is some 
$C=C(\gamma)>0$ such that
$$
|\eta'|^2(s)=\frac{|\gamma'|^2(s)}{4} \leq C,
\quad \forall\, s \in \textnormal{dom}(\eta),
$$
holds for every horizontal lift
$\eta$ of $\gamma\lfloor_{\textnormal{dom}(\eta)}$ w.r.t. $\pi$. Hence, there is some $\epsilon=\epsilon(F,\gamma)>0$, such that for every $\tilde s \in \sphere^1$ there is 
some $\tilde q \in \sphere^3$ such that
$\textnormal{trace}(\eta\lfloor_{\sphere^1 \cap B_{\epsilon}(\tilde s)}) \subset B_{\delta}(\tilde q)$. 
Now, using the fact that $F$ maps each connected component 
of $F^{-1}(B_{\delta}(\tilde q))$ diffeomorphically onto 
its image in $B_{\delta}(\tilde q) \cap \pi^{-1}(\textnormal{trace}(\gamma))$,
we obtain the existence of at least one smooth map
$\eta_F:\sphere^1 \cap B_{\epsilon}(\tilde s) \longrightarrow F^{-1}(B_{\delta}(\tilde q))$ satisfying $F\circ \eta_F = \eta$ on $\sphere^1 \cap B_{\epsilon}(\tilde s)$, and thus also
$\pi \circ F\circ \eta_F = \gamma$
on $\sphere^1 \cap B_{\epsilon}(\tilde s)$. Again in combination with the first part of the lemma, we finally infer from this construction, that for every $s^* \in \sphere^1$ and every $x^*\in (\pi \circ F)^{-1}(\gamma(s^*))$
there is a smooth map 
$\eta_F^{(s^*,x^*)}:\sphere^1 \cap B_{\epsilon}(s^*) 
\to \Sigma$, which possesses the two desired properties in (\ref{eta_F}). Moreover, it immediately follows from
the horizontal property of every constructed lift $\eta$ of $\gamma$ w.r.t. $\pi$ in the first part of this lemma, that every constructed lift $\eta_F$ of $\gamma$
w.r.t. $\pi \circ F$ intersects the fibers of $\pi \circ F$ perpendicularly w.r.t. the induced pullback metric $F^{*}(g_{\textnormal{euc}})$. Hence, the last assertion of the lemma now turns out to be evident. 
\qed
\noindent
\begin{remark} 
It is important to understand, that horizontal 
lifts w.r.t. $\pi$ of closed curves in $\sphere^2$ 
would in general not close up in $\sphere^3$.     
As explained in \cite{Pinkall}, p. 381, a horizontal lift $\eta$ w.r.t. $\pi$ of a simple, closed path $\gamma:\sphere^1 \longrightarrow \sphere^2$, 
which performs $k \geq 2$ loops and encloses the area 
$A$ on $\sphere^2$, closes up, if and only if there 
holds the relation $A=\frac{4\,\pi}{k}$. 
Consider for example the standard parametrization
$p(\varphi):= (\cos (2\varphi) + j \sin(2\varphi))$
of a great circle in $\sphere^2$.
Its preimage w.r.t. $\pi$ is the Clifford torus, 
and one can easily infer at first from Remark \ref{Surfaces.genus.1}, that the Clifford torus is conformally equivalent to the special parallelogram $D$ in $\com$ with vertices $(0,0)$, $(2\pi,0)$, $(\pi,\pi)$ and $(3\pi,\pi)$, and then secondly check that any horizontal lift $\eta$ 
of $p$ corresponds to a certain diagonal in $D$, 
which indeed closes up exactly when the parameter 
$\varphi$ reaches the value $2\pi$.  
\qed 
\end{remark}
\noindent 
Finally, we are going to employ the complete integrability of the Euler-Lagrange-equation of the elastic energy $\Wil$ -
see the right hand side of formula (\ref{first.variation}) - 
and the Theory of Elliptic Integrals and Jacobi Elliptic Functions, in order to precisely compute the critical values
of the elastic energy $\Wil$, in particular in order 
to exclude critical values of $\Wil$ in the surprisingly  large interval $\big{(}2\pi, \frac{8\pi}{\sqrt{2}} \big{]}$.
\begin{proposition} \label{Energy.gap}
\begin{itemize}	
	\item[1)] Up to isometries of $\sphere^2$, there are only countably many different smooth closed curves
	$\gamma:\sphere^1 \longrightarrow \sphere^2$, parametrized with constant speed, which are critical points of the elastic energy $\Wil$, called ``closed 
	elastic curves on $\sphere^2$''. Vice versa, for each pair of positive integers $(m,n)$ with $\textnormal{gcd}(m,n)=1$ and $\frac{m}{n} \in (0,2-\sqrt 2)$ there is - up to isometric equivalence - a unique arc-length parametrized elastic curve $\gamma_{(m,n)}$ in $\sphere^2$, which closes up after $n$ periods and traverses some fixed great circle exactly $m$ times.
	\item[2)] There are no critical values of the elastic energy $\Wil$ in the interval 
	$\big{(} 2\pi, \frac{8\pi}{\sqrt{2}} \big{]}$.
	\end{itemize}
\end{proposition}
\noindent
\underline{Proof of the first part of Proposition \ref{Energy.gap}:}
As explained in \cite{Langer.Singer}, every smooth closed stationary curve $\gamma:[a,b]/(a \sim b) \to \sphere^2$ of 
the elastic energy $\Wil$ satisfies 
the differential equation
\begin{equation}  \label{elastic.curve}
2 \, \Big{(} \nabla^{\perp}_{\frac{\gamma'}{|\gamma'|}}\Big{)}^2
(\vec{\kappa}_{\gamma}) + |\vec{\kappa}_{\gamma}|^2 \vec{\kappa}_{\gamma} + \vec{\kappa}_{\gamma}
\equiv \nabla_{L^2} \Wil(\gamma)  
\equiv 0  \quad \textnormal{on} \,\,\,[a,b].
\end{equation}
Now, Langer and Singer have pointed out in \cite{Langer.Singer} that equation \eqref{elastic.curve} holds for a non-geodesic closed curve $\gamma$, if and only if the signed curvature $\kappa_{\gamma}$ of its parametrization with speed $|\gamma'|\equiv 1$ satisfies the ordinary differential equation
\begin{equation}  \label{curve.ODE.4th.order}
\Big{(} \frac{d\kappa}{ds} \Big{)}^2 =
-\frac{1}{4} \, \kappa^4  - \frac{1}{2} \kappa^2 + A,  
\quad \textnormal{on} \,\,\,  \rel,
\end{equation}
for some integration constant $A \in \rel$, depending on the respective solution $\gamma$ of (\ref{elastic.curve}).
Moreover, for such non-geodesic solutions $\gamma$ of (\ref{elastic.curve}), i.e. having non-constant curvature $\kappa_{\gamma} \not \equiv 0$, equation (\ref{curve.ODE.4th.order}) is equivalent to the KdV-type-equation
\begin{equation}  \label{curve.ODE.3rd.order}
\Big{(}\frac{du}{ds} \Big{)}^2
= - \, u^3 \, - 2 \, u^2 + 4\,A\,u,  
\quad \textnormal{on} \,\,\, \rel,
\end{equation}
to be satisfied by the square $\kappa^2_{\gamma}$ of the curvature of the arc-length parametrized solution $\gamma$, which is simply obtained by pointwise multiplication of equation (\ref{curve.ODE.4th.order}) with $4\, \kappa^2_{\gamma}$.
One can easily verify that for every non-geodesic closed solution $\gamma$ of (\ref{elastic.curve}), the polynomial
$P(x):=x^3+2\,x^2-4A\,x$, occuring on the right hand side of equation (\ref{curve.ODE.3rd.order}), must have three 
different real roots $-\alpha_1< 0=\alpha_2<\alpha_3$
satisfying the algebraic relations
\begin{equation}  \label{relations.alphas}
\alpha_1 - \alpha_3 = 2  \quad \textnormal{and} \quad
\alpha_1 \, \alpha_3 = 4 \,A.
\end{equation}
Now, as explained in Section 2 of \cite{Langer.Singer} the
solutions $u$ of equation (\ref{curve.ODE.3rd.order}) are 
exactly given by the Jacobi Elliptic Functions
\footnote{See e.g. Chapter 2 in \cite{Lawden} or 
pp. 18--32 in \cite{Byrd} for an introduction to this subject.} 
of the particular type:
\footnote{See here the precise computation of the ``annual'' 
shift of the perihelion of a relativistic planetary orbit 
in Section 5.5 of \cite{Lawden} or also Appendix A.1 of \cite{Mueller.Spener.2020}, in order to obtain a clean derivation of formulae (\ref{Jacobi.function}) and (\ref{relations.r.p}).} 
\begin{equation}  \label{Jacobi.function}
u(s)= \alpha_3 \,\textnormal{cn}^2(r \cdot s\,;\,p),  
\quad \forall \, s \in \rel,
\end{equation}
with
\begin{equation}  \label{relations.r.p}
r:=\frac{1}{2} \sqrt{\alpha_1 + \alpha_3} \equiv
\frac{1}{\sqrt 2} \sqrt{\alpha_3 + 1}  \quad \textnormal{and} \quad
p:=\sqrt{\frac{\alpha_3}{\alpha_3 + \alpha_1}}
\equiv \frac{1}{\sqrt 2} \, 
\sqrt{\frac{\alpha_3}{\alpha_3+1}}.
\end{equation}
Hence, the modulus $p$ occurring in formula (\ref{Jacobi.function}) has to be contained in the 
open interval $\big{(}0,\frac{1}{\sqrt 2}\big{)}$, 
and combining formulae (\ref{relations.alphas}) and (\ref{relations.r.p}) one can express the roots $\alpha_3$ and $\alpha_1$ of the polynomial $P(x):=x^3+2\,x^2-4Ax$, occuring in equation (\ref{curve.ODE.3rd.order}), in terms of $p$:
$$
\alpha_3 = \frac{2p^2}{1-2p^2} \quad \textnormal{and} \quad
\alpha_1 = \frac{2p^2}{1-2p^2}+2 = \frac{2-2p^2}{1-2p^2}. 
$$
Moreover, combining these formulae again with formula (\ref{relations.alphas}) the modulus $p$ in formula (\ref{Jacobi.function}) automatically yields the integration constant $A$ appearing in equations (\ref{curve.ODE.4th.order}) and (\ref{curve.ODE.3rd.order}), and also the frequency $r$ 
in (\ref{Jacobi.function}):
\begin{equation} \label{formula.A}
	A = \frac{1}{4} \, \alpha_1 \, \alpha_3 
	= \frac{p^2-p^4}{(1-2p^2)^2}>0, \quad
	r=\frac{1}{\sqrt{2-4\,p^2}} 
	\in \big{(}\frac{1}{\sqrt{2}},\infty\big{)}.
\end{equation}
Formula (\ref{Jacobi.function}) particularly implies that every arc-length parametrized solution $\gamma$ of equation (\ref{elastic.curve}) performs a periodic path on 
$\sphere^2$ with
\begin{equation}  \label{period} 	
\textnormal{one period of}\, \gamma 
=4\, \frac{K(p)}{r} = 4 \,\sqrt{2-4\,p^2} \,K(p)
\end{equation} 
on account of formula (2.2.5) in \cite{Lawden} and 
formula (\ref{formula.A}) above, where
\begin{equation}  \label{first.elliptic.complete}
K(p):=\int_0^{\pi/2} \frac{1}{\sqrt{1-p^2 \sin^2(\varphi)}} \, d\varphi
\end{equation}
denotes the complete elliptic integral of the first kind with
parameter $p\in [0,\frac{1}{\sqrt 2}]$; see \cite{Lawden}, Sections 3.1 and 3.8, and \cite{Byrd}, pp. 8--17. Another 
important consequence of equation (\ref{Jacobi.function})
is, that every arc-length parametrized solution $\gamma$ of (\ref{elastic.curve}) possesses a well-defined wavelength $\Lambda(\gamma) \in \rel$
whose quotient $\frac{\Lambda}{2\pi}$ has to be rational, because $\gamma$ is supposed to be a closed curve on $[0,L]$, with $L:=\textnormal{length}(\gamma)$.
\footnote{See here also Sections 2, 3 and 5 in \cite{Langer.Singer} and Section 4 in \cite{Langer.Singer.1987} for further motivation, to attach an essentially 
unique ``wavelength'' $\Lambda$ to any 
closed solution $\gamma$ of equation (\ref{elastic.curve}).}
Hence, for every non-geodesic, closed and arc-length parametrized
solution $\gamma$ of equation (\ref{elastic.curve}) there is a unique pair $(m,n) \in \nat \times \nat$ of positive integers with $\textnormal{gcd}(m,n)=1$, such 
that $\Lambda = \frac{m}{n} \, 2\pi$, which means 
geometrically that this particular path $\gamma$ 
closes up after $n$ periods respectively ``lobes'' -  
whose common lengths are given by formula (\ref{period}) - 
and traverses some fixed great circle in $\sphere^2$\, 
$m$ times, while its arc-length parameter $s$  
runs from $0$ to $L$. 
\footnote{Compare here also with Section 5.4 in 
\cite{Mandel.2018} and with Section 3 in \cite{Mueller.Spener.2020} for closed elastica 
in the hyperbolic plane.} 
Hence, up to isometric equivalence we are able to 
``count'' non-geodesic, closed and arc-length parametrized 
elastica on $\sphere^2$ systematically, which has already 
proved the first assertion of the first part of  
Proposition \ref{Energy.gap}.            \\
Now, in order to prove the entire classification of closed elastica on $\sphere^2$, as asserted in the first part 
of Proposition \ref{Energy.gap}, we have to understand 
vice versa, for which quotients $\frac{m}{n}$ of 
coprime, positive integers there are actually  
closed elastic curves ``$\gamma_{(m,n)}$'', possessing the 
aforementioned two geometric properties, according to the 
respective pairs $(m,n)$. 
To this end, we employ the second quantitative ingredient 
of the proof of Proposition \ref{Energy.gap}, namely 
the following formula (\ref{wave.length}), taken directly 
from Section 4 of \cite{Langer.Singer.1987}, p. 148,
which expresses the wavelength $\Lambda$ of a non-geodesic, closed, arc-length parametrized solution $\gamma$ of equation (\ref{elastic.curve}) as a function of the above 
parameter $p \in \big{(}0,\frac{1}{\sqrt 2}\big{)}$:
\begin{equation} \label{wave.length}
\Lambda(\gamma(p)) =
2 \pi \, \varepsilon \,\Lambda_0(\psi(p),p)
- 2 \,(3-4\,p^2)\,\sqrt{1-(1- p^2) \sin^2(\psi(p))}\, 
\sin(\psi(p))\, K(p)
\end{equation}
where $\psi(p):=\arcsin\Big{(} \sqrt{8} \, \frac{\sqrt{1-2p^2}}{3-4p^2} \Big{)}$, $\varepsilon:=
\frac{4 \,p^2-1}{|4 \,p^2 - 1|}$, $\Lambda_0$ denotes the Heuman-Lambda-function - see \cite{Byrd}, pp. 35--37 and 
pp. 344--349 - and $K(p)$ had been introduced in (\ref{first.elliptic.complete}).
By means of the computations in Section 3.8 in \cite{Lawden} one can verify, that the first derivative w.r.t. $p$ of the function in (\ref{wave.length}) reads:
\footnote{Compare here also with Proposition 4.1 in \cite{Langer.Singer.1987} for further information.}
\begin{equation}  \label{Derivative.Lambda} 
\frac{d\Lambda(\gamma(p))}{dp} = \frac{2 \,\sqrt 8 \,
[(1- p^2) \,K(p) - E(p)]}{p \, (1-p^2)\, \sqrt{1-2p^2}  
\,(3 - 4 p^2) \, \sqrt{1 - p^2 \,\sin^2(\psi(p))}}
\end{equation}  
for $p\in \big{(}0,\frac{1}{\sqrt 2}\big{)}$, where
$E(p):=\int_0^{\pi/2}\sqrt{1-p^2 \sin^2(\varphi)}\, d\varphi$
denotes the complete elliptic integral of the second kind. 
Since we have $(1-p^2)\,K(p)-E(p)<0$ for $p\in \big{(}0,\frac{1}{\sqrt 2}\big{)}$, and 
furthermore $\sin(\psi(0))=\frac{\sqrt 8}{3}$, 
$\Lambda_0\big{(}\arcsin\big{(}
\frac{\sqrt 8}{3}\big{)},0\big{)}=\frac{\sqrt 8}{3}$, 
$\sin(\psi(\frac{1}{2}))=1$ and $\Lambda_0\big{(}\frac{\pi}{2},\frac{1}{2}\big{)}=1$,
formula (\ref{Derivative.Lambda}) shows us, that the function 
$p \mapsto \Lambda(\gamma(p))$ in (\ref{wave.length}) 
decreases \underline{strictly monotonically} from 
its initial value $\Lambda(0) = 
- 2 \pi \cdot \frac{\sqrt 8}{3} - 2\, \frac{\sqrt{2}}{3} \, \pi = -\sqrt 2 \,(2 \pi)$ to its minimal value
$\min_{p \in [0,\frac{1}{\sqrt 2}]} \Lambda(p)
= \lim_{p \nearrow \frac{1}{2}} \Lambda(p) = 
-2 \pi - 2\, K\big{(}\frac{1}{2}\big{)} \approx -9,65469$
as $p$ increases from $0$ to $\frac{1}{2}$, then
jumps at the point $p=\frac{1}{2}$ from its minimal value $\lim_{p \nearrow \frac{1}{2}} \Lambda(p)$ to its maximal value
\footnote{Note that the function $[p \mapsto \varepsilon(p)]$ 
jumps at $p=\frac{1}{2}$ from $-1$ to 
$1$ and that $\Lambda_0\big{(}\frac{\pi}{2},\frac{1}{2}\big{)}=1$.}   
$\max_{p \in [0,\frac{1}{\sqrt 2}]} \Lambda(p)
= \lim_{p \searrow \frac{1}{2}} \Lambda(p) = 
2 \pi - 2 \, K\big{(}\frac{1}{2}\big{)}
\approx 2,91169$, and finally decreases 
\underline{strictly monotonically} from this maximum 
to its final value $\Lambda\big{(}\frac{1}{\sqrt 2}\big{)} 
= 2\, \pi \Lambda_0\big{(}0,\frac{1}{\sqrt 2}\big{)}=0$. 
Since we are allowed to shift the value of 
the wavelength $\Lambda(p)$ in the first interval
$[-2 \pi-2 \,K\big{(}\frac{1}{2}\big{)},-\sqrt 2 \,(2 \pi)]$
about the height of the jump of $\Lambda(p)$, i.e. 
about $4 \pi$, to the right in $\rel$ into the interval 
$[2 \pi - 2\,K\big{(}\frac{1}{2}\big{)}, 
2 \pi \,(2-\sqrt 2)]$, we conclude that for each 
coprime pair $(m,n) \in \nat \times \nat$ 
satisfying
\begin{equation}   \label{correct.frequences}
\frac{m}{n} \in \Big{(} 0, 
1-\frac{K\big{(}\frac{1}{2}\big{)}}{\pi} \Big{]} \cup
\Big{[} 1-\frac{K\big{(}\frac{1}{2}\big{)}}{\pi},  
2-\sqrt 2 \Big{)}=(0,2-\sqrt 2)
\end{equation}
there is a ``unique'' arc-length parametrized elastic curve $\gamma_{(m,n)}$, which closes up after $n$ periods and traverses some fixed great circle $C$ in $\sphere^2$ exactly $m$ times - ``unique'' only up to the action of all those isometries of 
$\sphere^2$, which leave the great circle $C$ invariant, 
just as claimed in the first part of the proposition. 
\qed 
\noindent \\\\
\underline{Proof of the second part of Proposition \ref{Energy.gap}:} 
Relying on the proof of the first part of Proposition 
\ref{Energy.gap} we attempt to prove its second part 
in a ``straightforward manner'', which means that we fix 
a pair of coprime integers $(m,n) \in \nat \times \nat$ satisfying condition (\ref{correct.frequences}) 
and that we simply combine formulae (\ref{Jacobi.function}) and (\ref{period}) with the definition of the elastic energy in (\ref{elastic.energy.functional}), in order to compute both 
length and elastic energy of the unique solution 
$\gamma_{(m,n)}$ of equation (\ref{elastic.curve}) directly. 
This is actually possible, because the given data ``$(m,n)$''
determine the value of the wavelength 
$\Lambda=\frac{m}{n} \,2\pi$ and thus also a unique value of 
the parameter $p=p(\gamma_{(m,n)})$ in (\ref{relations.r.p}), 
inverting the strictly monotonic wavelength-function 
$p \mapsto \Lambda(\gamma(p))$ in (\ref{wave.length}). 
Hence, also recalling that we only work with arc-length parametrized elastic curves $\gamma_{(m,n)}$, we can firstly compute by means of formula (\ref{period}), abbreviating here $p=p(\gamma_{(m,n)})$:
\footnote{Compare here also with Proposition 9 in 
\cite{Mandel.2018} and with Propositions 3.3 and 3.4 in 
\cite{Mueller.Spener.2020}.}  
\begin{equation}  \label{length.directly}
\textnormal{length}(\gamma_{(m,n)}) = 
n \,\textnormal{periods of}\, \gamma_{(m,n)} 
= 4n\, \frac{K(p)}{r} = 4n \, \sqrt{2-4\,p^2} \,K(p), 
\end{equation} 
and together with formulae (\ref{elastic.energy.functional}) and  (\ref{Jacobi.function}) and with formulae 
(3.4.15) and (3.4.27) in \cite{Lawden} we obtain furthermore: 
\begin{eqnarray} \label{energy.directly}    
\Wil(\gamma_{(m,n)}) = \int_{0}^{4n \,\frac{K(p)}{r}}
1 + \kappa^2_{\gamma_{(m,n)}}(s) \, ds                 
= \int_0^{4n \,\frac{K(p)}{r}} 
1+\alpha_3 \,\textnormal{cn}^2(r \cdot s\,;\,p) \,ds \nonumber  \\ 
= 4n \,\sqrt{2-4\,p^2} \,K(p) 
+ \frac{2p^2}{1-2p^2} \,\frac{4n}{r} 
\int_0^{K(p)}  \textnormal{cn}^2(u;\,p) \, du   \qquad      \\
= 4n \,\sqrt{2-4\,p^2} \,K(p) 
+ \frac{2p^2}{1-2p^2} \, 4n \,\sqrt{2-4\,p^2} \,\,\frac{1}{p^2}\,
\Big{[} E(p) - (1-p^2) K(p)\, \Big{]}                  \nonumber    \\
=  4n \,\sqrt{2-4\,p^2} \,K(p) 
\Big{(} 1+ \frac{2p^2}{1-2p^2} \Big{)}  
+ \frac{16n}{\sqrt{2-4p^2}} \, (E(p)-K(p)) \nonumber   \\     
=  \frac{8n}{\sqrt{2-4p^2}} \,(2E(p)-K(p)),         \nonumber
\end{eqnarray} 
a formula which has essentially already appeared in the proof 
of Corollary 6.4 in \cite{Mueller.Spener.2020} for 
closed wavelike elastica in the hyperbolic plane, being 
similar to the slightly simpler formula (46) 
in \cite{Mandel.2018} for the elastic energy of closed 
orbitlike elastica in the hyperbolic plane.
\footnote{Compare here also with Lemma 4.1 in \cite{Mueller.Spener.2020} and with p. 19 in \cite{Langer.Singer}.}  	   
Now, formula (\ref{energy.directly}) is not only useful  
for numerical purposes, i.e. in order to compute elastic energies 
$\Wil(\gamma_{(m,n)})$ effectively, but it also 
meets the aim of the second part of Proposition 
\ref{Energy.gap}, to rigorously determine a ``rather accurate'' 
lower bound for all possible elastic energies $\Wil(\gamma_{(m,n)})$ of non-geodesic elastic curves. 
The key observation in this situation is, that the function 
$f(p):=\frac{1}{\sqrt{1-2p^2}} \,(2E(p)-K(p))$ - appearing  
on the right hand side of formula (\ref{energy.directly}) - is 
actually \underline{strictly monotonically increasing} 
on the entire open interval $\big{(}0,\frac{1}{\sqrt 2}\big{)}$. 
In order to prove this, we firstly compute by means of formulae 
(3.8.7) and (3.8.12) in \cite{Lawden}:    
\begin{eqnarray}   \label{derivative.2E-K}  
	\frac{d}{dp}\big{(}2E(p) - K(p)\big{)} 
	= 2 \,\frac{E(p)- K(p)}{p} 
	- \frac{E(p) - (1-p^2)\,K(p)}{p\,(1-p^2)}       \\  
	= \frac{1-2p^2}{p\,(1-p^2)} \,E(p) - \frac{1}{p}\, K(p)
	<\frac{1}{p} \, (E(p)-K(p))<0,    \nonumber
\end{eqnarray} 
for $p\in \big{(}0,\frac{1}{\sqrt 2}\big{)}$, 
showing first of all that the function $p\mapsto 2E(p) - K(p)$ 
decreases monotonically from $\frac{\pi}{2}$ in $p=0$ 
to $2E\big{(}\frac{1}{\sqrt 2}\big{)} 
- K\big{(}\frac{1}{\sqrt 2}\big{)}>0$ in 
$p=\frac{1}{\sqrt 2}$, and we can continue by deriving 
the function $f$ w.r.t. $p$:  
\begin{eqnarray*}  
\frac{df(p)}{dp} 
=  \frac{1}{\sqrt{1-2p^2}} \Big{(} \frac{1-2p^2}{p(1-p^2)} 
E(p) - \frac{1}{p} \,K(p) \Big{)} 
+ \frac{2p}{(1-2p^2)^{3/2}} (2E(p)-K(p))                 \\
=\frac{\sqrt{1-2p^2}}{p\,(1-p^2)} E(p) 
- \frac{1}{p\,\sqrt{1-2p^2}}  K(p) 
+ \frac{4p}{(1-2p^2)^{3/2}} E(p) 
- \frac{2p}{(1-2p^2)^{3/2}}  K(p)                     \\
= \frac{1}{p\,(1-p^2)\,(1-2p^2)^{3/2}} E(p) 
- \frac{1}{p\,(1-2p^2)^{3/2}} K(p)                     \\  
= \frac{1}{p\, (1-2p^2)^{3/2}} \,
\Big{(} \frac{1}{1-p^2}\, E(p)-K(p) \Big{)}. 
\end{eqnarray*} 	  
Moreover, we see that the function 
$g(p):=\frac{1}{1-p^2}\, E(p) - K(p)$ satisfies 
$g(0)=0$, and that by formulae 
(3.8.7) and (3.8.12) in \cite{Lawden} its derivative is: 
\begin{eqnarray*}  
	\frac{dg(p)}{dp}
	= \frac{E(p)- K(p)}{(1-p^2)\,p} 
	+ \frac{2p}{(1-p^2)^2} \, E(p) 
	- \frac{E(p) - (1-p^2)\,K(p)}{p\,(1-p^2)}              \\   
    =\frac{2p}{(1-p^2)^2} \,E(p) - \frac{p}{1-p^2} \,K(p) 
    > \frac{2p}{(1-p^2)} \,E(p) - \frac{p}{1-p^2} \,K(p)   \\
    = \frac{p}{(1-p^2)} \,(2 \,E(p) - K(p)) >0     
\end{eqnarray*} 
for every $p\in \big{(}0,\frac{1}{\sqrt 2}\big{)}$, 	  
using that here $(1-p^2) \in \big{(}\frac{1}{2},1\big{)}$ 
and $E(p)>0$, and furthermore that by (\ref{derivative.2E-K})    
$\min_{p\in \big{[}0,\frac{1}{\sqrt 2}\big{]}} 
(2 \,E(p) - K(p)) = 2E\big{(}\frac{1}{\sqrt 2}\big{)} 
- K\big{(}\frac{1}{\sqrt 2}\big{)}>0$. 
Hence, we can infer that $g(p)>0$ for 
$p\in \big{(}0,\frac{1}{\sqrt 2}\big{)}$ and that 
therefore $f$ increases strictly monotonically from 
$f(0)= \frac{\pi}{2}$ to $\infty$, as  
$p$ runs from $0$ to $\frac{1}{\sqrt 2}$.  
Hence, without even guessing for which pair $(m,n)$ 
of coprime integers - respectively for which value 
of the modulus $p=p(\gamma_{(m,n)})\in 
\big{(}0,\frac{1}{\sqrt 2}\big{)}$ - 
the elastic energy $\Wil(\gamma_{(m,n)})$ in 
(\ref{energy.directly}) might attain its minimal value 
among all non-geodesic elastica on $\sphere^2$, we can 
roughly, but rigorously estimate by means of formula (\ref{energy.directly}) and 
$\inf_{p\in [0,\infty)} f(p)=f(0)=\frac{\pi}{2}$:    
\begin{equation} \label{rough.estimate.1}  
\Wil(\gamma_{(m,n)}) 
> \frac{16}{\sqrt{2}} f(0) = \frac{8\pi}{\sqrt 2} 
\approx 17,771532  
\end{equation} 
for coprime pairs $(m,n)\in \nat \times \nat$ satisfying condition (\ref{correct.frequences}), where we have 
also used the fact that we must have ``$n\geq 2$'' on  
account of condition (\ref{correct.frequences}), 
that $p(\gamma_{(m,n)}) \in \big{(}0,\frac{1}{\sqrt 2}\big{)}$ 
and that $f$ increases strictly monotonically from $f(0)$ 
to $f(p(\gamma_{(m,n)}))$, for any fixed coprime pair 
$(m,n)$ satisfying condition (\ref{correct.frequences}).  
\qed 

\begin{remark}   \label{remark.for.proposition.6}
In order to assess the mathematical eloquence of our 
method leading to estimate (\ref{rough.estimate.1}) in the 
second part of Proposition \ref{Energy.gap}, 
we shall sketch here a completely different method, 
which yields the rigorous lower bound  
$16 \,\sqrt{\frac{\pi}{3}}$ for the elastic energies of 
all non-geodesic elastica $\gamma_{(m,n)}$ on $\sphere^2$, 
for any pair $(m,n)$ of positive, coprime integers 
with $\frac{m}{n} \in (0,2-\sqrt 2)$.
\footnote{Obviously, the number 
$16 \,\sqrt{\frac{\pi}{3}}$ is slightly smaller than our 
threshold $\frac{8\pi}{\sqrt 2}$ from (\ref{rough.estimate.1}).}  
To this end, we notice first of all, that equation (2.2) 
in \cite{Heller} is exactly equation (\ref{curve.ODE.4th.order}) above, simply with $\nu=-A$, $G=1$ and $\mu=-\frac{1}{2}$ in equation (2.2) of \cite{Heller}, and that the formulae on p. 350 in \cite{Heller} yield the coefficients
\begin{eqnarray}  \label{g2.g3}
a_2 = \frac{1}{48} - \frac{A}{4}  \qquad \textnormal{and}
\qquad  a_3 = \frac{1}{1728} + \frac{A}{48}
\end{eqnarray}
of a certain polynomial equation\,$y^2=4\,x^3-a_2 \,x-a_3$,
whose set of solutions $[x:y:1] \in \CP^2$ yields a particular elliptic curve $E_{(m,n)} \subset \mathbf{CP}^2$, 
which turns out to contain all relevant, computable information about the fixed elastic curve $\gamma_{(m,n)}$. One can easily derive from formula (\ref{g2.g3}) that the discriminant
\begin{equation} \label{discriminant}
D(F):= a_2^3 - 27 \, a_3^2
\end{equation}
of the polynomial $F(x) = 4\,x^3-a_2\,x-a_3$
is zero if and only if $A=0$ or $A=-\frac{1}{4}$. But
formula (\ref{formula.A}) tells us that $A$ actually has to
be a positive real number, which rules out both of these possibilities.
Hence, for any $p\in \big{(}0,\frac{1}{\sqrt{2}}\big{)}$, respectively for every $A>0$, the corresponding polynomial 
$F(x) = 4\, x^3 - a_2 \, x - a_3$ has a non-vanishing discriminant. Since formulae (\ref{g2.g3}) and
(\ref{discriminant}) show that the discriminant is negative
for large $A$, we conclude more precisely that for every $A>0$ there holds $D(F)<0$, and on account of $a_3>0$ in formula (\ref{g2.g3}) for every $A>0$, this finally implies that
\begin{equation} \label{Uniformization.trick}
\frac{a_2^3}{a_2^3-27 \,a_3^2}
\equiv \frac{a_2^3}{D(F)} \in (-\infty,1),
\end{equation}
for every $A>0$. More precisely, we have:
\begin{equation} \label{Uniformization.Trick.2}
\frac{a_2^3}{a_2^3-27 \,a_3^2} \searrow -\infty \,\,\, \textnormal{as} \,\,\,A \searrow 0 \quad \textnormal{and} \quad
\frac{a_2^3}{a_2^3-27 \,a_3^2} \nearrow  1  \,\,\, \textnormal{as} \,\,\, A \nearrow \infty.
\end{equation}
Now, already due to $D(F)\not =0$, for every fixed $A>0$, the Uniformization Theorem - see e.g. Theorem 2.9 in \cite{Apostol} - guarantees us the existence of a unique lattice
$\Omega\equiv \Omega_{(m,n)}:=\omega_1 \ganz 
\oplus \omega_2 \ganz$ in $\com$, with $\Im(\frac{\omega_2}{\omega_1})>0$, such that the lattice invariants
\begin{eqnarray*}
g_2(\Omega):= 60\, G_4(\Omega) \equiv 60 \, \sum_{\omega \in \Omega\setminus \{0\}} \frac{1}{\omega^4}  \,\, 
\textnormal{and} \,\,  
g_3(\Omega):= 140 \, G_6(\Omega) \equiv 140 \, \sum_{\omega \in \Omega\setminus \{0\}} \frac{1}{\omega^6}
\end{eqnarray*}
satisfy $g_2(\Omega)=a_2$ and $g_3(\Omega)=a_3$.
Moreover, the lattice $\Omega= \Omega_{(m,n)}$ yields a unique \underline{Weierstrass-$\wp$-function}
$\wp(z,\Omega):=\frac{1}{z^2} + 
\sum_{\omega \in \Omega\setminus \{0\}}
\Big{(} \frac{1}{(\omega-z)^2} - \frac{1}{\omega^2} \Big{)}$, 
defined in every $z\in \com\setminus \Omega$, 
which solves the complex differential equation
\begin{equation}  \label{Weierstrass.equation}  
(\wp'(z))^2= 4 \,(\wp(z))^3 - a_2 \, \wp(z) - a_3 
\equiv F(\wp(z)), \quad \forall \, z\in \com/\Omega, 
\end{equation}  
and therefore ``parametrizes'' the elliptic curve $E_{(m,n)}$.
\footnote{See here \cite{Silverman1}, p. 165--170, and 
formula (3.4) in \cite{Heller}.} 
Now, the exact relation between the elliptic curve $E_{(m,n)}$ 
and the considered solution $\gamma_{(m,n)}$ of equation 
(\ref{elastic.curve}) is expressed by Lemmata 1 and 2 in \cite{Heller} and by equation (\ref{Weierstrass.equation}) above.
A combination of these two lemmata shows, that there exists some 
$x_0 \in \com \setminus \Big{(}\frac{1}{2} \, \Omega \oplus \rel\Big{)}$, depending on the exact parametrization of $\gamma_{(m,n)}$, such that there holds:
\begin{equation} \label{Weierstrass.p.kappa}
\wp(x + x_0) =
- i \frac{\kappa'(x)}{4}-\frac{\kappa^2(x)}{8} 
- \frac{1}{24}, \quad \textnormal{for} \,\,\, x \in \rel,
\end{equation}
where $\kappa$ still denotes the signed curvature function of $\gamma_{(m,n)}$. In combination with equation (\ref{Jacobi.function}) this implies that $\wp$ has a real primitive period, say $\omega_1$, namely the real primitive period of the Jacobi-elliptic function 
$\textnormal{cn}(r \,\cdot \,;p)$:
\begin{equation} \label{first.period}
\omega_1(\wp) = 4\, \frac{K(p)}{r} 
=  4\, \sqrt{2 - 4\,p^2} \,K(p),
\end{equation}
where we have used formula (\ref{formula.A}), just as in 
(\ref{period}). Furthermore, it will turn out to be of key importance for this method, to locate - at least qualitatively - the value of the quotient $\tau \equiv \tau_{(m,n)}=\frac{\omega_2}{\omega_1} \in \quat$ 
of the primitive periods of the lattice $\Omega=\Omega_{(m,n)}$ respectively of the corresponding Weierstrass-$\wp$-function. 
To this end, we invoke 
the proof of the Uniformization Theorem
\footnote{Compare here with the proof of Theorem 2.9 in \cite{Apostol}.}, showing that this quotient 
$\tau = \tau_{(m,n)}$ has to satisfy the fundamental equation
\begin{equation} \label{J}
J(\tau)\stackrel{!}=\frac{a_2^3}{a_2^3 - 27 \,a_3^2}
\equiv \frac{a_2^3}{D(F)} \in (-\infty,1),
\end{equation}
where we also used formula (\ref{Uniformization.trick}).
Here, $J$ denotes \underline{Klein's modular function} - see 
pp. 15--17 in \cite{Apostol} - being    $\textnormal{PSL}(2,\ganz)$-invariant and mapping the closure of every fundamental domain $R_{\Gamma(1)}\subset \com$ of the modular group $\Gamma(1) \cong \textnormal{PSL}(2,\ganz)$ meromorphically onto the entire complex plane $\com$. Since Klein's modular function $J$, restricted to the closure of its typical fundamental domain
$R_{\Gamma(1)}:=\{\tau \in \com |\, \Im(\tau) >0, 
|\tau| > 1, \,\,|\tau + \bar \tau| < 1 \}$
- see Section 2.3 in \cite{Apostol} - maps exactly the 
two subarcs of $\sphere^1$ connecting the points $\exp(i\frac{2\pi}{3})$
and $i$ respectively the points $\exp(i\frac{\pi}{3})$ 
and $i$ monotonically onto the interval $[0,1]$ 
- see here \cite{Apostol}, p. 41 - and the two lines
$\{\Re(\tau) =\pm \frac{1}{2}\} \cap \partial R_{\Gamma(1)}$ monotonically onto the ray $(-\infty,0]$, we infer that equation (\ref{J}) has exactly two
solutions $\tau_{1,2} \in \partial R_{\Gamma(1)}$, 
which are either of the form 
$\tau_{1,2} =\exp(i \,\pi(\frac{1}{2} \pm \delta))$,
for exactly one angle $\delta \in (0,\frac{1}{6}]$, or
$\tau_{1,2}= \pm \frac{1}{2} + i \, h$ for exactly one 
$h > \frac{\sqrt 3}{2}$. Using formulae (\ref{g2.g3}), 
(\ref{Uniformization.Trick.2}) and (\ref{J}), one can 
easily see that the latter case corresponds to the case 
``$0<A<\frac{1}{12}$'', whereas the former case holds 
if and only if $A\geq\frac{1}{12}$, respectively if and only 
if $p \in \big{[}\frac{1}{2} \sqrt{2-\sqrt 3},
\frac{1}{\sqrt 2} \big{)}$. Moreover, because of
\begin{equation}  \label{symmetry.imag.axis}
\Re(\tau_1) = -\Re(\tau_2) \qquad \textnormal{and} \qquad
\Im(\tau_1) = \Im(\tau_2),
\end{equation}
there holds either $\Omega=\omega_1 \ganz \oplus \tau_1 \omega_1 \ganz$ or $\overline{\Omega} = \omega_1 \ganz \oplus \tau_1 \omega_1 \ganz$.
Now, again on account of (\ref{symmetry.imag.axis}) the
two nomes $q_{1,2}:=\exp(2\pi i\,\tau_{1,2})$ satisfy
$q_2=\overline{q_1}$, and therefore formulae (6.6.4) and 
(6.6.16) in \cite{Lawden} combined with Theorem 6.3 in
Chapter I of \cite{Silverman2} yield the value of the 
``quasi-period'' respectively of the ``jump'' $\eta(\omega_1,\Omega):=\zeta(z+\omega_1,\Omega) -\zeta(z,\Omega)$ 
of the \underline{Weierstrass-$\zeta$-function}
$\zeta(z,\Omega):= \frac{1}{z} + 
\sum_{\omega \in \Omega\setminus \{0\}}
\Big{(} \frac{1}{z-\omega} + \frac{1}{\omega}
+ \frac{z}{\omega^2} \Big{)}$, for 
$z \in \com \setminus \Omega$,
corresponding to the lattice $\Omega$ - respectively corresponding to the flipped lattice $\overline{\Omega}$ - in direction of its real period $\omega_1$ in terms of:
\begin{equation} \label{Silverman2}
\eta(\omega_1,\Omega) \,\, \textnormal{or} \,\,
\overline{\eta(\omega_1,\Omega)}
= \frac{\eta(1,\ganz \oplus \tau \ganz)}{\omega_1}
= \frac{1}{\omega_1} \,\Big{(}\frac{\pi^2}{3} - 8\, \pi^2 \,
\sum_{n=1}^{\infty} \frac{q^n}{(1-q^n)^2} \Big{)},    \quad
\end{equation}
for $\tau \in \{\tau_1,\tau_2\}$ respectively 
$q\in \{q_1,q_2\}$. Hence, in either case a combination of formulae (\ref{first.period}) and (\ref{Silverman2}) with 
Theorem 5 in \cite{Heller}\footnote{Theorem 5 in \cite{Heller}
essentially follows from an integration of the real part
of equation (\ref{Weierstrass.p.kappa}) and a comparison of 
that result with the meaning of the functional 
$\Wil$ itself.} 
finally yields the desired, new formulae for both, length and 
elastic energy of the elastic curve $\gamma_{(m,n)}$:
\begin{eqnarray}  \label{length.formula}
\textnormal{length}(\gamma_{(m,n)}) 
= n \cdot \omega_1  \qquad \textnormal{and} \\
\Wil(\gamma_{(m,n)}) = 8 \cdot n \cdot \frac{\Re(\eta(1,\ganz \oplus \tau \ganz))}{\omega_1} + \frac{2}{3} \cdot n \cdot \omega_1,   \label{elastic.energy}
\end{eqnarray}
for any pair $(m,n) \in \nat \times \nat$ 
of coprime integers satisfying $\frac{m}{n} \in (0,2-\sqrt 2)$, according to condition (\ref{correct.frequences}).\\
Now, in order to derive from formula (\ref{elastic.energy}) 
the asserted lower bound $16 \,\sqrt{\frac{\pi}{3}}$ for 
all possible energies $\Wil(\gamma_{(m,n)})$, we 
firstly recall \underline{Legendre's relation}:
\begin{equation} \label{Legendre}
\omega_2 \, \eta(\omega_1,\Omega) - 
\omega_1 \, \eta(\omega_2,\Omega) = 2\pi i,
\end{equation}
holding for every lattice 
$\Omega=\omega_1 \ganz \oplus \omega_2 \ganz$ with
$\Im(\omega_2/\omega_1)>0$; see \cite{Silverman1}, p. 179.
In combination with the relation
$\eta(\omega,\Omega)=2 \,\zeta(\frac{\omega}{2},\Omega)$, 
for every $\omega \in \Omega \setminus 2 \Omega$, and with 
the obvious identity
$$
\zeta(i z, \ganz \oplus i \ganz)
= \frac{1}{i} \, \zeta\Big{(}z,\frac{1}{i} \ganz \oplus \ganz \Big{)}
= -i \, \zeta(z, \ganz \oplus i \ganz), \,\,\, 
\forall \,z \in \com \setminus \Omega,
$$
formula (\ref{Legendre}) implies especially for the rectangular lattice $\Omega=\ganz \oplus i \ganz$:
\begin{eqnarray*}
2\pi i = i \, \eta(1,\ganz \oplus i \ganz) 
- 1 \, \eta(i,\ganz \oplus i \ganz)
=  2i \, \zeta\Big{(}\frac{1}{2},\ganz \oplus i \ganz\Big{)} 
-2 \, \zeta\Big{(}\frac{i}{2},\ganz \oplus i \ganz\Big{)}                    \\
= 4i \, \zeta\Big{(}\frac{1}{2},\ganz \oplus i \ganz\Big{)}
=2i \, \eta(1,\ganz \oplus i \ganz),
\end{eqnarray*}
and therefore:
\begin{equation}  \label{eta.1.i}
\eta(1,\ganz \oplus i \ganz) = \pi.
\end{equation}
Moreover, it can be easily derived from Eisenstein's 
classical formula
$$
\frac{\pi^2}{\sin^2(\pi z)}
=\sum_{\nu=-\infty}^{\infty} \frac{1}{(z+\nu)^2},
\quad \textnormal{for} \,\,\,z \in \mathbf{C} \setminus \ganz,
$$
that there holds, again writing $q:=\exp(2\pi i \,\tau)$:
\begin{eqnarray} \label{cosec.expansion}
	\eta(1,\ganz \oplus \tau \ganz)
	= \frac{\pi^2}{3} - 8 \pi^2 \, \sum_{n=1}^{\infty} \frac{q^n}{(1-q^n)^2}
	= 2\, \zeta(2) + 2 \, \sum_{\mu=1}^{\infty} \frac{\pi^2}{\sin^2(\mu \pi\tau)}   \quad \\
	\label{second.Eisenstein}
	=\sum_{\mu =-\infty}^{\infty}  \sum_{\nu =-\infty}^{\infty} \frac{1}{(\mu \tau + \nu)^2} 
	=:G_2(\tau),                  \quad
\end{eqnarray}
where in line (\ref{cosec.expansion}) the symbol $\zeta$ denotes Riemann's ``zeta-function'' and where in the double series in line (\ref{second.Eisenstein})
the index $\nu=0$ is omitted if the second summation index $\mu$ is zero.\footnote{The second equation in line (\ref{cosec.expansion}) can be proved via the product-expansion of the $\sigma$-function and its relation to the Weierstrass-$\zeta$-function via
$\zeta(z)=\frac{\partial}{\partial z}(\ln(\sigma(z)))$,
see \cite{Lawden}, Sections 6.5 and 6.6 and pp. 183--184.}
It is important to note here, that the 
\underline{second Eisenstein-series} $G_2$, as defined in line (\ref{second.Eisenstein}), is not absolutely convergent, 
which implies that the order of summation in line (\ref{second.Eisenstein}) matters and that $G_2$ is a holomorphic function, but not a modular form on $\quat:=[\Im(\tau)>0]$. Actually, it will turn out crucial for this proof that $G_2$ transforms w.r.t. the generator $\tau \mapsto -\frac{1}{\tau}$ of the modular group $\Gamma(1)$ exactly in the following way:
\begin{equation}  \label{G2.trafo}
	G_2\Big{(}-\frac{1}{\tau} \Big{)}
	= \tau^2 \, G_2(\tau) - 2\pi i \tau,  
	\quad \forall\, \tau \in \quat,
\end{equation}
see here \cite{Apostol}, Theorem 3.1 and pp. 69--71.
\footnote{On pp. 69--71 in \cite{Apostol} Apostol 
	shows the equivalence between formula 
	(\ref{G2.trafo}) and the famous functional equation 
	$\eta\big{(}\frac{-1}{\tau}\big{)}=
	(-i\tau)^{1/2} \, \eta(\tau)$ 
	of the Dedekind-$\eta$-function, 
	which is proved classically in Theorem 3.1 of \cite{Apostol}.} 
Combining now formulae (\ref{eta.1.i}) and (\ref{second.Eisenstein}) we obtain the equation 
$G_2(i)=\pi$, and the last step of the proof 
of the second part of this proposition consists in 
the verification of the assertion:
\begin{equation}  \label{The.inequality}
	\Re G_2(\tau) \geq \pi, \quad \forall \,\tau \in
	\partial R_{\Gamma(1)}.
\end{equation}
First of all, we note that lines (\ref{Silverman2}), (\ref{cosec.expansion}) and (\ref{second.Eisenstein}) immediately imply:
\begin{equation}  \label{conjugated.G2}
	G_2\Big{(}-\frac{1}{\tau}\Big{)} =  \overline{G_2(\tau)},   \quad \forall \,
	\tau \in \partial R_{\Gamma(1)} \cap \sphere^1.
\end{equation}
Now, multiplication of formula (\ref{G2.trafo}) with $-\frac{1}{\tau}$ and inserting equation (\ref{conjugated.G2}) into the resulting equation yields the identity
\begin{equation}  \label{Im.Re.Eisenstein}
	\Im(G_2(\tau)) \,\tau_1 +  \Re(G_2(\tau)) \,\tau_2 
	= \pi, \quad \forall \,
	\tau= \tau_1 + i\tau_2 \in \partial R_{\Gamma(1)} \cap \sphere^1,
\end{equation}
which will turn out to be one of the key tools in the following argument. Moreover, as in Sections 1.14 and 3.10 of \cite{Apostol} 
one can convert the second expression in line (\ref{cosec.expansion}) into an absolutely convergent 
Fourier expansion of $G_2$, namely:
\begin{equation}  \label{Fourier.expansion.G2}
	G_2(\tau) = \frac{\pi^2}{3} - 8 \pi^2 \, \sum_{n=1}^{\infty}
	\sigma_1(n) \, \exp(2\pi in\tau),  \quad \textnormal{for}  \,\, \tau \in \quat,
\end{equation}
where $\sigma_1(n):=\sum_{d|n} d$ denotes the divisor-sum 
with power $1$. On account of the absolute convergence of 
this Fourier series in any $\tau \in \quat$, one can differentiate both sides of formula (\ref{Fourier.expansion.G2}), in order to obtain the Fourier expansion of the complex derivative of $G_2$:
\begin{equation}  \label{Derivative.Eisenstein}
\frac{\partial G_2}{\partial \tau}(\tau)
=- 16 \pi^3 i \, \sum_{n=1}^{\infty} n\,\sigma_1(n)\, 
\exp(2\pi in\tau), \quad \textnormal{for}  \,\, 
\tau \in \quat.
\end{equation}
This implies in particular the formula
\begin{equation}  \label{Derivative.Eisenstein.revisited}
	\frac{\partial G_2}{\partial \varphi}(e^{i\varphi})
	= 16 \pi^3  \, \sum_{n=1}^{\infty} n\,\sigma_1(n) 
	\exp(- 2\pi n\sin(\varphi)) \,
	\exp(i(2\pi n\cos(\varphi) + \varphi)), \quad
\end{equation}
for $\varphi \in [\frac{\pi}{3}, \frac{2\pi}{3}]$.
First of all, we consider the imaginary part of formula
(\ref{Derivative.Eisenstein.revisited}):
\begin{equation}  \label{Imaginary.Derivative.Eisenstein}
	\frac{\partial \Im G_2}{\partial \varphi}(e^{i\varphi})
	= 16 \,\pi^3  \, \sum_{n=1}^{\infty} n\,\sigma_1(n) 
	\exp(- 2\pi n \sin (\varphi))\,
	\sin(2\pi n \cos(\varphi) + \varphi), \quad
\end{equation}
for $\varphi \in [\frac{\pi}{3}, \frac{2\pi}{3}]$, 
and we will prove now \underline{in an intermediate step}, 
that this angular derivative is \underline{positive} 
for all angles $\varphi \in [\arccos(\frac{1}{4}),\frac{\pi}{2}]$, 
i.e. that the sum (\ref{Imaginary.Derivative.Eisenstein}) 
has a positive value for every 
$\varphi \in [\arccos(\frac{1}{4}),\frac{\pi}{2}]$.
To this end, we invoke the elementary estimate 
$1 \leq n \, \sigma_1(n) \leq \frac{n^2\, (n+1)}{2}$ and 
the two geometric series
\begin{equation} \label{geometric.series}
	\sum_{n=1}^{\infty} n^2 \, q^n = \frac{q + q^2}{(1-q)^3} \qquad \textnormal{and}
	\qquad \sum_{n=1}^{\infty} n^3 \, q^n 
	=\frac{q + 4 q^2+ q^3}{(1-q)^4},
\end{equation}
for any fixed $q\in [0,1)$. Applied to $q=q(\varphi):=
\exp(-2\pi\,\sin(\varphi))$ they yield:
\begin{eqnarray}  \label{estim.geom.series}
	\sum_{n=2}^{\infty} n\,\sigma_1(n) 
	\exp(- 2\pi n \sin (\varphi)) \,
	|\sin(2\pi n \cos(\varphi) + \varphi)| \qquad \\
	\leq \frac{8q^2-5q^3+4q^4-q^5}{2 \,(1-q)^4} + \frac{4q^2-3q^3+q^4}{2\,(1-q)^3}
	\leq 1,142 \cdot 10^{-4}, \,\,\, \forall \,\varphi \in \big{[}\frac{\pi}{3}, \frac{\pi}{2}\big{]}.
	\nonumber
\end{eqnarray}
Moreover, the function $\varphi \mapsto \sin(2\pi \cos(\varphi) + \varphi)$ increases monotonically from $\frac{1}{4}$ to $1$, as $\varphi$ increases from $\arccos(\frac{1}{4})$ to $\frac{\pi}{2}$, whereas $\varphi \mapsto \exp(-2\pi \sin(\varphi))$ decreases monotonically from $\exp(-2\pi \sin(\arccos(\frac{1}{4})))$ to 
$\exp(-2\pi)$, yielding:
\begin{eqnarray*}
	\exp(- 2\pi \sin (\varphi)) \, \sin(2\pi \cos(\varphi) 
	+ \varphi)
	> \frac{\exp(-2\pi)}{4} \approx 4,669 \cdot 10^{-4},
	\\  \forall\, \varphi \in
	\Big{[}\arccos\Big{(}\frac{1}{4}\Big{)},\frac{\pi}{2}\Big{]},
\end{eqnarray*}
for the first summand in formula (\ref{Imaginary.Derivative.Eisenstein}).
In combination with estimate (\ref{estim.geom.series}) 
this shows on account of formula (\ref{Imaginary.Derivative.Eisenstein}):
$\frac{\partial \Im G_2}{\partial \varphi}(e^{i\varphi})>0$
for $\varphi \in
\Big{[}\arccos\Big{(}\frac{1}{4}\Big{)},\frac{\pi}{2}\Big{]}$.
Since we know already that $G_2(i)=\pi$, and thus 
$\Im G_2(i)=0$, this especially proves that
$\Im G_2(e^{i\varphi}) < 0$, $\forall \,\varphi \in
\Big{[}\arccos\Big{(}\frac{1}{4}\Big{)}, \frac{\pi}{2}\Big{)}$.
In combination with formula (\ref{Im.Re.Eisenstein}) and with the symmetry of $G_2$, stated in formula (\ref{conjugated.G2}), this implies:
\begin{equation}  \label{The.ineq.1}
	\Re(G_2(\tau)) > \frac{\pi}{\tau_2} > \pi  \quad
	\forall \,\tau = \tau_1 + i \tau_2  \in \partial R_{\Gamma(1)}\cap \sphere^1 \,\,
	\textnormal{satisfying} \,\,\, |\tau_1| \leq \frac{1}{4} \wedge \tau_1\not =0.
\end{equation}
In order to achieve this result also for
$\tau= \tau_1 +i \tau_2 \in \partial R_{\Gamma(1)}\cap \sphere^1$ with $|\tau_1| > \frac{1}{4}$ we invoke 
the real part of formula (\ref{Derivative.Eisenstein.revisited}) which reads:
\begin{equation}  \label{Real.Derivative.Eisenstein}
	\frac{\partial \Re G_2}{\partial \varphi}(e^{i\varphi})
	= 16 \,\pi^3  \, \sum_{n=1}^{\infty} n\,\sigma_1(n) 
	\exp(- 2\pi n\sin(\varphi)) \,
	\cos(2\pi n \cos(\varphi) + \varphi).  \quad
\end{equation}
For $\varphi \in [\frac{\pi}{3}, \arccos(\frac{1}{4})]$
there hold $\cos(2\pi\, \cos(\varphi) + \varphi) \leq -\frac{1}{2}$ and $\exp(-2\pi \sin (\varphi)) \geq
\exp(-2\pi \sin(\arccos(\frac{1}{4})))>0$, and therefore:
\begin{eqnarray}  \label{schmok}
	\exp(- 2\pi \sin (\varphi)) \, \cos(2\pi \cos(\varphi) + \varphi)
	\leq - \frac{\exp\big{(} - 2\pi \sin\big{(}\arccos\big{(}\frac{1}{4}\big{)} \big{)}\big{)}}{2}  \qquad \\
	\approx - 1,1399 \cdot 10^{-3}, \,\,\, \forall \, 
	\varphi \in \Big{[}\frac{\pi}{3}, \arccos\Big{(}\frac{1}{4}\Big{)}\Big{]},     \nonumber
\end{eqnarray}
for the first summand in (\ref{Real.Derivative.Eisenstein}).
Since we can again estimate as in (\ref{estim.geom.series}):
\begin{eqnarray*}
	\sum_{n=2}^{\infty} n\,\sigma_1(n) 
	\exp( - 2\pi n \sin (\varphi)) \,
	|\cos(2\pi n \cos(\varphi) + \varphi)|
	\leq 1,142 \cdot 10^{-4}, \,\, \forall \,\varphi \in \Big{[}\frac{\pi}{3}, \frac{\pi}{2}\Big{]}
\end{eqnarray*}
via the geometric series in (\ref{geometric.series}) with 
$q=q(\varphi):=\exp(-2\pi\, \sin (\varphi))$, we can conclude together with formula (\ref{Real.Derivative.Eisenstein}) and estimate (\ref{schmok}), that
$$
\frac{\partial \Re G_2}{\partial \varphi}(e^{i\varphi}) <0,
\quad \forall \, \varphi \in
\Big{[}\frac{\pi}{3}, \arccos\Big{(}\frac{1}{4}\Big{)}\Big{]},
$$
i.e. that the function $\varphi \mapsto \Re G_2(e^{i\varphi})$
is monotonically decreasing, as $\varphi$ increases from
$\frac{\pi}{3}$ to $\arccos(\frac{1}{4})$.
Since we know already on account of formula (\ref{The.ineq.1})
that $\Re G_2(e^{i\varphi})> \pi$ holds for $\varphi \in
[\arccos(\frac{1}{4}), \frac{\pi}{2})$, we can conclude that
\begin{equation}  \label{The.ineq.2}
	\Re(G_2(\tau)) \geq \Re G_2\Big{(} \frac{1}{4}
	+ i \sin\Big{(}\arccos\Big{(}\frac{1}{4}\Big{)}\Big{)}\Big{)} 
	> \pi,  \,\,\forall \,\tau \in \partial R_{\Gamma(1)}\cap \sphere^1 \,\,\textnormal{with} \, |\tau_1|  \geq \frac{1}{4}.
\end{equation}
Moreover, using the Fourier expansion (\ref{Derivative.Eisenstein})
of the complex derivative of $G_2$ together with the Cauchy-Riemann equations for $G_2$ on $\quat$ we easily obtain:
\begin{eqnarray}  \label{derivative.eisenstein.vertical}
	\frac{d}{dt} \Big{(}\Re G_2\Big{(} \pm \frac{1}{2} 
	+ i t\Big{)} \Big{)}
	=\frac{\partial \Re G_2}{\partial \tau_2}\Big{(} \pm \frac{1}{2} + i t\Big{)}
	= - \frac{\partial \Im G_2}{\partial \tau_1} \Big{(} \pm \frac{1}{2} + i t \Big{)}
	\nonumber \\
	= 16 \,\pi^3 \, \sum_{n=1}^{\infty} n \, \sigma_1(n) (-1)^n \,\exp(-2\pi n \,t) <0, \,\,\, 
	\forall \, t \geq \frac{\sqrt{3}}{2}.      \quad
\end{eqnarray}
The value of the previous sum is negative for every 
$t \geq \frac{\sqrt{3}}{2}$, because one can estimate  
exactly as in estimate (\ref{estim.geom.series}) by 
means of the geometric series in 
(\ref{geometric.series}) - but here 
with $q=q(t):=\exp(-2\pi\,t)$ - that the modulus of 
the sum $\sum_{n=2}^{\infty} n\,\sigma_1(n) \, 
(-1)^n \, \exp(-2\pi n\,t)$
is at least $37,96$ times smaller than the modulus of the first negative summand $-\exp(-2\pi\,t)$ of the entire sum 
in (\ref{derivative.eisenstein.vertical})
for every $t \geq \frac{\sqrt{3}}{2}$. Hence, the function
$[t \mapsto \Re G_2\big{(}\frac{1}{2} + i t\big{)}]$ is strictly monotonically decreasing as $t \nearrow \infty$. 
Furthermore, one immediately infers from the Fourier expansion (\ref{Fourier.expansion.G2}) of $G_2$, combined with the first geometric series in (\ref{geometric.series}) - here again with 
$q=q(t):=\exp(-2\pi\,t)$ - that there holds:
$$
\Re G_2\big{(} \pm \frac{1}{2} + it\big{)}
= \frac{\pi^2}{3} - 8 \pi^2 \, 
\sum_{n=1}^{\infty} \sigma_1(n) \,(-1)^n \exp(-2\pi nt)  \longrightarrow \frac{\pi^2}{3} \quad \textnormal{as} 
\,\,\,t \nearrow \infty,
$$
which finally shows that 
$\Re G_2\Big{(}\pm \frac{1}{2}+it\Big{)}> 
\frac{\pi^2}{3} \,\,\forall \, t \geq \frac{\sqrt{3}}{2}$.
Together with the fact that $G_2(i)=\pi$ and with formulae (\ref{The.ineq.1}) and (\ref{The.ineq.2}) this proves 
the asserted inequality (\ref{The.inequality}) and therefore on account of formula (\ref{cosec.expansion}):
\begin{eqnarray}  \label{Nurd}
	\Re(\eta(1,\ganz \oplus \tau \ganz))=\Re G_2(\tau) \geq \pi,  \quad \forall \, \tau \in \partial R_{\Gamma(1)}.
\end{eqnarray}
Hence, on account of formula (\ref{elastic.energy}) and inequality (\ref{Nurd}) we are led to minimize the function
$g(\omega):=\frac{8 \pi}{\omega} + \frac{2}{3}\, \omega$\, 
for \,$\omega >0$. Differentiation gives the equation
$0 \stackrel{!} = g'(\omega) = 
-\frac{8 \pi}{\omega^2} + \frac{2}{3}$, 
which has the unique positive solution 
$\omega^* =\sqrt{12\, \pi}$. This solution is in fact the global minimizer of $g$ on $\rel_{>0}$,
because of $g''(\omega)= \frac{16 \pi}{\omega^3}> 0$ for any $\omega>0$. Now taking the simple fact into account, 
that for any ``admissible'' pair 
$(m,n) \in \nat \times \nat$  
the number $n$ of periods of the corresponding 
elastic curve $\gamma_{(m,n)}$ must 
be bigger than $1$ according to formula (\ref{correct.frequences}), we can finally combine formula (\ref{elastic.energy}) with estimate (\ref{Nurd}) and infer, that
$$
\Wil(\gamma_{(m,n)}) \geq 2\,\min_{\omega>0} g(\omega)
= 2\,\Big{(}\frac{8\, \pi}{\omega^*} + \frac{2}{3}\, \omega^*\Big{)}
=2\,\Big{(}\frac{8\, \pi}{\sqrt{12\, \pi}} + \frac{2}{3}\, 
\sqrt{12\, \pi}\Big{)} = 16 \, \sqrt{\frac{\pi}{3}}
$$
holds for any non-geodesic, closed elastic curve 
$\gamma_{(m,n)}$ in $\sphere^2$, just as asserted. 
\qed  
\end{remark} 
\noindent  
Another important method of checking both, correctness and 
accuracy of the second statement of our Proposition \ref{Energy.gap} is to compute some numerical values of elastic energies and lengths of elastica $\gamma_{(m,n)}$ and to compare 
the minimum of the values in our first table below 
with our threshold in Proposition \ref{Energy.gap}: $\frac{8\pi}{\sqrt 2}$. 
To this end, we prefer to apply formulae (\ref{length.directly}) and (\ref{energy.directly}) - and not formulae (\ref{length.formula}) and (\ref{elastic.energy}) - 
for coprime integers $1 \leq m \leq 4$ and $1 \leq n \leq 7$, satisfying $\frac{m}{n} \in (0,2-\sqrt 2)$ according to condition (\ref{correct.frequences}).
The interested reader might want to compare these values with the 
corresponding ones for closed, orbitlike elastica in the hyperbolic plane, collected in Table 1 of \cite{Mandel.2018}.
\noindent\\\\
$\underline{\Wil(\gamma_{(m,n)})}$:

\begin{tabular}{*{6}{p{2cm}}}
	
m/n &    m=1     &    m=2    &    m=3    &   m=4      \\
	
n=1 &   ---      &    ---    &    ---    &   ---      \\
	
n=2 &  19,17     &    ---    &    ---    &   ---      \\
	
n=3 &  38,38     &    ---    &    ---    &   ---      \\
	
n=4 &  62,88     &    ---    &    ---    &   ---      \\
	
n=5 &  96,62     &   55,01   &    ---    &   ---      \\
	
n=6 &  134,95    &    ---    &    ---    &   ---      \\
	
n=7 &  192,23    &   98,87   &   74,97   &   62,89    \\\\

\end{tabular}
\noindent
$\underline{\textnormal{length}(\gamma_{(m,n)})}$:  

\begin{tabular}{*{6}{p{2cm}}}
		
m/n &    m=1   &    m=2     &    m=3     &  m=4       \\
	
n=1 &   ---    &    ---     &    ---     &  ---       \\
	
n=2 &  14,68   &    ---     &    ---     &  ---       \\
	
n=3 &  13,68   &    ---     &    ---     &  ---       \\
	
n=4 &  13,98   &    ---     &    ---     &  ---       \\
	
n=5 &  13,77   &    28,51   &    ---     &  ---       \\
	
n=6 &  13,99   &    ---     &    ---     &  ---       \\
	
n=7 &  13,15   &    27,95   &   41,63    &  60,22     \\	

\end{tabular}\\\\\\
\noindent
\underline{Acknowledgements:}\\\\
The author would like to thank Professor Itai Shafrir 
and Professor Yehuda Pinchover for their strong support 
and hospitality at the Mathematics Department of the 
``Israel Institute of Technology''. The author is also indebted to the kind referee on account of his smart and helpful suggestions. The author was funded by the Ministry of Absorption 
of the State of Israel in the academic years 
2019/2020, 2020/2021 and 2021/2022.


\end{document}